\documentclass[11pt]{amsart}
\usepackage{amsmath,amsfonts,amscd,amssymb,epsf, euscript, eufrak}
\usepackage{amsxtra,mathrsfs} 
\newtheorem{theorem}{Theorem}[section]
\newtheorem{proposition}[theorem]{Proposition}
\newtheorem{lemma}[theorem]{Lemma}
\newtheorem{corollary}[theorem]{Corollary}
\theoremstyle{definition}
\newtheorem{definition}[theorem]{Definition}

\theoremstyle{remark} \newtheorem{remark}[theorem]{Remark}
\numberwithin{equation}{section}
\newcommand{\field}[1]{\ensuremath{\mathbb{#1}}}
\newcommand{\CC}{\field{C}}

\newcommand{\HH}{\field{H}}
\newcommand{\PP}{\field{P}}
\newcommand{\RR}{\field{R}}

\newcommand{\ZZ}{\field{Z}}

\newcommand{\complex}[1]{\mathsf{#1}} 

\newcommand{\SSS}{\complex{S}}

\newcommand{\KKK}{\complex{K}}

\newcommand{\FFF}{\complex{F}}

\newcommand{\Ka}{K\"{a}hler\:}
\newcommand{\Te}{Teichm\"{u}ller\:}
 \DeclareMathOperator{\Ker}{Ker}
\DeclareMathOperator{\Tr}{Tr} \DeclareMathOperator{\Diff}{Diff}
\DeclareMathOperator{\Mob}{M\ddot{o}b}

 \DeclareMathOperator{\PSL}{PSL}
\DeclareMathOperator{\PSU}{PSU} \DeclareMathOperator{\re}{Re}
\DeclareMathOperator{\Sp}{Sp}
\DeclareMathOperator{\Det}{Det}

\newcommand{\Del}{\mathbb{D}}
\newcommand{\del}{\partial}

\newcommand{\Z}{{\mathbb{Z}}}

\newcommand{\N}{{\mathbb{N}}}
\newcommand{\R}{{\mathbb{R}}}
\newcommand{\C}{{\mathbb{C}}}

\newcommand{\curly}[1]{\mathscr{#1}}

\newcommand{\cB}{\curly{B}}

\newcommand{\cE}{\curly{E}}

\newcommand{\cP}{\curly{P}}

\newcommand{\cS}{\curly{S}}

\newcommand{\cV}{\curly{V}}

\newcommand{\bk}{\backslash}
\newcommand{\pa}{\partial}
\newcommand{\la}{\langle}
\newcommand{\ra}{\rangle}

\newcommand{\ov}{\overline}

\newcommand{\vep}{\varepsilon}
\newcommand{\f}{\mathrm{f}}
\newcommand{\g}{\mathrm{g}}
\newcommand{\z}{\bar{z}}

\begin{document}

\title[Potential of the Weil-Petersson metric on $T(1)$]
{Weil-Petersson metric on the universal
Teichm\"{u}ller space II. K\"{a}hler potential and period mapping} 
\author{Leon A. Takhtajan} \address{Department of Mathematics \\
Stony Brook University\\ Stony Brook, NY 11794-3651 \\ USA}
\email{leontak@math.sunysb.edu}
\author{Lee-Peng Teo} \address{Department of Applied Mathematics \\
National Chiao Tung University, 1001 \\ Ta-Hsueh Road, Hsinchu
City, 30050 \\ Taiwan, R.O.C.} \email{lpteo@math.nctu.edu.tw}

\begin{abstract} 
We study the Hilbert manifold structure on $T_{0}(1)$ --- the connected component 
of the identity of the Hilbert manifold $T(1)$. We characterize points on $T_{0}(1)$
in terms of Bers and pre-Bers embeddings, and prove that the Grunsky operators 
$B_{1}$ and $B_{4}$, associated with the points in $T_{0}(1)$ via conformal welding,
are Hilbert-Schmidt. We define a ``universal Liouville action'' --- a real-valued 
function $\SSS_{1}$ on $T_{0}(1)$, and prove that it is a K\"{a}hler potential of 
the Weil-Petersson metric on $T_{0}(1)$. We also prove that $\SSS_{1}$ is 
$-\tfrac{1}{12\pi}$ times the logarithm of the Fredholm determinant of associated 
quasi-circle, which generalizes classical results of Schiffer and Hawley. We define 
the universal period mapping $\hat{\cP}: T(1)\rightarrow\cB(\ell^{2})$ of $T(1)$ 
into the Banach space of bounded operators on the Hilbert space $\ell^{2}$, prove 
that $\hat{\cP}$ is a holomorphic mapping of Banach manifolds, and show that 
$\hat{\cP}$ coincides with the period mapping introduced by Kurillov and Yuriev and 
Nag and Sullivan. We prove that the restriction of $\hat{\cP}$ to $T_{0}(1)$ is an 
inclusion of $T_{0}(1)$ into the Segal-Wilson universal Grassmannian, which is a 
holomorphic mapping of Hilbert manifolds. We also prove that the image of the 
topological group $S$ of symmetric homeomorphisms of $S^{1}$ under the mapping 
$\hat{\cP}$ consists of compact operators on $\ell^{2}$.
\end{abstract}
\maketitle
\tableofcontents
\section{Introduction} 
This is the second part of our paper \cite{TT-I}, which will be referred to as Part I. Here we continue our investigation of the Weil-Petersson metric on the universal \Te space $T(1)$. Namely, we study in detail the Hilbert manifold structure of $T(1)$ and establish relations between the Hilbert submanifold $T_{0}(1)$ --- the connected component of the identity in $T(1)$, and classical Grunsky operators $B_{l}$, $l=1,2,3,4$, associated with the conformal welding. In Part I, we have described the image of $T_{0}(1)$ under the Bers embedding $\beta: T(1)\rightarrow A_{\infty}(\Del)$. Here we characterize $T_{0}(1)$ in terms of the pre-Bers embedding
$\hat{\beta}: T(1)\rightarrow A_{\infty}^{1}(\Del)$ and prove that
the Grunsky operators $B_{1}$ and $B_{4}$ associated with the points in $T_{0}(1)$ are Hilbert-Schmidt. We establish the relation between eigenvalues of Grunsky operators and classical Fredholm eigenvalues, generalizing Schiffer's result for $C^{3}$ curves \cite{Schiffer2}. We prove that  the logarithm of the Fredholm determinant of the operator $I-B_{1}B_{1}^{\ast}$ associated with points in $T_{0}(1)$ (or, which is the same, of the Fredholm determinant of $I-B_{4}B_{4}^{\ast}$) is, up to a constant, a \Ka potential for the Weil-Petersson metric on $T_{0}(1)$. We prove the explicit formula for this Fredholm determinant, expressing it as the ``universal Liouville action''.  Using Grunsky operators, we define the universal period mapping $\cP$ of $T_{0}(1)$ into the Hilbert space $\cS_{2}$ of Hilbert-Scmidt operators on the Hilbert space $\ell^{2}$, as well as the mapping $\hat{\cP}$ of $T(1)$ into the Banach space $\cB(\ell^{2})$ of bounded operators on $\ell^{2}$. We prove that $\cP$ and $\hat{\cP}$ are holomorphic mappings of Hilbert and Banach manifolds respectively. We show that the mapping $\hat{\cP}$ coincides with the period mapping, first introduced by Kirillov and Yuriev \cite{KY2} for the homogenous space $\Mob(S^1)\bk\Diff_+(S^1)$, studied in detail by Nag \cite{Nag92}, and then extended to $T(1)$ by Nag and Sullivan \cite{NS}\footnote{It is explained in \cite{Nag92} and \cite{NS} it what sense the mapping $\hat{\cP}$ generalizes the classical period mapping of compact Riemann surfaces.}. Finally, we prove that the image of the topological group $S$ of symmetric homeomorphisms of $S^{1}$ under the period mapping $\hat{\cP}$ is $\cS_{\infty}\cap\hat{\cP}(T(1))$, where $\cS_{\infty}$ is the ideal of the Banach algebra $\cB(\ell^{2})$ consisting of compact operators on $\ell^{2}$.  

Below is the detailed description of the paper. In what follows 
we are using notations and results from Part I; in particular, the normalization of the conformal welding for $T(1)$, described in Section 2.2.1, Part I. Namely, for every $[\mu]\in T(1)$ we consider the q.c.~mapping $w_{\mu}$ that fixes $-1, -i, 1$ as an element of
$S^{1}\bk\mathrm{Homeo}_{qs}(S^{1})$, which admits a conformal welding
\begin{align*}
w_{\mu} =\g_{\mu}^{-1}\circ\f^\mu,
\end{align*}
where $\f^{\mu}$ and $\g_{\mu}$ are q.c.~mappings whose restrictions on $\Del$ and $\Del^{\ast}$, respectively, are holomorphic functions satisfying
$\f^{\mu}(0)=0, (\f^{\mu})'(0)=1$ and $\g_{\mu}(\infty)=\infty$.

In Section 2 we characterize the univalent functions associated with the Hilbert manifold $T_{0}(1)$ in terms of the Hilbert spaces $A_{2}^{1}(\Del)$ and $A_{2}^{1}(\Del^{\ast})$ of holomorphic functions on $\Del$ and $\Del^{*}$ respectively, square integrable with respect to the Lebesgue measure. Using the embedding $A_{2}^{1}(\Del)\hookrightarrow A_{\infty}^{1}(\Del)$
into the Banach space of holomorphic functions on $\Del$, the Becker-Pommerenke theorem \cite{BecPom}, and the characterization of the topological group $S$ of symmetric homeomorphisms of $S^{1}$ given by Gardiner and Sullivan \cite{GaSu}, we prove that $T_{0}(1)$ is a subgroup of $S$.
The main result of this section is Theorem 2.12, which states that $[\mu]\in T_{0}(1)$ if and only if one of the following conditions holds: (i) $\mathcal{S}(\f^{\mu})\in A_{2}(\Del)$; (ii) $\mathcal{A}(\f^{\mu})\in A_{2}^{1}(\Del)$; (iii)  $\mathcal{S}(\g_{\mu})\in A_{2}(\Del^{\ast})$; (iv)  $\mathcal{A}(\g_{\mu})\in A_{2}^{1}(\Del^{\ast})$. 
Here $\mathcal{S}(f)$ is the Schwarzian derivative of the univalent function $f$, and 
\begin{align*}
\mathcal{A}(f)=\frac{f''}{f'}.
\end{align*}
This theorem allows us to introduce the ``universal Liouville action'' --- the function $\SSS_{1}: T_{0}(1)\rightarrow \RR$, defined by
\begin{align}
\SSS_1([\mu]) =
\iint\limits_{\Del}\left|\mathcal{A}(\f^{\mu})\right|^2d^2z+
\iint\limits_{\Del^*} \left|\mathcal{A}(\g_{\mu})\right|^2 d^2z
-4\pi\log|\g_{\mu}'(\infty)|.
\end{align}
In Section 3 to every $[\mu] \in T(1)$ we assign the Grunksy operators $B_{1},B_{2},B_{3}$ and $B_{4}$, associated with the corresponding pair $(\f^{\mu},\g_{\mu})$ of univalent functions. The Lebesgue measure of the quasi-circle $\CC\setminus\{f(\Del)\cup g(\Del^{\ast})\}$ is zero, so that the generalized Grunsky inequality \cite{Hummel,Pom} can be succinctly formulated as the unitarity of the operator
$\mathbf{B}= \begin{pmatrix}
     B_{1} & B_{2}\\
     B_{3} & B_{4}
     \end{pmatrix}$
on $\ell^{2}\oplus\ell^{2}$. The main result of Section 3.1 is Theorem 3.6, which states (see Corollary 3.9) that $[\mu]\in T_{0}(1)$ if and only if the corresponding Grunsky operators $B_{1}(\f^{\mu}), B_{4}(\g_{\mu})\in\cS_{2}$ --- the Hilbert space of Hilbert-Schmidt operators on $\ell^{2}$. In Theorem 3.10 we prove that the mapping $\cP: T_{0}(1)\rightarrow \cS_{2}$, defined by $ \cP([\mu])=B_{1}(\f^{\mu})$, is a holomorphic mapping of Hilbert manifolds. Extended to the universal \Te space $T(1)$, this defines a holomorphic mapping $\hat{\cP}: T(1)\rightarrow \cB(\ell^{2})$ of Banach manifolds, which we prove in Appendix B. In Section 3.2 we show that for $[\mu]\in T_{0}(1)$ the eigenvalues of the corresponding trace class operators
$B_{1}B_{1}^{\ast}$ and $B_{4}B_{4}^{\ast}$ are related to the eigenvalues of the classical Poincar\'{e}-Fredholm integral operator associated with the quasi-circle $\mathcal{C}=\f^{\mu}(S^{1})=\g_{\mu}(S^{1})$. Since for $[\mu]\in T_{0}(1)$ these quasi-circles contain all $C^{3}$ curves, this generalizes Schiffer's result \cite{Schiffer2}. Extending \cite{Schiffer-multiple}, we introduce the Fredholm determinant $\Det_{F}(\mathcal{C})$ of the  quasi-circle $\mathcal{C}$ as the Fredholm determinant $\det(I-B_{1}B_{1}^{\ast})=\det(I-B_{4}B_{4}^{\ast})$, and define the function $\SSS_{2}: T_{0}(1)\rightarrow\RR$ by
\begin{equation}
\SSS_{2}([\mu])=\log\Det_{F}(\f^{\mu}(S^{1})),\;\;[\mu]\in T(1).
\end{equation}

In Section 3.3 we define the semi-infinite period matrices of $1$-forms for natural bases of $A_{2}^{1}(\Del)$ and $A_{2}^{1}(\Del^{*})$, which generalize imaginary parts of the classical period matrices for
compact Riemann surfaces, and show that they correspond to the operators $B_{2}B_{2}^{\ast}$ and $B_{3}B_{3}^{\ast}$. 

In Section 4 we compute the ``first variations'' of the functions $\SSS_{1}$ and $\SSS_{2}$ --- the $(1,0)$-forms $\del\SSS_{1}$ and
$\del\SSS_{2}$, where $\del$ is the $(1,0)$-component of the de Rham differential on the Hilbert manifold $T_{0}(1)$. Namely, we show in Theorems 4.5 and 4.1  (see Corollaries 4.9 and 4.2) that
\begin{align}
\del \SSS_{1}=2\boldsymbol{\vartheta}\quad\text{and}\quad
\del\SSS_{2}=-\frac{1}{6\pi}\boldsymbol{\vartheta},
\end{align}
where the $(1,0)$-form $\boldsymbol{\vartheta}$ on $T_{0}(1)$, under the natural isomorphism $T_{[\mu]}^{\ast}T_{0}(1)\simeq A_{2}(\Del^{\ast})$, is given by
\begin{equation}
\boldsymbol{\vartheta}_{[\mu]}=\mathcal{S}(\g_{\mu}).
\end{equation}
The proof of Theorem 4.1 is rather standard, whereas the proof of Theorem 4.5 relies heavily on the identity given in Lemma 4.6. The latter can be interpreted as an extension of the generalized Grunsky equality to pairs of univalent functions $(\f^{\mu},\g_{\mu})$ for
$[\mu]\in T_{0}(1)$, which we consider quite interesting. Since the functions $\SSS_{1}$ and $\SSS_{2}$ on $T_{0}(1)$ both vanish at $0\in T_{0}(1)$, from (1.3) we immediately obtain that
$$\SSS_{2} = -\frac{1}{12\pi}\SSS_{1},$$
thus expressing the Fredholm determinant as the universal Liouville action. In Corollary 4.12 and Remark 4.13 we interpret this relation as a surgery type formula for the determinants of elliptic operators on domains on the Riemann sphere $\mathbb{P}^{1}$.

In Section 5 we show that the relation (1.3) implies that the function $\SSS_{1}$ is a K\"{a}hler potential of the Weil-Petersson metric on $T_{0}(1)$. The proof goes along the same lines as in the case of finite-dimensional \Te spaces \cite{LT}. This explains why the function $\SSS_{1}$ is called the universal Liouville action. In Section 6 we study the period mapping $\hat{\cP}: T(1)\rightarrow\cB(\ell^{2})$. We prove that it coincides with the Kirillov-Yuriev-Nag-Sullivan mapping of $T(1)$  into the infinite-dimensional analog of Siegel disk $\mathfrak{D}_{\infty}$. We also show that the period mapping $\cP: T_{0}(1)\rightarrow\cS_{2}$ gives an embedding of $T_{0}(1)$ into the Segal-Wilson universal Grassmannian. 

In Appendix A we study the Hilbert manifold structure on the topological group $\mathcal{T}_{0}(1)$ --- the pre-image of the Hilbert manifold $T_{0}(1)$ under the canonical projection $\pi:\mathcal{T}(1)\rightarrow T(1)$. We prove in Theorem A.3 that the Bers embedding $\beta:\mathcal{T}_{0}(1)\rightarrow A_{2}(\Del)\oplus\CC$ and the pre-Bers embedding 
$\hat{\beta}:\mathcal{T}_{0}(1)\rightarrow A_{2}^{1}(\Del)$ induce the same Hilbert manifold structure on $\mathcal{T}_{0}(1)$. This result is parallel to the one proved in the Appendix of \cite{Teo}. We also prove Corollaries A.4 and A.6, characterizing convergence in the Hilbert manifold topology of $\mathcal{T}_{0}(1)$, which were used in the proof of Lemma 4.6. 
Finally, in Appendix B we show that $\hat{\cP}: T(1)\rightarrow \cB(\ell^{2})$ is a holomorphic mapping of Banach manifolds and prove that the image of the topological group $S$ under the map $\hat{\cP}$ is the subset $\cS_{\infty}\cap\hat{\cP}(T(1))$ of $\cB(\ell^{2})$. The properties of the tower of embedded manifolds $T_{0}(1)\hookrightarrow S \hookrightarrow T(1)$ are summarized in a commutative diagram at the end of Appendix B.
\vspace{3mm}

\noindent \textbf{Acknowledgments.} We appreciate useful
discussions with M.~Luybich. The second author would like to thank P.Y.~Wu for helpful discussions about operator theory. The work of the first author was
partially supported by the NSF grant DMS-0204628. The work of the
second author was partially supported by the grant NSC
92-2115-M-009-017. The second author also thanks CTS for the
fellowship to visit Stony Brook University in the Summer of 2003,
where a part of this work was done.
\section{Hilbert spaces of univalent functions} 
\label{Hilbertunivalent}
It is well-known (see, e.g., Section 2.2 in Part I) that the universal
Teichm\"uller space $T(1)$ is isomorphic to the space $\mathcal{D}$
of univalent functions on $\Del$. Here we characterize the univalent
functions associated with the Hilbert manifold $T_0(1)$.

In addition to the Hilbert spaces $A_2(\Del)$ and $A_2(\Del^*)$, introduced in Section 3, Part I, we define the following Hilbert spaces of holomorphic functions,
\begin{align*}
A^1_{2}(\Del)&= \left\{\psi \;\text{holomorphic on}\,\, \Del: \,
\|\psi\|_{2}^2 = \iint\limits_{ \Del}
\left|\psi(z)\right|^2d^2z < \infty \right\},\\
A^1_{2}(\Del^*)&= \left\{\psi \;\text{holomorphic on}\,\, \Del^*:
\, \|\psi\|_{2}^2 = \iint\limits_{ \Del^*}
\left|\psi(z)\right|^2d^2z < \infty \right\}.
\end{align*}
We denote by $\ov{A^{1}_{2}(\Del)}$ and $\ov{A^{1}_{2}(\Del^{*})}$ the corresponding Hilbert spaces of antiholomorphic functions. 
\begin{remark}
Every $\psi\in A_2^1(\Del)$ corresponds to a holomorphic $1$-form $\omega=\psi(z) dz$ on $\Del$ (or on $\Gamma\bk\Del$ for a cofinite Fuchsian group $\Gamma$) such that the $(1,1)$-form $\omega\wedge \bar{\omega}$ is integrable. Similarly, every $\phi\in A_2(\Del)$ corresponds to a holomorphic quadratic
differential $q=\phi(z) (dz)^2$ on $\Del$ (or on $\Gamma\bk\Del)$ such that the $(1,1)$-form $ (|\phi(z)|^2/\rho(z))dz\wedge d\z$ is integrable, so that the latter space could be also denoted by $A_2^2(\Del)$. We will use the same notation $\Vert~\Vert_2$ for the norms in these Hilbert spaces. To avoid confusion, in the main text we always denote elements in the spaces $A_2$ by $\phi$,
and elements in the spaces $A_2^1$ by $\psi$. 
\end{remark}
In addition to the Banach spaces $A_{\infty}(\Del)$ and $A_{\infty}(\Del^{\ast})$ introduced in Section 2.1, Part I, we define the following Banach spaces of holomorphic functions,
\begin{align*}
A^1_{\infty}(\Del)&= \left\{\psi \;\text{holomorphic on}\,\, \Del:
\, \|\psi\|_{\infty} = \sup_{z\in \Del}
\left|(1-|z|^2)\psi(z)\right| < \infty\right\},\\
A^1_{\infty}(\Del^*)&= \left\{\psi \;\text{holomorphic on}\,\,
\Del^*: \, \|\psi\|_{\infty} = \sup_{z\in \Del^*}
\left|(1-|z|^2)\psi(z)\right| < \infty \right\}.
\end{align*}
For a holomorphic function $f:\Omega\rightarrow\CC$ such that
$f^\prime\neq 0$ on $\Omega$ we set
\begin{align*}
\mathcal{A}(f)= \frac{f^{\prime\prime}}{f^{\prime}}.
\end{align*}
\begin{remark}
Classical distortion theorem (see e.g., \cite{Pom, Duren}) implies
that if $f:\Del\rightarrow \C$ and $g:\Del^*\rightarrow \C$ are
univalent functions, then $\mathcal{A}(f)\in A^1_{\infty}(\Del)$ and $\mathcal{A}(g)\in A^1_{\infty}(\Del^*)$. In \cite{Teo}, it was shown that
the Bers embedding of the universal Teichm\"uller curve $\mathcal{T}(1)$ into $A_{\infty}(\Del)\oplus \C$ can be factorized as the composition of two holomorphic embeddings 
$$\mathcal{T}(1)\rightarrow
A^{1}_{\infty}(\Del)\rightarrow A_{\infty}(\Del)\oplus \C.$$ 
Here the map $\mathcal{T}(1)\rightarrow A^{1}_{\infty}(\Del)$ is given by
$\gamma=g^{-1}\circ f \mapsto \mathcal{A}(f)$
and the map $A^{1}_{\infty}(\Del)\rightarrow A_{\infty}(\Del)\oplus
\C$ is defined as
\begin{align*}
\psi\mapsto \left( \psi_z-\tfrac{1}{2}\psi^2,
\tfrac{1}{2}\psi(0)\right).
\end{align*}
\end{remark}
Similar to Lemma 3.1 in Part I, we have
\begin{lemma}\label{subspace1}
The vector spaces $A^1_{2}(\Del)$ and $A^1_{2}(\Del^*)$ are
subspaces of $A^1_{\infty}(\Del)$ and $A^1_{\infty}(\Del^{*})$ 
respectively. The natural inclusion maps 
$A^{1}_{2}(\Del) \hookrightarrow A^{1}_{\infty}(\Del)$ and 
$A^{1}_{2}(\Del^*) \hookrightarrow A^{1}_{\infty}(\Del^*)$
are bounded linear mappings of Banach spaces.
\end{lemma}
\begin{proof}
It is sufficient to consider only the spaces of holomorphic functions on $\Del$. For every $\psi \in
A^{1}_{2}(\Del)$ let $\psi(z)= \sum_{n=1}^{\infty}na_n
z^{n-1}$ be the power series expansion. Then
\begin{align*}
\Vert\psi\Vert_2^2=\iint\limits_{\Del} |\psi(z)|^2d^2z = \pi\sum_{n=1}^{\infty} n
|a_n|^2,
\end{align*}
and by Cauchy-Schwarz inequality, we have
\begin{align*}
|\psi(z)| \leq \sum_{n=1}^{\infty} n|a_n||z|^{n-1} \leq
\left(\sum_{n=1}^{\infty} n|a_n|^2\right)^{1/2}
\left(\sum_{n=1}^{\infty} n|z|^{2n-2}\right)^{1/2}
\end{align*}
for every $z\in\Del$.
Using
\[
\sum_{n=1}^{\infty} n|z|^{2n-2} = \frac{1}{(1-|z|^2)^2},
 \]
we get
\begin{align*}
\Vert \psi \Vert_{\infty}=\sup\limits_{z\in\Del}|(1-|z|^2)\psi(z)| \leq\frac{1}{\sqrt{\pi}}\Vert
\psi\Vert_{2}.
\end{align*}
\end{proof}
Similar to Remark 3.2 in Part I, we get
\begin{corollary}
For every $\psi\in A^1_{2}(\Del)$, 
\begin{align*}
\lim_{|z|\rightarrow 1^-} (1-|z|^2)\psi(z)=0.
\end{align*}
Similar statement holds for every $\psi\in A^1_2(\Del^*)$.
\end{corollary}
For a holomorphic function $f:\Omega\rightarrow\CC$ set 
\begin{equation*}
\Psi(f)=f_z -\tfrac{1}{2}f^2.
\end{equation*}
If $f^{\prime}\neq 0$ on $\Omega$, then $\mathcal{S}(f)=(\Psi\circ\mathcal{A})(f)$,
where $\mathcal{S}(f)$ is the Schwarzian derivative of $f$. 
In \cite{Teo} it was proved that 
$\Psi\left(A^1_{\infty}(\Del)\right)\subset A_\infty(\Del)$ and 
$\Psi\left(A^1_{\infty}(\Del^*)\right)\subset A_\infty(\Del^*)$.
Similarly, we have the following result.
\begin{lemma}\label{sub1} 
$\Psi\left(A^1_2(\Del)\right)\subset A_2(\Del)$ and 
$\Psi\left(A^1_2(\Del^*)\right)\subset A_2(\Del^*)$.
\end{lemma}
\begin{proof}
Again it is sufficient to consider functions on $\Del$. For
$\psi=\sum_{n=1}^{\infty} na_nz^{n-1} \in A^1_{2}(\Del)$ we have
\begin{align*}
\iint\limits_{\Del} |\Psi(\psi)|^2 \rho(z)^{-1} d^2z \leq
2\iint\limits_{\Del}  |\psi_z(z)|^2 \rho(z)^{-1} d^2z +
\frac{1}{2} \iint\limits_{\Del}
|\psi(z)|^4\rho(z)^{-1} d^2z.
\end{align*}
For the first term, a straightforward computation gives
\begin{align*}
\iint\limits_{\Del} |\psi_z(z)|^2 \rho(z)^{-1} d^2z=
\frac{\pi}{2}\sum_{n=2}^{\infty} \frac{n(n-1)}{n+1} |a_n|^2<
\tfrac{1}{2}\Vert \psi \Vert_{2}^2<\infty.
\end{align*}
For the second term, since $\psi\in A^1_{\infty}(\Del)$, we have
\begin{align*}
\iint\limits_{\Del} \rho(z)^{-1} |\psi(z)|^4d^2z \leq\tfrac{1}{4}\Vert\psi\Vert_{\infty}^{2}
\Vert \psi\Vert_{2}^2<\infty.
\end{align*}
\end{proof}

The following theorem of Becker and Pommerenke \cite{BecPom} characterizes univalent functions on $\Del$ that admit a q.c.~extension to a larger domain such that the complex dilation is continuous on $S^{1}$.
\begin{theorem}\label{BecPomThm}
Let $f: \Del\rightarrow \C$ be a univalent function such that
$f(\Del)$ is a Jordan domain. Then the following conditions are
equivalent.
\begin{itemize}
\item[(i)] $f$ has a q.c.~extension $F$ to
$\{z : |z|<R, R>1\}$ such that the complex dilation $\mu(z) =
F_{\z}/F_z$ satisfies
\begin{align*}
\lim_{|z|\rightarrow 1^+} \mu(z)=0.
\end{align*}
\item[(ii)]
\begin{align*}
\lim_{|z|\rightarrow 1^-} (1-|z|^2)^2 \mathcal{S}(f)(z) =0.
\end{align*}
\item[(iii)]
\begin{align*}
\lim_{|z|\rightarrow 1^-} (1-|z|^2) \mathcal{A}(f)(z) =0.
\end{align*}
\end{itemize}
\end{theorem}
In \cite{GaSu}, Gardiner and Sullivan have studied the subgroup $$S=\Mob(S^{1})\bk\text{Homeo}_{s}(S^{1})$$
of symmetric homeomorphisms in $QS=\Mob(S^{1})\bk\text{Homeo}_{qs}(S^{1})\simeq T(1)$. They proved that
as a Banach submanifold of $T(1)$,  $S$ is a topological group, 
and that univalent functions $f$ associated to elements in $S$ are precisely the functions satisfying condition (ii)
of Theorem \ref{BecPomThm}. 
\begin{remark} \label{ga-su} Actually in \cite{GaSu} this condition is stated as follows: for every $\vep>0$ there is a compact subset $K$ of $\Del$ such that $|(1-|z|^2)^2 \mathcal{S}(f)(z)|<\vep$ for $z\in\Del\setminus K$, which is clearly equivalent to (ii).
\end{remark}
Using this remark and Remark 3.2 in Part I, we get the following statement. 
\begin{corollary}
The group $T_0(1)$ is a subgroup of $S$.
\end{corollary}
\begin{remark} It is known \cite{GaSu} that the topological group $S$ contains the subgroup of $C^{1}$-homeomorphisms of $S^{1}$. Similarly, the topological group $T_{0}(1)$ contains the subgroup of 
$C^{3}$-homeomorphisms. Indeed, it is known (see, e.g., \cite{Hamilton}) that
if $\gamma\in QS$ is $C^{3}$ then corresponding $\f$ and $\g$ are of $C^{2}$ class on the boundary and all their derivatives are Holder continuous with $\alpha<1$. From here it follows that $\mathcal{S}(f)\in A_{2}(\Del)$.
\end{remark}
\begin{remark}
For $[\mu]\in T_0(1)$ it is an interesting open problem to characterize
intrinsically the corresponding map $\left.w_{\mu}\right|_{S^1}$ and the
quasi-circle $\f(S^1)$, as it was done for by Gardiner and Sullivan in \cite{GaSu} for $[\mu]\in S$. 
\end{remark}
Another important consequence of Becker-Pommerenke Theorem is the following result.
\begin{lemma}\label{sub2} Let $f$ and $g$ be univalent functions on
$\Del$ and $\Del^*$ such that $\mathcal{S}(f)\in A_2(\Del)$ and
$\mathcal{S}(g)\in A_2(\Del^*)$. Then $\mathcal{A}(f)\in A^1_2(\Del)$
and $\mathcal{A}(g)\in A^1_2(\Del^*)$.
\end{lemma}
\begin{proof} It is sufficient to consider functions on $\Del$.
If $\mathcal{S}(f)\in A_2(\Del)$, then by Remark 3.2 in Part I $f$ satisfies the condition (ii) in Theorem \ref{BecPomThm} and hence it satisfies the condition (iii). In particular, there exists $r'>0$ such that 
$$(1-|z|^2)|\mathcal{A}(f)(z)|\leq 1/2\quad\text{for all}\quad r'<|z|<1.$$ 
By triangle and geometric mean inequalities,  
\begin{align*}
\left|\mathcal{S}(f)(z)\right|^{2}\geq &\left(|\mathcal{A}(f)^{\prime}(z)| -\tfrac{1}{2}|\mathcal{A}(f)(z)|^{2}\right)^{2}\\
=&|\mathcal{A}(f)^{\prime}(z)|^{2} + \tfrac{1}{4}|\mathcal{A}(f)(z)|^{4} -
|\mathcal{A}(f)^{\prime}(z)||\mathcal{A}(f)(z)|^{2}\\
\geq &\tfrac{1}{2}\left(|\mathcal{A}(f)^{\prime}(z)|^{2} -|\mathcal{A}(f)(z)|^{4}\right),
\end{align*}
so that for $r'<|z|<1$,
\begin{equation} \label{ineq1}
2(1-|z|^{2})^{2}|\mathcal{S}(f)(z)|^{2}\geq 
(1-|z|^{2})^{2}|\mathcal{A}(f)^{\prime}(z)|^{2} -\tfrac{1}{4}|\mathcal{A}(f)(z)|^{2}.
\end{equation}
Let $\mathcal{A}(f)(z)=\sum_{n=1}^{\infty}na_nz^{n-1}$ be the
power series expansion of $\mathcal{A}(f)$ and let $\Del_{r}$ be the disk of radius $r$. We have,
\begin{align*}
\iint\limits_{\Del_{r}}(1-|z|^2)^2|\mathcal{A}(f)^{\prime}(z)|^2d^2z
=&\pi \sum_{n=2}^{\infty}n^2(n-1)^2 |a_n|^2 r^{2n}
\left(\frac{r^{-2}}{n-1}-\frac{2}{n}+\frac{r^{2}}{n+1}\right)
\end{align*}
and
\begin{align*}
\iint\limits_{\Del_{r}}|\mathcal{A}(f)(z)|^2d^2z=\pi\sum_{n=1}^{\infty} n|a_n|^2 r^{2n}.
\end{align*}
Using the elementary inequality
$$ (n-1)^2\left(\frac{r^{-2}}{n-1}-\frac{2}{n}
+\frac{r^{2}}{n+1}\right)\geq \frac{1}{2n}$$
for all $n\geq 2$ and $0<r<1$, we get
\begin{align}\label{inequality2}
\iint\limits_{\Del_{r}}(1-|z|^2)^2|\mathcal{A}(f)'(z)|^2d^2z \geq \tfrac{1}{2} \iint\limits_{\Del_{r}}|\mathcal{A}(f)(z)|^2d^2z
 -\tfrac{\pi}{2} |a_{1}|^{2}r^{2}.
\end{align}
Integrating the inequality \eqref{ineq1} over $\Del_{r}\setminus\Del_{r'}$, and using \eqref{inequality2}, we get for
$r>r^{\prime}$,
\begin{gather*}
2\iint\limits_{\Del_{r}\setminus\Del_{r'}}(1-|z|^{2})^{2}|\mathcal{S}(f)(z)|^{2}d^{2}z\geq \iint\limits_{\Del_{r}\setminus\Del_{r'}}\left((1-|z|^{2})^{2}|\mathcal{A}(f)^{\prime}(z)|^{2} - \tfrac{1}{4}|\mathcal{A}(f)(z)|^{2}\right)
d^{2}z \\
=\iint\limits_{\Del_{r}}(1-|z|^{2})^{2}|\mathcal{A}(f)^{\prime}(z)|^{2}
d^{2}z - \tfrac{1}{4}\iint\limits_{\Del_{r}}|\mathcal{A}(f)(z)|^{2}d^{2}z \\
-\iint\limits_{\Del_{r'}}(1-|z|^{2})^{2}|\mathcal{A}(f)^{\prime}(z)|^{2}
d^{2}z + \tfrac{1}{4}\iint\limits_{\Del_{r'}}|\mathcal{A}(f)(z)|^{2}d^{2}z \\
\geq \tfrac{1}{4}\iint\limits_{\Del_{r}}|\mathcal{A}(f)(z)|^{2}d^{2}z +\tfrac{1}{4}\iint\limits_{\Del_{r'}}|\mathcal{A}(f)(z)|^{2}d^{2}z \\
-\iint\limits_{\Del_{r'}}(1- |z|^{2})^{2}|\mathcal{A}(f)^{\prime}(z)|^{2}d^{2}z -\tfrac{\pi}{2}|a_{1}|^{2} r^{2}.
\end{gather*}
Since $\mathcal{S}(f)\in A_{2}(\Del)$, from this inequality we conclude that there exists $C>0$ such that
$$\iint\limits_{\Del_{r}}|\mathcal{A}(f)(z)|^{2}d^{2}z <C$$
for all $0<r<1$, i.e., $\mathcal{A}(f)\in A^{1}_{2}(\Del)$.
\end{proof}
The following statement is the main result of this section.
\begin{theorem}\label{hilbertspaces}
Let $w_{\mu} = \g_{\mu}^{-1}\circ \f^{\mu}$ be the conformal 
welding corresponding to $[\mu]\in T(1)$. Then $[\mu]\in T_{0}(1)$ if and only if one of the following conditions holds.
\begin{itemize}
\item[(i)] $\mathcal{S}(\f^\mu)\in A_2(\Del)$.
\item[(ii)] $\mathcal{A}(\f^\mu) \in A_2^1(\Del)$.
\item[(iii)] $\mathcal{S}(\g_\mu)\in A_2(\Del^*)$.
\item[(iv)] $\mathcal{A}(\g_\mu)\in A_2^1(\Del^*)$.
\end{itemize} 
\end{theorem}
\begin{proof} 
Since under the Bers embedding $\beta(T_{0}(1))=\beta(T(1))\cap A_{2}(\Del)$, it follows that  
if $w_{\mu}= \g_{\mu}^{-1}\circ \f^{\mu}$ is the
conformal welding associated to $[\mu]\in T(1)$, then
$\mathcal{S}(\f^{\mu})\in A_2(\Del)$ if and only if $[\mu]\in T_{0}(1)$. 
Let $j$ be the antiholomorphic inversion  $z\mapsto 1/\z$.
Since q.c.~mapping $w_{\mu}$ on $\hat{\C}$ satisfies $j\circ w_{\mu}\circ j = w_{\mu}$, we have
\begin{align*}
w_{\mu}^{-1} = j\circ w_{\mu}^{-1}\circ j = (j \circ
(\f^{\mu})^{-1}\circ j)\circ (j \circ \g_{\mu}\circ j).
\end{align*}
Thus 
\begin{equation} \label{inversion}
\f^{\mu^{-1}} = r\circ j \circ \g_{\mu}\circ j\quad\text{and}\quad
\g_{\mu^{-1}}= r \circ j \circ \f^{\mu}\circ j,
\end{equation}
where $r$ is the dilation $z\mapsto \ov{g_{\mu}'(\infty)}\,z$. Since $[\mu^{-1}]\in T_{0}(1)$ if and only $[\mu]\in T_{0}(1)$, we have $\mathcal{S}(\f^{\mu})\in A_{2}(\Del)$ if and only if $\mathcal{S}(\f^{\mu^{-1}})\in A_2(\Del)$, and hence if and only if
$$\mathcal{S}(\g_{\mu})=\ov{\mathcal{S}(\f^{\mu^{-1}})\circ j\,
j_{\bar{z}}^{2}}\in A_2(\Del^*).$$
The statement of the theorem now follows from Lemmas \ref{sub1} and \ref{sub2}.
\end{proof}

Let $\mathcal{T}_{0}(1)$ be the Teichm\"uller curve of
$T_{0}(1)$, i.e., the inverse image of $T_{0}(1)$ under the
fibration $\mathcal{T}(1)\rightarrow T(1)$ of Hilbert manifolds. It was proved in Appendix A of Part I that $\mathcal{T}_{0}(1)$ is a topological group. Using proofs of Lemma \ref{sub1} and Theorem \ref{sub2}, we can easily modify the proof in the Appendix of \cite{Teo} to show that $A^{1}_{2}(\Del)$ and $A_2(\Del)\oplus \C$ induce the same Hilbert manifold structure on $\mathcal{T}_{0}(1)$. We leave the details to Appendix A.

Results of this section justify the following
\begin{definition} The ``universal Liouville action'' $\SSS_{1}: T_{0}(1)\rightarrow \R$ is defined by
\begin{align*}
\SSS_1([\mu]) =
\iint\limits_{\Del}\left|\mathcal{A}(\f^{\mu})\right|^2d^2z+
\iint\limits_{\Del^*} \left|\mathcal{A}(\g_{\mu})\right|^2 d^2z
-4\pi\log|\g_{\mu}'(\infty)|,
\end{align*}
where $w_{\mu}=\g_{\mu}^{-1}\circ\f^{\mu}$ is the conformal welding corresponding to $[\mu]\in T_{0}(1)$.
\end{definition}
We will prove in Section 5 that the universal Liouville action is a K\"{a}hler potential of the Weil-Petersson metric on $T_{0}(1)$.
\begin{remark}
When $\g^{\prime}_{\mu}$ is continuous on $S^1$, the last term in the definition of $\SSS_{1}$ can be written as
\begin{align*}
-2\oint_{S^1}\log|\g_{\mu}'(e^{i\theta})|d\theta.
\end{align*}
When the quasicircle $\g_{\mu}(S^{1})=\f^{\mu}(S^{1})$ is of $C^{3}$ class, functionals of this type were studied by Schiffer and Hawley in  \cite{Schiffer3}. Here we extend the definition to quasicirles for the Hilbert manifold $T_{0}(1)$.
\end{remark}

\section{Grunsky operators for $T_{0}(1)$}
\label{GIGOFD}
\subsection{Grunsky coefficients and operators} \label{GCAP}
Here we prove that Grunsky operators associated to a point in $T_{0}(1)$ are Hilbert-Schmidt.
Suppose that $f:\Del\rightarrow \C$ and $g:\Del^*\rightarrow \C$
are univalent functions on $\Del$ and $\Del^*$ such that
$f(0)=0, f'(0)=1$, $g(\infty)=\infty$, and $f(\Del)\cap g(\Del^*)=\emptyset$. Such univalent functions are said to form a 
normalized disjoint pair.
The generalized Grunsky coefficients $b_{n,m}$, $n,m\in \Z$ of a normalized disjoint pair $(f,g)$ are 
defined as follows (see e.g., \cite{ Pom})
\begin{align*}
\log \frac{g(z)-g(\zeta)}{z-\zeta}&=
b_{00}-\sum_{m=1}^{\infty}\sum_{n=1}^{\infty} b_{mn} z^{-m}\zeta^{-n},\\
\log \frac{g(z) - f(\zeta)}{bz} &= -\sum_{m=1}^{\infty}
\sum_{n=0}^{\infty} b_{m,-n} z^{-m} \zeta^n,\\
\log \frac{f(z)-f(\zeta)}{z-\zeta} &= b_{00}-\sum_{m=0}^{\infty}
\sum_{n=0}^{\infty} b_{-m, -n} z^{m} \zeta^n.
\end{align*}
By definition, $b_{00} = \log b$, where
$b=g'(\infty)$. Grunsky coefficients $b_{n,m}$ are symmetric in
$n, m$ when $n,m \geq 1$ or $n,m \leq 0$, so for $n\geq 0, m
\geq 1$, we define $b_{-n, m}= b_{m,-n}$.  It is also clear that 
coefficients $b_{n, m}$, $|n|\geq 1$ and $|m|\geq 1$ do not changed when
$f$ and $g$ are simultaneosly post-composed with a dilation 
$z\mapsto rz$.

Grunsky coefficients satisfy the generalized Grunsky inequality, due to 
Hummel \cite{Hummel} (see also \cite{Pom}).
\begin{theorem} Let $(f,g)$ be a normalized disjoint pair of univalent functions. Then for every $\lambda_{-m},\dots,\lambda_{m}\in\CC$, 
\begin{align*}
\sum_{k=-\infty}^{\infty} |k| \left|\sum_{l=-m}^m b_{kl}
\lambda_l\right|^2  \leq \sideset{}{^{\prime}}\sum_{k=-m}
^{m}\frac{|\lambda_k|^2}{|k|} +2\re \left[ \bar{\lambda}_0 \sum_{l=-m}^m
b_{0l}\lambda_l\right],
\end{align*}
where the prime over the sum indicates that the term $k=0$ is omitted.
The equality for all $\lambda_{-m},\dots,\lambda_{m}$ holds if and only if the set
$F=\C\setminus\{f(\Del)\cup g(\Del^*)\}$ has Lebesgue measure zero.
\end{theorem}
\begin{remark}
For $\gamma\in\mathcal{T}(1)$ let $\gamma =g^{-1}\circ f$ be 
the corresponding conformal welding. 
Since $(f,g)$ is a normalized disjoint pair of univalent functions
and the quasicircle $\mathcal{C}=f(S^1)=g(S^1)$ has
Lebesgue measure zero, corresponding Grunsky coefficients $b_{mn}$ satisfy the generalized Grunsky equality.  
Setting $\lambda_{0}=1$ and $\lambda_{k}=0, k\neq 0$, we get
\begin{displaymath}
2\re b_{00}=\sum_{k=-\infty}^{\infty} |k||b_{k0}|^2.
\end{displaymath}
Since
\begin{displaymath}
\log \frac{g(z)}{z}=b_{00}-\sum_{k=1}^{\infty} b_{k0}z^{-k},
\hspace{1cm}\text{and}\hspace{1cm}\log
\frac{f(z)}{z}=-\sum_{k=1}^{\infty} b_{-k,0}z^{k},
\end{displaymath}
and $\re b_{00}= \log|g'(\infty)|$, we have
\begin{align*}
2\pi \log |g'(\infty)| = \iint\limits_{\Del} \left|
\frac{f'(z)}{f(z)}-\frac{1}{z}\right|^2d^2z+\iint\limits_{\Del^*}
\left| \frac{g'(z)}{g(z)}-\frac{1}{z}\right|^2d^2z.
\end{align*}
According to Theorem 5.3 in Part I, this gives an 
integral formula for the K\"{a}hler potential of the Velling-Kirillov 
metric on $\mathcal{T}(1)$.
\end{remark}
Now let $(f,g)$ be a normalized disjoint pair of univalent functions such that the corresponding set $F$ has Lebesgue measure zero. Putting in the generalized Grunsky 
equality $\lambda_0=0$ and rescaling $\lambda_l\mapsto \sqrt{|l|}\lambda_l$, 
we obtain the following equality
\begin{align*}
\sum_{k=-\infty}^{\infty}  \left|\sideset{}{^{\prime}}\sum_{l=-m }^m \sqrt{|kl|}b_{kl}
\lambda_l\right|^2 
= \sideset{}{^{\prime}} \sum_{k=-m}^{m} |\lambda_k|^2.
\end{align*}
By polarization, we get
\begin{align}\label{identities1}
\sideset{}{^{\prime}}\sum_{k=-\infty}^{\infty}\,\,&\sideset{}{^{\prime}}
\sum_{l=-m}^{m}\,\,\sideset{}{^{\prime}}\sum_{l'=-m}^{m}\sqrt{|kl|} b_{kl} 
\sqrt{|kl'|}\ov{b_{kl'}} \lambda_{l}\bar{\eta}_{l'} =\sideset{}{^{\prime}}\sum_{k=- m
}^{m} \lambda_k\bar{\eta}_k,
\end{align}
where $\lambda_k, \eta_k$ are arbitrary complex numbers. Grunsky coefficients $b_{mn}$ 
give rise to semi-infinite matrices $B_l$, $l=1,2,3,4$, defined by
\begin{align*}
(B_1)_{mn}=&\sqrt{mn}\;b_{-m,-n}, \hspace{0.3cm}
(B_2)_{mn}=\sqrt{mn}\;b_{-m,n},\\
(B_3)_{mn}=&\sqrt{mn}\;b_{m,-n},
\hspace{0.5cm}(B_4)_{mn}=\sqrt{mn}\;b_{mn}.
\end{align*}
From generalized Grunsky equality it immediately follows that
matrices $B_l$ define bounded linear operators on 
the separable Hilbert space
\begin{displaymath}
\ell^{2}=\left\{ x=\{x_{n}\}_{n=1}^{\infty}\,\,:\,\,
\sum_{n=1}^{\infty}|x_n|^2 <\infty\right\}
\end{displaymath}
which we continue to denote by $B_{l}$, $l=1,2,3,4$.
Here a linear operator $A$ on $\ell^{2}$ associated with the matrix $\{a_{mn}\}_{m,n=1}^{\infty}$ is given by $y=Ax$, where
$y_{m}=\sum_{n=1}^{\infty}a_{mn}x_{n}$.

In terms of the operators $B_{l}$, generalized Grunsky equality \eqref{identities1} is equivalent to
\begin{align}\label{relation1}
B_1 B_1^* + B_2 B_2^* = I,\hspace{1cm} B_3 B_1^* + B_4 B_2^*=0,\\
B_1 B_3^* + B_2 B_4^* = 0,\hspace{1cm}  B_3 B_3^* + B_4B_4^*=I,\nonumber
\end{align}
where $I$ is the identity operator on $\ell^{2}$ and $B_{l}^{*}$ stands for the adjoint operator to $B_{l}$.
These identities immediately imply that $\Vert B_l\Vert\leq 1,\,l=1,2,3,4.$
\begin{remark}\label{classicalgrunsky}
The operator $B_4$ is the Grunsky operator associated
to the univalent function $g$. The classical Grunsky inequality
(see e.g. \cite{Pom}) can be succintly stated as $I-B_4B_4^*\geq 0$,
and $I - B_{4}B_{4}^{\ast}$ is a  positive-definite operator if
and only if the complement of $g(\Del^*)$ has positive Lebesgue measure. 
Similarly, $B_{1}$ is the Grunsky operator associated 
to the univalent function $f$ and the classical Grunsky inequality
is equivalent to $I-B_{1}B_{1}^*\geq 0$.
For the pair $(\f^{\mu},\g_{\mu})$ associated to a 
point $[\mu]\in T(1)$, the operators $I-B_{1}B_{1}^{\ast}$ and 
$I-B_4B_4^*$ are positive-definite, so that $\Vert B_{1}\Vert, \Vert B_4\Vert <1$ and $\Ker B^{*}_{2}=\Ker B^{*}_{3} =\{0\}$. Moreover, it follows from symmetry property of Grunsky coefficients that also
$\Ker B_{2}=\Ker B_{3} =\{0\}$, so that the operators $B_{2}, B_{3}:\ell^{2}\rightarrow \ell^{2}$ are topological isomorphisms.
\end{remark}

The operators $B_{l}$ define
a bounded linear operator $\mathbf{B}$ on the Hilbert space $\ell^{2}\oplus\ell^{2}$ by
\begin{align*}
\mathbf{B}= \begin{pmatrix}
     B_1 & B_2\\
     B_3 & B_4
     \end{pmatrix}.
\end{align*}
Since
\begin{align*}
\mathbf{B}^{\ast}= \begin{pmatrix}
     B_{1}^{\ast} & B_{3}^{\ast}\\
     B_{2}^{\ast} & B_{4}^{\ast}
     \end{pmatrix},
\end{align*}
the generalized Grunsky equality can be succinctly rewritten as
\begin{align*}
 \mathbf{B}\mathbf{B}^{\ast} =\mathbf{I},
\end{align*}
where $\mathbf{I}=\left(\begin{smallmatrix} 
I & 0 \\ 0 & I \end{smallmatrix}\right)$ is the identity operator on $\ell^{2}\oplus\ell^{2}$. 
Let $J$ be the complex-conjugation operator on $\ell^{2}$ defined by
\begin{equation} \label{conjugation}
(Jx)_{n}=\bar{x}_{n},\quad x=\{x_{n}\}_{n=1}^{\infty}\in\ell^{2}.
\end{equation}
Setting  $\mathbf{J}=\left(\begin{smallmatrix} 
J & 0 \\ 0 & J \end{smallmatrix}\right)$, we can express symmetry property of Grunsky coefficients as
$$\mathbf{B}^{\ast} =\mathbf{J}\mathbf{B}\mathbf{J}.$$
Thus
$$\mathbf{B}^{*}\mathbf{B}=\mathbf{J} \mathbf{B}\mathbf{J}\mathbf{B}=\mathbf{J}\mathbf{B}\mathbf{B}^{*}\mathbf{J} =\mathbf{I},$$
so that $\mathbf{B}$ is a unitary operator on $\ell^{2}\oplus\ell^{2}$.

The operators $B_l $ can be also realized as linear  operators from the Hilbert spaces of antiholomorphic functions to the Hilbert spaces of holomorphic functions. Namely, the kernels 
\begin{align*}
K_{1}(z,w) &= \frac{1}{\pi}\left(\frac{1}{(z-w)^2}-
\frac{f'(z)f'(w)}{(f(z)-f(w))^2}\right)=\frac{1}{\pi}
\sum_{n,m=1}^{\infty}nmb_{-n,-m}z^{n-1}w^{m-1},\\
K_2(z,w)&= \frac{1}{\pi}\frac{f'(z)
g'(w)}{(f(z)-g(w))^2}=\frac{1}{\pi}
\sum_{n,m=1}^{\infty}nmb_{-n,m}z^{n-1}w^{-m-1},\\
K_{3}(z,w) &= \frac{1}{\pi}
\frac{g'(z)f'(w)}{(g(z)-f(w))^2}=\frac{1}{\pi}
\sum_{n,m=1}^{\infty}nmb_{n,-m}z^{-n-1}w^{m-1},\\
K_4(z,w)&=\frac{1}{\pi}\left( \frac{1}{(z-w)^2}-\frac{g'(z)
g'(w)}{(g(z)-g(w))^2}\right)=\frac{1}{\pi}
\sum_{n,m=1}^{\infty}nmb_{n,m}z^{-n-1}w^{-m-1},
\end{align*}
define the linear operators $K_l$ as follows,
\begin{alignat*}{3}
K_1 &: \ov{A^{1}_{2}(\Del)} \rightarrow A^{1}_{2}(\Del),
&\qquad (K_1 \psi)(z) &= \iint\limits_{\Del}
K_1(z,w)\ov{\psi(w)}d^2w,\\
K_2&: \ov{A^{1}_{2}(\Del^*)}\rightarrow A^{1}_{2}(\Del),
&\qquad (K_2\psi)(z)&= \iint\limits_{\Del^*}
K_2(z,w)\ov{\psi(w)}d^2w,\\
K_3&: \ov{A^{1}_{2}(\Del)} \rightarrow A^{1}_{2}(\Del^*),
&\qquad (K_3\psi)(z)&= \iint\limits_{\Del}
K_3(z,w)\ov{\psi(w)}d^2w,\\
K_4&: \ov{A^{1}_{2}(\Del^*)} \rightarrow A^{1}_{2}(\Del^*),
&\qquad (K_4\psi)(z)&= \iint\limits_{\Del^*}
K_4(z,w)\ov{\psi(w)}d^2w.
\end{alignat*}
\begin{remark}\label{singular} It is well-known that if $\phi$ is a 
holomorphic function on $\Del$, then 
\begin{align*}
\iint\limits_{\Del} \frac{\ov{\phi(w)}}{(z-w)^2}d^{2}w=0,
\end{align*}
where the integral is understood in the principal value sense. 
Hence we can also represent operators $K_1$ and $K_4$ by the singular kernels
\begin{align*}
-\frac{1}{\pi}\frac{f'(z)f'(w)}{(f(z)-f(w))^2}\quad\text{and}\quad
-\frac{1}{\pi}\frac{g'(z)g'(w)}{(g(z)-g(w))^2}.
\end{align*}
\end{remark}
The Hilbert spaces $A^{1}_{2}(\Del)$ and $A^{1}_{2}(\Del^*)$ 
have standard orthonormal bases $\{e_{n}\}_{n=1}^{\infty}$ 
and $\{f_{n}\}_{n=1}^{\infty}$, given respectively by
\begin{align*}
e_n(z) = \sqrt{\frac{n}{\pi}}\,z^{n-1}\quad\text{and}\quad 
f_n(z) = \sqrt{\frac{n}{\pi}}\,z^{-n-1},\quad n\in\mathbb{N}.
\end{align*}
These bases define isomorphisms  $A^{1}_{2}(\Del)\simeq\ell^{2}$ and  $A^{1}_{2}(\Del^{\ast})\simeq\ell^{2}$. The operators 
$K_l$ and their adjoints $K_{l}^{*}$ ---  integral operators with the kernels $K_{l}^{*}(z,w)=\ov{K_{l}(w,z)}$, correspond respectively to the operators $B_l$ and $B_{l}^{\ast}$, $l=1,2,3,4$.
Similarly, positive self-adjoint operators
$\KKK_l = K_l  K_l^*$ are integral operators 
which correspond to the operators $B_{l}B^{*}_{l}$, and
we denote the kernels of the operators $\KKK_{l}$ by 
$\KKK_{l}(z,w)$.
Due to the relations \eqref{relation1},
\begin{equation} \label{K=K} 
\KKK_{2} = I-\KKK_{1},\,\,\,\KKK_{3}= I-\KKK_{4}.
\end{equation}
\begin{lemma} \label{H-S} 
The kernel $K_1(z,w)$ of the operator $K_1: \ov{A^{1}_{2}(\Del)}\rightarrow A^{1}_{2}(\Del)$ satisfies
\begin{align} \label{hs}
\iint\limits_{\Del}\iint\limits_{\Del}|K_1(z,w)|^{2}d^2zd^{2}w<\infty
\end{align}
if and only if the operator $K_1$ is Hilbert-Schmidt, i.e., if and only if the operator
$\KKK_{1}=K_{1}K_{1}^{*}$ on $A^{1}_{2}(\Del)$ is of trace class.
In this case,
$$\Tr \KKK_{1} = \iint\limits_{\Del}\iint\limits_{\Del}|K_1(z,w)|^{2}d^2zd^{2}w =\iint\limits_{\Del}\KKK_1(z,z) d^2z,$$
and $\mathcal{S}(f)\in A_2(\Del)$, where $f$ is the univalent function associated with the kernel $K_{1}(z,w)$.
Similar statements hold for the operators $K_{4}$ and $\KKK_4$.
\end{lemma}
\begin{proof}
It is sufficient to prove the lemma for the operator $\KKK_{1}$.
For the basis $\{e_{n}\}_{n\in\mathbb{N}}$ of the Hilbert space $A^{1}_{2}(\Del)$ we have
\begin{align*}
\Tr\KKK_{1} = & \sum_{n=1}^{\infty}\la\KKK_{1}e_{n},e_{n}\ra =\sum_{n=1}^{\infty}\Vert K^{*}_{1}e_{n}\Vert^{2} = \sum_{n,m=1}^{\infty}nm|b_{-n,-m}|^2 \\
= & \iint\limits_{\Del}\iint\limits_{\Del} \left|K_1(z,w)\right|^2d^2zd^2w =
\iint\limits_{\Del}\KKK_1(z,z) d^2z.
\end{align*}
Since the operator $\KKK_{1}$ is positive, it is of trace class if and only if the inequality \eqref{hs} holds.
On the other hand, we have
\begin{align*}
\mathcal{S}(f)(z) =-6\pi \lim_{w\rightarrow z} K_1(z,w) = -6\sum_{n=2}^{\infty} \left(\sum_{k+l=n} kl
b_{-k,-l}\right)z^{n-2}.
\end{align*}
Hence if the inequality \eqref{hs} holds, 
\begin{align*}
\Vert \mathcal{S}(f)\Vert_2^2 = &18\pi \sum_{n=2}^{\infty}
\frac{1}{n^3-n} \left|\sum_{k=1}^{n-1} k(n-k) b_{-k,-(n-k)}\right|^2 \\
\leq & 18\pi \sum_{n=2}^{\infty}\frac{1}{n^3-n}
\left(\sum_{k=1}^{n-1}k(n-k)\right)\left(\sum_{k=1}^{n-1}
k(n-k)|b_{-k,
-(n-k)}|^2\right)\\
=&3\pi \sum_{n=2}^{\infty} \sum_{k=1}^{n-1} k(n-k)|b_{-k,
-(n-k)}|^2=3\pi \sum_{n,m=1}^{\infty} nm|b_{-n,-m}|^2 <\infty.
\end{align*}
\end{proof}
\begin{theorem}\label{traceclass}
If the pair $(\f^{\mu}, \g_{\mu})$ corresponds to a point $[\mu]\in T_{0}(1)$, then the operators
$\KKK_1$ and $\KKK_{4}$ associated to $\f^{\mu}$ and $\g_{\mu}$ respectively, are of trace class.
\end{theorem}
\begin{proof}
According to Lemma \ref{H-S}, it is sufficient to show that
\begin{align*}
\iint\limits_{\Del}\KKK_1(z,z) d^2z<\infty\,\,\,\text{and}\,\,\,\iint\limits_{\Del^{*}}\KKK_4(z,z) d^2z<\infty.
\end{align*}
For $[\mu]\in T_{0}(1)$ choose a representative $\mu\in L^2(\Del^*, \rho(z)d^2z)\cap \mathcal{O}(\Del^*)_1$. It follows from Lemma 3.9 in Part I that the path $[t\mu]$ connecting $0$ to $[\mu]$ in $T(1)$ lies on $T_{0}(1)$. Let $w_{t\mu}=
\g_{t\mu}^{-1}\circ \f^{t\mu}$ be the corresponding conformal welding and
denote by $(K_{1})_{t}(z,w)$ the kernel $K_{1}(z,w)$ associated with the univalent function $\f^{t\mu}$.  
We have the following lemma.
\begin{lemma}\label{varyK1}
\begin{gather}
\left.\frac{d}{d s}\right\vert_{s=0} (K_1)_{s+t}\left(\f_t^{-1}(z),
\f_t^{-1}(w)\right)\left(\f_t^{-1}\right)'(z)\left(\f_t^{-1}\right)'(w)  \nonumber \\
= \frac{1}{\pi^2}\iint\limits_{ \Omega_t^*}
\frac{\mu_t(u)}{(u-z)^2(u-w)^2}d^2u , \label{vary-K1}
\end{gather} 
where $\Omega_{t}^*=\f^{t\mu}(\Del^*)=\g_{t\mu}(\Del^{*})$,
\begin{align*}
(\mu_t\circ \g_{t\mu}) \frac{\ov{\g'_{t\mu}}}{\g'_{t\mu}}= D_{t\mu}R_{(t\mu)^{-1}}(\mu),
\end{align*}
and the integral \eqref{vary-K1} is understood in the principal value sense.
\end{lemma}
\begin{proof} Set $w_{t}=w_{t\mu}, \f_{t}=\f^{t\mu}, \g_{t}=\g_{t\mu}$ and $v_{s}=\f_{s+t}\circ \f_{t}^{-1}$. We have
$$v_{s}\circ \g_{t} = \g_{s+t}\circ w_{s+t}\circ w_{t}^{-1},$$
so that $v_{s}$ is a q.c. mapping which is holomorphic on $\Omega_{t}=\f_{t}(\Del)$  and has Beltrami differential $\mu_{s,t}$ on $\Omega_{t}^{*}$ with
$$(\mu_{s,t}\circ \g_{t})\frac{\ov{\g^{\prime}_{t}}}{\g^{\prime}_{t}}=
\frac{(w_{s+t}\circ w_{t}^{-1})_{\z}}{(w_{s+t}\circ w_{t}^{-1})_{z}}.$$
It follows from the standard variational formula for q.c. mappings that
\begin{align}\label{variation}
\left.\frac{d}{ds}\right\vert_{s=0} v_s(z)
=-\frac{1}{\pi}\iint\limits_{\Omega_{t}^*}
\frac{\mu_t(u)z(z-1)}{(u-z)u(u-1)}d^2u+ p(z),
\end{align}
where $p(z)$ is a degree two polynomial. 
We have
\begin{align*}
&(K_1)_{s+t}\left(\f_t^{-1}(z),
\f_t^{-1}(w)\right)\left(\f_t^{-1}\right)'(z)\left(\f_t^{-1}\right)'(w)\\
=&\frac{1}{\pi}\frac{\left(\f_t^{-1}\right)'(z)\left(\f_t^{-1}\right)'(w)
}{\left(\f_t^{-1}(z)-\f_t^{-1}(w)\right)^2}-\frac{1}{\pi}\frac{v_s'(z)v_s'(w)}{(v_s(z)-v_s(w))^2},
\end{align*}
and
\begin{align*}
\left.\frac{d}{ds}\right|_{s=0}
\frac{v_s'(z)v_s'(w)}{(v_s(z)-v_s(w))^2}=-\frac{1}{\pi}
\iint\limits_{\Omega^*}\frac{\mu_t(u)}{(u-z)^2(u-w)^2}d^2u,
\end{align*}
so that the result follows.
\end{proof}

Now we use the fundamental theorem of calculus to estimate
\begin{align*}
\iint\limits_{\Del}\KKK_1(z,z)d^2z=
&\iint\limits_{\Del}\iint\limits_{\Del} \bigl|(K_1)_1(z,w)\bigr|^2 d^2zd^2w\\
=&\iint\limits_{\Del}\iint\limits_{\Del}\left|
 \int_0^1 \frac{d}{dt}(K_1)_t(z,w)dt\right|^2d^2zd^2w \\
\leq & \int_0^1 \iint\limits_{\Del}\iint\limits_{\Del}\left|
 \frac{d}{dt}(K_1)_t(z,w)\right|^2d^2zd^2wdt \\
=&\int_0^1 \iint\limits_{\Del}\iint\limits_{\Del}\left|
\frac{d}{ds}\Bigr\vert_{s=0}(K_1)_{t+s}(z,w)\right|^2\!d^2zd^2wdt\\ 
=&\int_{0}^{1}I(t)dt.
\end{align*}
Making a change of variables $z\mapsto \f^{-1}_{t}(z), w\mapsto \f^{-1}_{t}(w)$ in the inner integral $I(t)$, we get 
\begin{align*}
I(t)=
 & \iint\limits_{\Omega_t}\iint\limits_{\Omega_t} \left|
\left.\frac{d}{ds}\right|_{s=0} (K_1)_{t+s}\left(\f_t^{-1}(z),
\f_t^{-1}(w)\right)\left(\f_t^{-1}\right)'(z)\left(\f_t^{-1}\right)'(w)\right|^2
\!d^2zd^2w\\
=&\frac{1}{\pi^4}
\iint\limits_{\Omega_t}\iint\limits_{\Omega_t} \left|\,
\iint\limits_{ \Omega_t^*}
\frac{\mu_t(u)}{(u-z)^2(u-w)^2}d^2u\,\right|^2\!d^2zd^2w.
\end{align*}
Using the inequality 
\begin{equation} \label{inequality}
\iint\limits_{\Omega_{t}}\frac{d^{2}w}{|w-z|^{4}}\leq 4\pi(\rho_{2})_{t}(z),\,\,z\in\Omega^{*}_{t},
\end{equation}
where $(\rho_{2})_{t}(z)$ is the density of the hyperbolic metric on $\Omega^{*}_{t}$ (see the proof of Theorem 3.3 in Part I), and the fact that the Hilbert transform is an isometry on
$L^2(\C, d^2z)$, we obtain
\begin{align*}
I(t)\leq &\frac{1}{\pi^{2}}\iint\limits_{\Omega_{t}}\iint\limits_{\Omega_t^*} \frac{ |
\mu_t(z)|^2}{|z-w|^4} d^2zd^{2}w\leq \frac{4}{\pi} \iint\limits_{\Omega_t^*}
|\mu_t(z)|^2 (\rho_2)_t(z)d^2z\\=&\frac{4}{\pi}
\iint\limits_{\Del^*}|\tilde{\mu}_t(z)|^2 \rho(z)d^2z =\frac{4}{\pi} \Vert
\tilde{\mu}_t\Vert_2^2,
\end{align*}
where $\tilde{\mu}_t = D_{t\mu}R_{(t\mu)^{-1}}(\mu)$.
Now it follows from Remark 3.8 in Part I that there exists a constant $C$ such that 
\begin{align*}
\Vert \tilde{\mu}_t\Vert_2 \leq C \Vert \mu \Vert_2
\end{align*}
for all $0\leq t\leq 1$, so that 
$$\iint\limits_{\Del}\KKK_1(z,z) d^2z<\infty.$$
The corresponding estimate for the kernel $\KKK_{4}(z,w)$ is proved similarly. Alternatively, using the relation \eqref{inversion}
we get
\begin{equation} \label{symmetry1-4}\KKK_{1}([\mu^{-1}])(z,w)=
\ov{\KKK_4([\mu])\left(\frac{1}{\z},
\frac{1}{\bar{w}}\right)}\,\frac{1}{z^2}\frac{1}{w^2}.
\end{equation}
Since $T_{0}(1)$ is a group, the inequality for the kernel $\KKK_{4}$ follows from the corresponding inequality for the kernel $\KKK_{1}$.
\end{proof}
\begin{remark} Actually using the generalized Grunsky equality one can prove an estimate sharper than \eqref{inequality}. Just observe that
for $z\in\Omega^{\ast}_{t}$
$$\iint\limits_{\Omega_{t}}\frac{d^{2}w}{|z-w|^{4}}=\pi^{2}\KKK_{3}(g^{-1}(z),g^{-1}(z))|(g^{-1})'(z)|^{2}$$
and that $\tfrac{1}{\pi(1-z\bar{w})^{2}}$ is the kernel of the identity operator on $A_{2}^{1}(\Del^{\ast})$.
Hence the second equation in \eqref{K=K} gives,
$$\KKK_{3}(z,z)= \frac{1}{\pi(1-|z|^{2})^{2}}-\KKK_{4}(z,z)\leq  \frac{1}{\pi(1-|z|^{2})^{2}},$$
and we get
$$\iint\limits_{\Omega_{t}}\frac{d^{2}w}{|z-w|^{4}}\leq\frac{\pi}{4}(\rho_{2})_{t}(z).$$
\end{remark}
\begin{corollary} \label{trace-grunsky} Grunsky operators $B_{1}$ and $B_{4}$ associated with the pair $(\f^{\mu},\g_{\mu}),\,[\mu]\in T(1)$, are Hilbert-Schmidt operators on $\ell^{2}$ if and only if $[\mu]\in T_{0}(1)$.
\end{corollary}
\begin{proof} Under the isomorphisms  $A^{1}_{2}(\Del)\simeq\ell^{2}$ and  $A^{1}_{2}(\Del^{\ast})\simeq\ell^{2}$, the operators $\KKK_{1}$ and $\KKK_{4}$ correspond to the operators $B_{1}B_{1}^{\ast}$ and $B_{4}B_{4}^{\ast}$ respectively. Since $\beta(T_{0}(1))=A_{2}(\Del)\cap\beta(T(1))$, the ``only if'' part of the statement follows from Lemma \ref{H-S}.
\end{proof}

As an application, consider the Hilbert space $\cS_{2}$ of
Hilbert-Schmidt operators on $\ell^{2}$,
\begin{align*}
\cS_{2}=\left\{T: \ell^{2} \rightarrow \ell^{2}\;\text{a
bounded operator}\;\;\Bigr|\; \Vert T \Vert_2^2 = \Tr TT^*
<\infty\right\},
\end{align*}
and define the mapping $\cP: T_{0}(1)\rightarrow \cS_{2}$ by
\begin{align*}
\cP([\mu])= B_1(\f^{\mu}), \,\,[\mu]\in T_{0}(1).
\end{align*}
Since Grunsky coefficients characterize univalent functions up to a post-composition with M\"{o}bius transformation, the mapping $\cP$ is one to one. In fact, we have a stronger result.
\begin{theorem}\label{holomorphicembedding}
The mapping $\cP$ is a holomorphic inclusion
of the Hilbert manifold $T_{0}(1)$ into the Hilbert space
$\cS_{2}$. 
\end{theorem}
\begin{proof}
We need to show that for every $[\nu]\in T_{0}(1)$ and $\mu\in
H^{-1,1}(\Del^*)$, the map $\C\ni t \mapsto
B_1(t)=B_1(f^{\nu+t\mu})$ is holomorphic in a neighbourhood of
$t=0$ in $\C$.  For this aim, since the mapping $[\mu]\rightarrow \f^{\mu}(z)$ is holomorphic for fixed $z\in\Del$, for every $z,w\in\Del$ the map
\begin{align*}
t \mapsto
K_1^{\nu+t\mu}(z,w)=\frac{1}{\pi}\left(\frac{1}{(z-w)^2}-\frac{(\f^{\nu+
t\mu})'(z)(\f^{\nu+t\mu})'(w)}{(\f^{\nu+t\mu}(z)-\f^{\nu+t\mu}(w))^2}\right)
\end{align*}
is holomorphic in a neighbourhood of $t=0$ in $\C$. We choose
$\delta>0$ so that $\Vert \nu+t\mu\Vert_{\infty}<1$ for all
$|t|<\delta$. For every $t_0$ such that $|t_0|<\delta$, let
$\delta_1$ be such that $0<\delta_1< \delta-|t_0|$. Then for all
$|t-t_0|<\delta_1$, we have by Cauchy integral formula,
\begin{align*}
&\left(K_1^{\nu+t\mu}-K_1^{\nu+t_0\mu} -(t-t_0)
\left.\frac{d}{dt}\right|_{t=t_0}K_1^{\nu+t\mu}\right)(z,w)\\
=&\frac{(t-t_0)^2}{2\pi i}
\oint\limits_{|\zeta-t_0|=\delta_1}\frac{K_1^{\nu+\zeta\mu}(z,w)}{(\zeta-t)(\zeta-t_0)^2}d\zeta.
\end{align*}
Hence
\begin{align}\label{derivative}
&\left\Vert \frac{B_1(\f^{\nu+t\mu})-B_1(\f^{\nu+t_0\mu})}{t-t_0} -
\frac{d}{dt}\Bigr\vert_{t=t_0}B_1(\f^{\nu+t\mu})\right\Vert_2^2\\
=&\iint\limits_{\Del}\iint\limits_{\Del}\left|\left(\frac{K_1^{\nu+t\mu}-K_1^{\nu+t_0\mu}}
{t-t_0}-
\left.\frac{d}{dt}\right|_{t=t_0}K_1^{\nu+t\mu}\right)(z,w)\right|^2d^2zd^2w\nonumber\\
\leq & \frac{|t-t_0|^2}{4\pi^2}
\oint\limits_{|\zeta-t_0|=\delta_1}\iint\limits_{\Del}\iint\limits_{\Del}
\left|K_1^{\nu+\zeta\mu}(z,w)\right|^2d^2zd^2w |d\zeta|
\oint\limits_{|\zeta-t_0|=\delta_1}\frac{|d\zeta|}{|\zeta-t|^2|\zeta-t_0|^4}.\nonumber
\end{align}
We have from the proof of Theorem \ref{traceclass},
\begin{align*}
\iint\limits_{\Del}\iint\limits_{\Del}
\left|K_1^{\nu+\zeta\mu}(z,w)\right|^2d^2zd^2w\leq C \Vert
\nu+\zeta\mu\Vert_2^2\leq C(\Vert \nu\Vert_2+\delta_1\Vert
\mu\Vert_2)^2,
\end{align*}
so that \eqref{derivative} tends to $0$ as $t\rightarrow t_0$, which
proves the assertion.
\end{proof}
\begin{remark}\label{Gperiod}
Since the classical Grunsky operator $B_1$ is bounded, the mapping $\cP$ extends to the whole Banach manifold $T(1)$. Let  $\cB(\ell^2)$ be the space of bounded linear
operators on $\ell^2$,
\begin{align*}
\cB(\ell^2)=\left\{ T:\ell^2\rightarrow \ell^2 \;\text{a
linear operator}\; : \Vert T\Vert=\sup_{\Vert u
\Vert=1} \Vert Tu\Vert<\infty.\right\},
\end{align*}
and define the mapping $\hat{\cP}: T(1)\rightarrow
\cB(\ell^2)$ by
\begin{displaymath}
\hat{\cP}([\mu])=B_1(\f^{\mu}),\quad [\mu]\in T(1).
\end{displaymath}
Analogous to Theorem \ref{holomorphicembedding}, we
show in Appendix B that the mapping $\hat{\cP}$ is a
holomorphic inclusion. 
\end{remark}
\subsection{Fredholm eigenvalues and Fredholm determinant}
In \cite{Schiffer1}, Schiffer has studied the eigenvalues of the classical Poincar\'{e}-Fredholm boundary value problem of potential theory on a $C^3$ curve. Here we show how Fredholm eigenvalues for a quasi-circle $\mathcal{C}=\f^{\mu}(S^{1})=\g_{\mu}(S^{1})$, 
associated with $[\mu]\in T_{0}(1)$, are related to the eigenvalues of trace class operators $\KKK_1$ and $\KKK_{4}$. 

Let $\mathfrak{h}$ be a separable Hilbert space with the inner product $\la~,~\ra$. A conjugation operator $J$ on $\mathfrak{h}$ is an $\RR$-linear operator satisfying $J^{2}=I$ and
$$(Jx,Jy)=\ov{(x,y)}\quad\text{for all}\quad x,y\in\mathfrak{h}.$$
Conjugation operator is necessarily complex anti-linear. For every bounded linear operator $T$ on $\mathfrak{h}$,
$$\la JTJx,y\ra = \ov{\la TJx,Jy\ra}=\ov{\la Jx,T^{*}Jy\ra} =\la x,JT^{*}Jy\ra\quad\text{for all}\quad x,y\in \mathfrak{h},$$
so that
$$(JTJ)^{*} =JT^{*}J.$$
In particular, if $U$ is a unitary operator on $\mathfrak{h}$, then $JUJ$ is also a unitary operator.
For a bounded linear operator $T$ on $\mathfrak{h}$ its transpose is defined as
$$T^{t} =JT^{*}J.$$
Generalizing the notion of symmetric complex-valued matrix,
a bounded operator $T$ on $\mathfrak{h}$ is called symmetric
with respect to the conjugation $J$, if
$$T=T^{t}.$$
The Hilbert space $\mathfrak{h}=\ell^{2}$ carries a standard conjugation operator $J$, defined by \eqref{conjugation}. 
The following statement is a generalization of Schur's Lemma (see, e.g., \cite[Sect.~3.6]{Pom}) to the case of compact operators on $\ell^{2}$.
\begin{lemma}
Let $T$ be a compact operator on $\ell^{2}$, symmetric with respect to the standard conjugation operator $J$. Then there exist a unitary operator $U$ on $\ell^{2}$ and an operator $D\geq 0$ on $\ell^{2}$,  diagonal with respect to the standard basis for $\ell^{2}$, such that
$$T=UDU^{t}.$$
\end{lemma}
\begin{proof}
As in \cite{Pom}, consider the decomposition
$$T=\frac{T+JTJ}{2} +i\frac{T-JTJ}{2i} =A + iB,$$
where $A$ and $B$ are self-adjoint compact operators satisfying $AJ=JA$ and
$BJ=JB$. Let $\mathbf{T}$ be the self-adjoint operator on the Hilbert space $\ell^{2}\oplus\ell^{2}$ defined by
\begin{align*}
\mathbf{T}= \begin{pmatrix}
     A & B\\
     B & -A
     \end{pmatrix}.
\end{align*}
The operator $\mathbf{T}$ is compact and satisfies
\begin{equation}\label{properties}
\mathbf{T}\mathbf{E}= - \mathbf{E} \mathbf{T}\quad\text{and}\quad
\mathbf{T}\mathbf{J}=\mathbf{J}\mathbf{T}, 
\end{equation}
where
\begin{align*}
\mathbf{E}= \begin{pmatrix}
     0 & - I\\
     I & 0
     \end{pmatrix}\quad\text{and}\quad\mathbf{J}=\begin{pmatrix}
     J & 0\\
     0 & J
     \end{pmatrix}.
\end{align*}
From the first equation in \eqref{properties} it follows that if 
$\mathbf{u}\in\ell^{2}\oplus\ell^{2}$
is an eigenvector for $\mathbf{T}$ with eigenvalue $\lambda$, then 
$\mathbf{v} =\mathbf{E}\mathbf{u}$ 
is also an eigenvector for $\mathbf{T}$ with eigenvalue $-\lambda$. It follows from Hilbert-Schmidt theorem on canonical form of compact self-adjoint operator that there
exist a unitary operator  $\mathbf{U}$ on $\ell^{2}\oplus\ell^{2}$ of the form
\begin{align*}
\mathbf{U}= \begin{pmatrix}
     U_{1} & U_{2}\\
     U_{2} & -U_{1}
     \end{pmatrix},
\end{align*}
and an operator $\mathbf{D}$ on $\ell^{2}\oplus\ell^{2}$ of the form
\begin{align*}
\mathbf{D}= \begin{pmatrix}
     D & 0\\
     0 & -D
     \end{pmatrix},
\end{align*}
where $D$ is diagonal with non-negative entries, such that
$$\mathbf{T}=\mathbf{U}\mathbf{D}\mathbf{U}^{\ast}.$$
From the second equation in \eqref{properties} it follows that $\mathbf{T}(\mathbf{J}\mathbf{U}\mathbf{J})=(\mathbf{J}\mathbf{U}\mathbf{J})\mathbf{D}$. Since $\mathbf{J}\mathbf{U}\mathbf{J}$ is also a unitary operator, we have
$$\mathbf{T}=(\mathbf{J}\mathbf{U}\mathbf{J})\mathbf{D}(\mathbf{J}\mathbf{U}\mathbf{J})^{*}.$$
Consequently, we can choose $\mathbf{U}$ so that 
$\mathbf{U}=\mathbf{J}\mathbf{U}\mathbf{J}$.
Now it follows from the canonical form that
$$T=A+iB =(U_{1} + i U_{2})D(U_{1}^{\ast} + i U_{2}^{\ast}).$$
Let $U=U_{1} + i U_{2}:\ell^{2}\rightarrow\ell^{2}$.
Since $\mathbf{U}$ is a unitary operator, $U$ is also unitary, and the property $\mathbf{U}^{*}=\mathbf{J}\mathbf{U}^{*}\mathbf{J}$ implies that
$$JU^{\ast}J=J(U_{1}^{*} - i U_{2}^{*})J =U_{1}^{*} + iU_{2}^{*},$$
since $J$ is complex anti-linear.
\end{proof}
\begin{corollary} The non-zero entries of the operator $D$ are singular values of the operator $T$.
\end{corollary}
\begin{proof} Since the operator $U^{t}$ is unitary,
$$TT^{\ast}=UD^{2}U^{\ast}=UD^{2}U^{-1},$$
so that the entries of $D^{2}$ are the eigenvalues of $TT^{\ast}$.
\end{proof}

Now let $(f,g)$ be a normalized disjoint pair of univalent functions such that the corresponding set $F$ has Lebesgue measure zero and the Grunsky operator $B_{1}$ is compact. We apply Schur's Lemma to the operator $B_1$ on $\ell^{2}$. It follows from the symmetry property of Grunsky coefficients that 
$$B_{1}^{\ast}=JB_{1}J.$$
Thus there exist a unitary operator $U$ on $\ell^{2}$ and a diagonal operator  $D$ with non-negative entries such that
$$B_{1}=UDU^{t}.$$
From the first identity in \eqref{relation1}, we obtain
\begin{align*}
U^{-1}B_2B_2^*U= I -D^2.
\end{align*}
Since $\Vert B_{1}\Vert<1$, the operator $I-D^{2}$ is
positive-definite and hence invertible, so that the operator
$$V= B_{2}^{\ast} U(I- D^2)^{-1/2}$$
is also unitary. Using the property $B_{3}^{t}=B_{2}$, which follows from the symmetry of Grunsky coefficients, and the third identity in \eqref{relation1}, we obtain
$$V^{t}B_{4} V =- D.$$ 
Collecting everything together, we get the following identities:
\begin{alignat*}{3}
B_1 JUJ& =  UD, & \qquad B_3 JUJ & = J VJ(I-D^2)^{1/2},\\
B_2 V & = U(I-D^2)^{1/2}, & \qquad B_4V  & =  -JV JD.\nonumber
\end{alignat*}
Letting
\begin{align*}
\lambda_n &= (D)_{nn},\hspace{2cm}
\rho_n=((1-D^2)^{1/2})_{nn}=\sqrt{1-\lambda_n^2}\\
\mathfrak{u}_n
(z)&=\sum_{m=1}^{\infty}\sqrt{\frac{m}{\pi}}U_{mn}z^{m-1},
\hspace{0.5cm}\mathfrak{v}_n
(z)=\sum_{m=1}^{\infty}\sqrt{\frac{m}{\pi}}(JV J)_{mn}z^{-m-1},
\end{align*}
and realizing $B_l$'s as linear operators $K_l$'s, we obtain for $n\in\mathbb{N}$,
\begin{alignat*}{3}
\iint\limits_{\Del} K_1(z,w) \ov{\mathfrak{u}_n(w)}d^2w&=\lambda_n
\mathfrak{u}_n(z), &\quad \iint\limits_{\Del} K_3(z,w)
\ov{\mathfrak{u}_n(w)}d^2w&=\rho_n
\mathfrak{v}_n(z)\\
\iint\limits_{\Del^*} K_2(z,w) \ov{\mathfrak{v}_n(w)}d^2w&=\rho_n
\mathfrak{u}_n(z), &\quad \iint\limits_{\Del^*} K_4(z,w)
\ov{\mathfrak{v}_n(w)}d^2w&=-\lambda_n \mathfrak{v}_n(z).
\end{alignat*}
Setting
\begin{align*}
u_n = \mathfrak{u}_n\circ f^{-1}
(f^{-1})'\hspace{1cm}\text{and}\hspace{1cm}
v_n=\mathfrak{v}_n\circ g^{-1}(g^{-1})',
\end{align*}
we get,
\begin{alignat}{3} \label{Fredholmintegral}
\frac{1}{\pi}\iint\limits_{\Omega}\frac{\ov{u_n(w)}}{(z-w)^2}d^2w & =
-\lambda_n u_n(z), &\quad & z\in \Omega,\\
\frac{1}{\pi}\iint\limits_{\Omega}\frac{\ov{u_n(w)}}{(z-w)^2}d^2w & =
\rho_n v_n(z), & \quad & z\in \Omega^*, \nonumber\\
\frac{1}{\pi}\iint\limits_{\Omega^*}\frac{\ov{v_n(w)}}{(z-w)^2}d^2w & =
\rho_n u_n(z), & \quad & z\in \Omega,\nonumber\\
\frac{1}{\pi}\iint\limits_{\Omega^*}\frac{\ov{v_n(w)}}{(z-w)^2}d^2w & =
\lambda_n v_n(z), &\quad & z\in \Omega^*.\nonumber
\end{alignat}
Comparing equations \eqref{Fredholmintegral} with corresponding formulas in \cite{Schiffer1}, we find that
$\{\pm\lambda_n^{-1}\}_{n=1}^{\infty}$ are Fredholm eigenvalues
associated to the quasi-circle $\mathcal{C}=f(S^{1})=g(S^{1})$. 
\begin{remark}
The relation between the Fredholm eigenvalues and the eigenvalues
of the Grunsky operator for a $C^3$ curve was first obtained by
Schiffer in \cite{Schiffer2}.
Specifically, in \cite{Schiffer2} Schiffer has shown that Fredholm eigenvalues, defined as the eigenvalues of classical Poincar\'{e}-Fredholm integral operator on $C^{3}$ curve,  satisfy \eqref{Fredholmintegral}. Furthermore, using completeness of
the bases $\{u_n\}$, $\{v_n\}$ in corresponding Hilbert spaces, he proved the relation \eqref{relation1}, which is equivalent to the generalized Grunsky
equality with $\lambda_0=0$. Here we use the opposite approach. We
start from the generalized Grunsky equality for the pair 
$(\f^{\mu},\g_{\mu})$ for $[\mu]\in T_{0}(1)$ and use it for deriving all
necessary properties of the Grunsky operators.
In particular, we prove that the Grunsky operators $B_{1}$ and $B_{4}$ associated with $[\mu]\in T_{0}(1)$ are Hilbert-Schmidt.  The case we consider is more general than in \cite{Schiffer2} since the set of all quasi-circles $\f^{\mu}(S^{1})$ for $[\mu]\in T_{0}(1)$ contains the set of all $C^3$ curves as a proper subset. In fact, we prove in Appendix B that the Grunsky operators $B_{1}$ and $B_{4}$ associated with $[\mu]\in T(1)$ are compact if and only if $[\mu]\in S$, the subgroup of symmetric homeomorphisms of $S^{1}$. Our analysis of the relation between singular values of Grunsky operators and Fredholm eigenvalues still holds for this case. 
\end{remark}

As in \cite{Schiffer-multiple}, for a pair $(f,g)$ such that the corresponding operators $\KKK_{1}$ and $\KKK_{4}$ are of trace class, we define the Fredholm determinant for the corresponding quasi-circle 
$\mathcal{C}=f(S^{1})$ by
\begin{align*}
\Det_{F}(\mathcal{C}) = \prod_{n=1}^{\infty} \rho_n^2= \det(I-\KKK_{1})=\det(I-\KKK_{4}).
\end{align*}
Theorem \ref{traceclass} justifies the following definition.
\begin{definition}
The real-valued function $\SSS_2: T_{0}(1)\rightarrow \R$ is
defined as
\begin{align*}
\SSS_{2}([\mu])=\log\Det_F(\f^{\mu}(S^1)),\quad [\mu]\in T_{0}(1).
\end{align*}
\end{definition}
It follows from \eqref{symmetry1-4} that 
\begin{equation} \label{symmetry-2}
\SSS_2([\mu])=\SSS_2([\mu]^{-1}),\quad [\mu]\in T_{0}(1).
\end{equation}
\subsection{Period matrix of $1$-forms}
For a normalized disjoint pair $(f,g)$ of univalent functions we set
$\Omega=f(\Del), \,\Omega^*=g(\Del^{*})$, and define the Hilbert spaces
\begin{align*}
A^{1}_{2}(\Omega) & = \left\{\psi \;\text{holomorphic on}\,\, \Omega:
\, \|\psi\|_{2}^2 = \iint\limits_{ \Omega}
\left|\psi(z)\right|^2d^2z < \infty \right\},\\
A^{1}_{2}(\Omega^*) & = \left\{\psi \;\text{holomorphic on}\,\,
\Omega^*: \, \|\psi\|_{2}^2 = \iint\limits_{ \Omega^*}
\left|\psi(z)\right|^2d^2z < \infty \right\}.
\end{align*}
The Hilbert spaces  $A^{1}_{2}(\Omega)$ and $A^{1}_{2}(\Omega^{*})$ --- the Hilbert spaces of holomorphic $1$-forms on corresponding domains, are, respectively, naturally isomorphic to the Hilbert spaces $A^{1}_{2}(\Del)$ and $A^{1}_{2}(\Del^*)$.

Consider generalized Faber polynomials of $g$ and $f$ defined, respectively, by \cite{Pom, Teo2}
\begin{align*}
\log \frac{g(z) -w}{bz} &= -\sum_{n=1}^{\infty} \frac{P_n(w)}{n}
z^{-n} ,\\
\log \frac{w- f(z)}{w} &= \log \frac{f(z)}{ z}-\sum_{n=1}^{\infty}
\frac{Q_n(w)}{n} z^{n}.
\end{align*}
Here $P_n(w)$ is a polynomial of degree $n$ in $w$ and $Q_n(w)$ is a polynomial of degree $n$ in $1/w$. Specifically,
$$P_n(w) = (g^{-1}(w))^n_{\geq 0},$$
the polynomial part of the $n$-th power of the inverse function $g^{-1}$, and 
$$Q_n(w) = (f^{-1}(w))^{-n}_{\leq 0},$$ 
the principal part of the negative $n$-th power of the inverse function $f^{-1}$. Here  for $S\subset\ZZ$ and a formal power series $A(w)=\sum_{n\in\ZZ} A_nw^n$ 
we denote $(A(w))_{S} = \sum_{n\in S} A_n w^n$.

Comparing the definition of Faber polynomials with the definition of Grunsky coefficients, we obtain the following relations (see, e.g. \cite{Pom, Teo2})
\begin{alignat*}{3}
P_n (g(z)) &= z^n + n\sum_{m=1}^{\infty} b_{nm}
z^{-m}, &\quad P_n (f(z)) & = nb_{n,0}+n\sum_{m=1}^{\infty}
b_{n, -m} z^m, \nonumber\\
Q_n(g(z)) &= -nb_{-n,0} + n\sum_{m=1}^{\infty} b_{m,-n}
z^{-m}, &\quad Q_n (f(z)) &= z^{-n} + n\sum_{m=1}^{\infty}
b_{-n,-m} z^{m}.\nonumber
\end{alignat*}
Now assume that the pair $(f,g)$ is such that the corresponding set $F=\CC\setminus\{f(\Del)\cup g(\Del^{*})\}$ has Lebesgue measure zero. Then it follows from the above formulas and Remark \ref{classicalgrunsky} that the Hilbert spaces $A^{1}_{2}(\Omega)$ and $A^{1}_{2}(\Omega^{*})$ have natural bases $\{\alpha_{n}\}_{n=1}^{\infty}$ and $\{\beta_{n}\}_{n=1}^{\infty}$, 
given respectively by the polynomials 
$$\alpha_n(z) =\frac{P_n^{\prime}(z)}{\sqrt{\pi n}},\quad n\in\mathbb{N},$$ 
and by the Laurent
polynomials 
$$\beta_n(z)=\frac{Q_n^{\prime}(z)}{\sqrt{\pi n}},\quad n\in\mathbb{N}.$$
Indeed, we have
\begin{align*}
\alpha_{n}\circ f f^{\prime} & =\sum_{m=1}^{\infty}(B_{3})_{nm}e_{m}\quad
\text{and}\quad
\beta_{n}\circ g\,g^{\prime} =\sum_{m=1}^{\infty}(B_{2})_{nm}f_{m},
\end{align*}
and the inner products are given by
\begin{align*}
\langle \alpha_n, \alpha_m\rangle  & = \iint\limits_{\Omega} \alpha_n(z) \ov{\alpha_m(z)}d^2z =  \iint\limits_{\Del} \alpha_n(f(z))f^{\prime}(z) \ov{\alpha_m(f(z))f^{\prime}(z)}d^2z \\
& =\sum_{k=1}^{\infty}(B_{3})_{nk}\ov{(B_{3})}_{mk}.
\end{align*}
Hence the period matrix of $A^{1}_{2}(\Omega)$ with respect to the
basis $\{\alpha_n\}_{n=1}^{\infty}$ of holomorphic $1$-forms on $\Omega$ (the Gram matrix of the basis) is given by
\begin{align*}
N_{\Omega}  = \left\{\la \alpha_n, \alpha_m\ra\right\}_{m,n=1}^{\infty}=B_3B_3^*.
\end{align*}
Similarly, the period matrix of the basis $\{\beta_n\}_{n=1}^{\infty}$ of 
holomorphic $1$-forms on $\Omega^{*}$ is given by
\begin{align*}
N_{\Omega^*}  = \{ \la \beta_n, \beta_m\ra\}_{m,n=1}^{\infty}=B_2B_2^*.
\end{align*}
We just proved the following result.
\begin{corollary} \label{det=period} Let $(f,g)$ be a normalized disjoint pair of univalent fuctions such that the set $F=\CC\setminus\{f(\Del)\cup g(\Del^{*})\}$ has Lebesgue measure zero and the corresponding Grunsky operators $B_{1}$ and $B_{4}$ are Hilbert-Schmidt. Then for $\mathcal{C}=f(S^{1})$,
\begin{align*}
\Det_F(\mathcal{C}) = \det N_{\Omega} =\det N_{\Omega^*}
\end{align*}
\end{corollary}
\section{Variations of the functions $\SSS_1$ and $\SSS_2$}
Let $\partial$ and $\bar{\partial}$ be $(1,0)$ and $(0,1)$ components of de Rham differential $d$ on the complex manifold $T_{0}(1)$. Here we compute the ``first variations'' of the functions $\SSS_{1}$ and $\SSS_{2}$ --- the $(1,0)$-forms $\partial\SSS_{1}$ and $\partial\SSS_{2}$ on $T_{0}(1)$.
\subsection{The first variation of $\SSS_2$}
\begin{theorem}\label{varyS2}
The real-valued function
$\SSS_2:T_{0}(1)\rightarrow \R$ is differentiable at every point $[\nu]\in T_{0}(1)$. In terms of the Bers coordinates $\vep_{\mu}$ on the chart $V_{\nu}$,
\begin{displaymath}
\frac{\partial\SSS_{2}}{\partial\vep_{\mu}}([\nu])=
-\frac{1}{6\pi}
\iint\limits_{\Del^*}\mathcal{S}(\g_{\nu})(z)\mu(z)d^2z.
\end{displaymath}
Here $w_{\nu}=\g_{\nu}^{-1}\circ \f^{\nu}$ is the conformal welding corresponding to $[\nu]\in T_{0}(1)$.
\end{theorem}
\begin{proof} By definition of the Bers coordinates (see Section 3.3. in Part I), for $\mu\in H^{-1,1}(\Del^{*})$ 
$$\frac{\pa\SSS_2}{\pa \vep_{\mu}}([\nu])=\left.\frac{d}{d\vep}\right|_{\vep=0}\SSS_{2}([\vep\mu\ast\nu]).$$
Set
$w_{\vep\mu}\circ w_{\nu}=\g_{\vep}^{-1}\circ \f^{\vep}$, 
$\f=\f^{0}=\f^{\nu}$, $\g=\g_{0}=\g_{\nu}$ and $\KKK_{1}(\vep)=
\KKK_{1}(\f^{\vep})$. Since $\KKK_{1}(\vep)$ is a holomorphic family, we have
\begin{align} \label{var-determinant}
\frac{\pa S_2}{\pa \vep_{\mu}}([\nu]) =\left.\frac{\pa}{\pa\vep}\right|_{\vep=0}\det(I-\KKK_{1}(\vep))=
 -\Tr
\left((I-\KKK_1)^{-1}\frac{\pa \KKK_1}{\pa
\vep}(0)\right)
\end{align}
(see, e.g., \cite[Ch. IV.1, Property 9]{gohberg-krein}).
Now using Lemma \ref{varyK1}, we have
\begin{align*}
\frac{\pa \KKK_1}{\pa \vep_{\mu}}([\nu])(z,w)=&
\frac{1}{\pi^2}\iint\limits_{\Del}\iint\limits_{\Del^*}\frac{\mu(u)\f'(z)\g'(u)^2
\f'(\zeta)}{(\f(z)-\g(u))^2(\g(u)-\f(\zeta))^2}K_1^*(\zeta,
w)d^2ud^2\zeta\\
=&\iint\limits_{\Del}\iint\limits_{\Del^*}\mu(u)K_2(z,u)K_3(u,\zeta)K_1^*(\zeta,
w)d^2ud^2\zeta\\
=&-\iint\limits_{\Del^*}\iint\limits_{\Del^*}\mu(u)K_2(z,u)K_4(u,\zeta)K_2^*(\zeta,
w)d^2ud^2\zeta.
\end{align*}
Here in the last line, we have used the second relation in \eqref{relation1},
$$K_3K_1^*=-K_4K_2^*.$$
Let $R_{2}(z,w)$ be the kernel of the
inverse operator $K_2^{-1}$ --- the anti-holomorphic function on
$\Del^*\times \Del$ satisfying
\begin{align*}
\iint\limits_{\Del^*}K_2(z,\zeta)R_{2}(\zeta,
w)d^2\zeta=I_1(z,w)=\frac{1}{\pi(1-z\bar{w})^2},\\
\iint\limits_{\Del}R_{2}(z,\zeta)K_2(\zeta,
w)d^2\zeta=I_2(z,w)=\frac{1}{\pi(1-\z w)^2}.
\end{align*}
Here $I_1(z,w)$ and $I_2(z,w)$ are the kernels of the identity operators on $A_2^1(\Del)$
and $\ov{A_2^1(\Del^*)}$ respectively. Similarly, let $R_{2}^{*}(z,w)$
be the kernel of the inverse operator $(K_{2}^{*})^{-1}$. We have
$$(I-\KKK_1)^{-1}=\KKK_2^{-1}=(K_2^*)^{-1}K_2^{-1},$$
so that
\begin{align*}
\frac{\pa S_2}{\pa
\vep_{\mu}}([\nu])=&\iint\limits_{\Del}\iint\limits_{\Del^*}\iint\limits_{\Del}
\iint\limits_{\Del^*}\iint\limits_{\Del^*}\mu(u)R_{2}^{*}(w,\eta)R_{2}(\eta,z)\\
&\hspace{3cm} K_2(z,u)K_4(u,\zeta)K_2^*(\zeta, w)d^2ud^2\zeta
d^2zd^2\eta d^2w\\
=&\iint\limits_{\Del^*}\iint\limits_{\Del^*}\iint\limits_{\Del^*}\mu(u)K_4(u,\zeta)I_2(\eta,u)I_2(\zeta,\eta)d^2u
d^2\eta d^2\zeta\\
=&\iint\limits_{\Del^*}\iint\limits_{\Del^*}\mu(u)K_4(u,\zeta)I_2(
\zeta, u)d^2ud^2\zeta\\
=&\iint\limits_{\Del^*}\mu(u)K_4(u,u)d^2u=-\frac{1}{6\pi}\iint\limits_{\Del^{*}}\mathcal{S}(\g_{\nu})(u)\mu(u)d^2u.
\end{align*}
Here in the last line we have used
\begin{align*}
K_4(u,u)=-\frac{1}{\pi}\lim_{\zeta\rightarrow
u}\left(\frac{\g'(u)\g'(\zeta)}{(\g(\zeta)-\g(u))^2}-\frac{1}{(\zeta-u)^2}\right)
=-\frac{1}{6\pi}\mathcal{S}(\g)(u).
\end{align*}
\end{proof}
Denote by $T^{*}_{[\mu]}T_{0}(1)$ the holomorphic cotangent space to $T_{0}(1)$ at a point $[\mu]\in T_{0}(1)$. The natural isomorphism
$T_{[\mu]}T_{0}(1)\simeq H^{-1,1}(\Del^{*})$ induces the isomorphism
$T^{*}_{[\mu]}T_{0}(1)\simeq A_{2}(\Del^{*})$. Define a holomorphic $1$-form $\boldsymbol{\vartheta}$ on $T_{0}(1)$ by
\begin{align*}
\boldsymbol{\vartheta}_{[\mu]} = \mathcal{S}(\g_{\mu})\in A_{2}(\Del^{\ast}),
\end{align*}
where $w_{\mu}=g_{\mu}^{-1}\circ \f^{\mu}\in T_{0}(1)$.
\begin{corollary} \label{1-form1} On $T_{0}(1)$,
$$\partial\SSS_{2} =-\frac{1}{6\pi}\boldsymbol\vartheta.$$
\end{corollary}
\begin{remark}
For $C^{3}$ curves the statement of Theorem \ref{varyS2} was obtained by Schiffer in
\cite{Schiffer-multiple}. The derivation in \cite{Schiffer-multiple} uses the variational theory of Fredholm eigenvalues and the exterior variation of the domain. Our proof is different from Schiffer's: we use general formula \eqref{var-determinant} 
and the quasi-conformal variation.
\end{remark}
\subsection{The first variation of $\SSS_{1}$}
In addition to $\SSS_{1}$, we introduce another function
$\tilde\SSS_{1}: T_{0}(1)\rightarrow \RR$ defined by
$$\tilde\SSS_{1}([\mu]) = \SSS_{1}([\mu^{-1}]).$$ Using \eqref{inversion}, we get
\begin{align*}
\tilde\SSS_{1}([\mu])& =
\iint\limits_{\Del}\left|\mathcal{A}(\f^{\mu})-2\frac{(\f^{\mu})'}{\f^{\mu}}+
\frac{2}{z}\right|^2d^2z
+\iint\limits_{\Del^*}
\left|\mathcal{A}(\g_{\mu})-2\frac{\g_{\mu}'}{\g_{\mu}}+\frac{2}{z}\right|^2\
d^2z \\ &-4\pi\log |(\g_{\mu})'(\infty)| \\
&=\iint\limits_{\Del}
\left|\mathcal{A}(\tilde{\g}_{\mu})\right|^2 d^2z+ \iint\limits_{\Del^*}\left|\mathcal{A}(\tilde{\f}^{\mu})\right|^2d^2z+4\pi\log|\tilde\g_{\mu}'(0)|,
\end{align*}
where $\tilde\f^{\mu}=\imath\circ\f^{\mu}\circ \imath$, $\tilde{\g}_{\mu}=\imath\circ \g_{\mu} \circ\imath$ and $\imath(z)=\tfrac{1}{z}$. The functions $\tilde\f^{\mu}$ and $\tilde\g_{\mu}$ are univalent, respectively, on the domains $\Del^*$ and $\Del$ and are normalized as $\tilde\f^{\mu}(\infty)=\infty,\,
(\tilde\f^{\mu})'(\infty)=1$ and $\tilde\g_{\mu}(0)=0$. They satisfy the factorization
\begin{equation} \label{tilde-welding}
\tilde{w}_{\mu}=\tilde\g_{\mu}^{-1}\circ\tilde\f^{\mu},
\end{equation}
where $\tilde{w}_{\mu}=\imath\circ w_{\mu}\circ \imath$. 

This description corresponds to the realization of $T(1)$ associated
with the model $\HH^{2}\simeq\Del$. Namely, due to the canonical isomorphism
\begin{align*}
\mu \in L^{\infty}(\Del^*) \mapsto \tilde{\mu}=
\imath^*\mu=\mu\left(\frac{1}{z}\right) \frac{z^2}{\z^2}\in
L^{\infty}(\Del),
\end{align*}
we have  $T(1) \simeq L^{\infty}(\Del)_1/\sim$. If $w_{\mu}$ is a q.c.~mapping associated with
$\mu \in L^{\infty}(\Del^{*})_1$, then $\tilde{w}_{\mu} $ is
the q.c.~ mapping associated with $\tilde\mu\in L^{\infty}(\Del)_{1}$, and corresponding conformal wielding is given by \eqref{tilde-welding}.

In this section, we will also use the model $T(1)\simeq L^{\infty}(\Del)_{1}$. To simplify the notations, for 
$\mu\in L^{\infty}(\Del)_{1}$ we will denote corresponding q.c.~mapping
by $w_{\mu}=\g_{\mu}^{-1}\circ\, \f^{\mu}$, where $\f^{\mu}$ and $\g_{\mu}$ are univalent on the domains $\Del^*$ and $\Del$ and are normalized as $\f^{\mu}(\infty)=\infty,\,(\f^{\mu})'(\infty)=1$ and $\g_{\mu}(0)=0$. Correspondingly, for $\gamma=g^{-1}\circ f\in\mathcal{T}(1)$ we would have the normalization $f(\infty)=\infty,\,f^{\prime}(\infty)=1$ and $g(0)=0$. To avoid confusion with the notations for our primary model $T(1)=L^{\infty}(\Del^{*})_{1}/\sim$,
we will always specify explicitly in the main text when we are using the model
$T(1)\simeq L^{\infty}(\Del)_{1}/\sim$.

The function $\SSS_{1}$ on $T_{0}(1)$ naturally extends to a function
$\hat\SSS$ on $\mathcal{T}_{0}(1)$, defined by
$$\hat\SSS(\gamma) = \iint\limits_{\Del}\left|\mathcal{A}(f)\right|^2d^2z+\iint\limits_{\Del^{*}}
\left|\mathcal{A}(g)\right|^2 d^2z-4\pi\log|g'(\infty)|,$$
where $\gamma=g^{-1}\circ f\in\mathcal{T}_{0}(1)$. 
For $\tilde\SSS(\gamma)=\hat\SSS(\gamma^{-1})$ we have
\begin{align*}
\tilde\SSS(\gamma)=
\iint\limits_{\Del}
\left|\mathcal{A}(\tilde{g})\right|^2 d^2z+ \iint\limits_{\Del^*}\left|\mathcal{A}(\tilde{f})\right|^2d^2z+4\pi\log|\tilde{g}'(0)|,
\end{align*}
where $\tilde{f}=\imath\circ f\circ \imath$ and $\tilde{g}=\imath\circ g \circ\imath$.
\begin{lemma} The function $\tilde{\SSS}$ is constant along the fibers of the canonical projection $\pi: \mathcal{T}_{0}(1)\rightarrow T_{0}(1)$, $\tilde{\SSS}=\tilde{\SSS}_{1}\circ \pi$.
 \end{lemma}
\begin{proof} We are using the model $T(1)\simeq L^{\infty}(\Del)_{1}/\sim$. For $\mu\in L^{\infty}(\Del)_{1}$ let $\gamma=g^{-1}\circ f,\, \gamma_{\mu}=g_{\mu}^{-1}\circ f^{\mu} \in\mathcal{T}_{0}(1)$ be such that $\pi(\gamma)=\pi(\gamma_{\mu})=[\mu]$.
Comparing the normalization for $f$ and $f^{\mu}$ at $\infty$, we get
\begin{align*}
f=\sigma\circ f^{\mu}\quad\text{and}\quad g= \sigma\circ g_{\mu}\circ
\alpha^{-1},
\end{align*}
 for some $\alpha\in \PSU(1,1)$ and $\sigma(z)=z+b_0$. Since $f\mapsto\mathcal{A}(f)$ is invariant if $f$ is
post-composed with a translation\footnote{This is why it is more convenient to use the model $T(1)\simeq L^{\infty}(\Del)_{1}/\sim$.}, to 
prove that $\tilde\SSS(\gamma)=\tilde\SSS(\gamma_{\mu})$ we need only to check that for
$\alpha\in \PSU(1,1)$,
\begin{align*}
\iint\limits_{\Del}|\mathcal{A}(g\circ
\alpha^{-1})|^2d^2z+4\pi\log |(g\circ
\alpha^{-1})'(0)|=\iint\limits_{\Del}|\mathcal{A}(g)|^2d^2z+4\pi\log
|g'(0)|.
\end{align*}Let
\begin{align*}
\alpha(z)=e^{i\theta} \frac{z-w}{1-z\bar{w}}
\end{align*}
and set $\log g'(z)=\sum_{n=0}^{\infty}a_nz^{n}$. Then
$\mathcal{A}(g)=\sum_{n=1}^{\infty}na_nz^{n-1}$ and
\begin{gather*}
\iint\limits_{\Del}|\mathcal{A}(g\circ
\alpha^{-1})|^2d^2z=\iint\limits_{\Del}|\mathcal{A}(g)\circ
\alpha^{-1}(\alpha^{-1})'+ \mathcal{A}(\alpha^{-1})|^2d^2z\\
=\iint\limits_{\Del}|\mathcal{A}(g)-\mathcal{A}(\alpha)|^2d^2z
=\iint\limits_{\Del}\left|\mathcal{A}(g)-\frac{2\bar{w}}{1-z\bar{w}}\right|^2d^2z\\
=\iint\limits_{\Del}|\mathcal{A}(g)|^2d^2z
-4\re\left(w\iint\limits_{\Del}\mathcal{A}(g)(z)\sum_{n=1}^{\infty}
(w\z)^{n-1} d^2z\right) \\ + 4|w|^2
\iint\limits_{\Del}\left|\sum_{n=1}^{\infty}
(w\z)^{n-1}\right|^2d^2z.
\end{gather*}
The last two terms give
\begin{gather*}
-4\pi\re\left( \sum_{n=1}^{\infty} a_n w^n \right)+
4\pi\sum_{n=1}^{\infty} \frac{|w|^{2n}}{n} \\
=-4\pi \log |g'(w)|+4\pi\log|g'(0)| -4\pi\log (1-|w|^2).
\end{gather*}
On the other hand, we have
\begin{align*}
(g\circ \alpha^{-1})'(0)=g'(\alpha^{-1}(0))(\alpha^{-1})'(0)=
(1-|w|^2)g'(w).
\end{align*}
This concludes the proof.
\end{proof}

\begin{theorem}\label{varyS3}
The real-valued function
$\tilde\SSS_1:T_{0}(1)\rightarrow \R$ is differentiable at every point $[\nu]\in T_{0}(1)$.
In terms of the Bers coordinates $\vep_{\mu}$ on the chart $V_{\nu}$,
\begin{displaymath}
\frac{\partial\tilde\SSS_{1}}{\partial\vep_{\mu}}([\nu]) = 2
\iint\limits_{\Del^*}\mathcal{S}(\g_{\nu})(z)\mu(z)d^2z.
\end{displaymath}
\end{theorem}
\begin{proof} We are using the model $T(1)\simeq L^{\infty}(\Del)_{1}$.
For $[\nu]\in T_{0}(1)$ choose a representative $\nu\in L^{\infty}(\Del)_{1}$ which is a product of elements in $H^{-1,1}(\Del)_{1}$, and let $w_{\vep}=w_{\vep\mu}\circ w_{\nu}=\g_{\vep}^{-1}\circ \f^{\vep}$.
It follows from Lemma 2.5 in Part I that corresponding $\gamma_{\vep}=g_{\vep}^{-1}\circ f^{\vep}$ fixes $0, 1, \infty$.  
By the above lemma,
$$\tilde\SSS_{1}([\vep\mu\ast\nu])=\tilde\SSS(\gamma_{\vep}).$$
We have
$\gamma_{\vep}\circ\gamma^{-1}_{\nu}=\gamma_{\vep\kappa}$,
where $\kappa = (\alpha^{-1})^*(\mu)$ and $\alpha=\gamma_{\nu}\circ w_{\nu}^{-1}\in \PSU(1,1)$. Set $f=f^{0}$, $g=g^{0}$, so that $g=\sigma\circ g_{\nu}\circ \alpha^{-1}$ for some $\sigma\in\PSL(2,\CC)$, and define
$v_{\vep}=f^{\vep}\circ f^{-1}$.
Since $f^{\vep}$ is normalized, it's Laurent expansion at $\infty$ has the form
\begin{align*}
f^{\vep}(z) = z + \frac{b_1}{z} +\frac{b_2}{z^2}+\ldots.
\end{align*}
Hence  
$$\left.\frac{\pa}{\pa\vep}\right\vert_{\vep=0}v_{\vep}(z) = O(z^{-1})\quad\text{as}\quad z\rightarrow\infty,$$
and the first variations of $v_{\vep}$ have the form
\begin{align*}
\left.\frac{\pa}{\pa\vep}\right\vert_{\vep=0}
v_{\vep}(z)=-\frac{1}{\pi}\iint\limits_{\Omega}\frac{((g^{-1})^*\kappa)(w)}{w-z}d^2w,\quad
\left.\frac{\pa}{\pa\bar{\vep}}\right\vert_{\vep=0} v_{\vep}(z)=0,
\end{align*}
where $\Omega = g(\Del)$. Since
$\gamma_{\vep\kappa}$ fixes $0,1,\infty$, we also have
\begin{align*}
\left.\frac{\pa}{\pa\vep}\right\vert_{\vep=0} \gamma_{\vep\kappa} (z)=&
-\frac{1}{\pi}\iint\limits_{\Del} \frac{z(z-1)\kappa(w)}{(w-z)w(w-1)}d^2w,\\
 \left.\frac{\pa}{\pa\bar{\vep}}\right\vert_{\vep=0}
\gamma_{\vep\kappa} (z)=& -\frac{1}{\pi}\iint\limits_{\Del}
\frac{z(z-1)\ov{\kappa(w)}}{(1-\bar{w}z)\bar{w}(1-\bar{w})}d^2w.
\end{align*}

Using $f^{\vep}=v_{\vep}\circ f$, we obtain
\begin{align*}
\mathcal{A}(f^{\vep}) =\mathcal{A}(v_{\vep})\circ f f' +
\mathcal{A}(f).
\end{align*} 
Applying the variational formulas for $v_{\vep}$, we have
\begin{align*}
\left.\frac{\pa}{\pa\vep}\right\vert_{\vep=0}\mathcal{A}(v_{\vep})(z)=\left.\frac{\pa^2}{\pa
z^2}\frac{\pa}{\pa\vep}\right\vert_{\vep=0}
v_{\vep}(z)=-\frac{2}{\pi}\iint\limits_{\Omega}\frac{((g^{-1})^*\kappa)(w)}{(w-z)^3}d^2w,
\end{align*}
and hence
\begin{align*}
\left.\frac{\pa}{\pa
\vep}\right\vert_{\vep=0}\iint\limits_{\Del^*}\left|\mathcal{A}(f^{\vep})\right|^2
d^2z & = -\frac{2}{\pi}\iint\limits_{\Del^*}\iint\limits_{\Del}
\frac{\kappa(w)g'(w)^2f'(z)}{(g(w)-f(z))^3}\ov{\mathcal{A}(f)(z)}d^2wd^2z =I_1.
\end{align*}
Similarly, using
\begin{align*}
g_{\vep}\circ \gamma_{\vep\kappa}=v_{\vep}\circ g,
\end{align*}
we have
\begin{align*}
g_{\vep}'\circ \gamma_{\vep\kappa}(\gamma_{\vep\kappa})_z &=v_{\vep}'\circ g\,g',\\
\intertext{and}
\mathcal{A}(g_{\vep})\circ
\gamma_{\vep\kappa}(\gamma_{\vep\kappa})_z+\mathcal{A}(\gamma_{\vep\kappa})&=\mathcal{A}(v_{\vep})\circ
g\,g' + \mathcal{A}(g),
\end{align*}
where
$$\mathcal{A}(\gamma_{\vep\kappa})=\frac{(\gamma_{\vep\kappa})_{zz}}{(\gamma_{\vep\kappa})_{z}}.$$
Hence we have
\begin{align*}
\left.\frac{\pa}{\pa\vep}\right|_{\vep=0}g_{\vep}'(0) & =-
g''(0)\left(\left.\frac{\pa}{\pa\vep}\right\vert_{\vep=0}\gamma_{\vep\kappa}\right)(0)-g'(0)
\frac{\pa}{\pa z}\left(\left.\frac{\pa}{\pa\vep}\right\vert_{\vep=0}
\gamma_{\vep\kappa}\right)(0) \\
&\,\,\,\;\;+g'(0)\frac{\pa}{\pa
z}\left(\left.\frac{\pa}{\pa\vep}\right|_{\vep=0}v_{\vep}\right)(0)\\
& =\frac{g'(0)}{\pi}
\iint\limits_{\Del}\kappa(w)\left(\frac{1}{w^2}
-\frac{1}{w(w-1)} - \frac{g'(w)^{2}}{g(w)^{2}}\right)d^2w, \\
\intertext{and}
\left.\frac{\pa}{\pa\vep}\right|_{\vep=0}\ov{g_{\vep}'(0)} & =-
\ov{g''(0)\left(\left.\frac{\pa}{\pa\bar{\vep}}\right|_{\vep=0}\gamma_{\vep\kappa}\right)(0)}-\ov{g'(0)
\frac{\pa}{\pa
z}\left(\left.\frac{\pa}{\pa\bar{\vep}}\right|_{\vep=0}\gamma_{\vep\kappa}\right)(0)}\\& =
\frac{\ov{g'(0)}}{\pi}\iint\limits_{\Del}
\frac{\kappa(w)}{w(w-1)}d^2w,
\end{align*}
as well as
\begin{align*}
\frac{\pa}{\pa\vep}\Bigr|_{\vep=0}\mathcal{A}(g_{\vep})\circ\gamma_{\vep\kappa}(\gamma_{\vep\kappa})_z
& =\left(\left.\frac{\pa}{\pa
\vep}\right|_{\vep=0}\mathcal{A}(v_{\vep})\right)\circ g
g'-\left.\frac{\pa}{\pa
\vep}\right|_{\vep=0}\mathcal{A}(\gamma_{\vep\kappa})\\
&=-\frac{2}{\pi}\iint\limits_{\Del}\kappa(w)\left(\frac{g'(w)^2g'(z)}{(g(w)-g(z))^3}-\frac{1}{(w-z)^3}\right)
d^2w,\\
\intertext{and}
\left.\frac{\pa}{\pa\vep}\right|_{\vep=0}\ov{\mathcal{A}(g_{\vep})\circ\gamma_{\vep\kappa}
(\gamma_{\vep\kappa})_z}  & =-\ov{\left.\frac{\pa}{\pa
\bar{\vep}}\right|_{\vep=0}\mathcal{A}(\gamma_{\vep\kappa})}
=\frac{2}{\pi}\iint\limits_{\Del}\frac{\kappa(w)}{(1-w\z)^3
w}d^2w.
\end{align*}
From here we get
\begin{align*}
2\pi\left.\frac{\pa}{\pa\vep}\right |_{\vep=0}\log|g_{\vep}'(0)|^2 & =
-2 \iint\limits_{\Del}
\kappa(w)\left(\frac{g'(w)^2}{g(w)^2}-\frac{1}{w^2}\right)d^2w=I_2,
\end{align*}
and
\begin{gather*}
\left.\frac{\pa}{\pa\vep}\right|_{\vep=0}\iint\limits_{\Del}\left|\mathcal{A}(g_{\vep})\right|^2
d^2z =\left.\frac{\pa}{\pa\vep}\right|_{\vep=0}\iint\limits_{\Del}\left|\mathcal{A}(g_{\vep})\circ
\gamma_{\vep\kappa}
(\gamma_{\vep\kappa})_z\right|^2(1-|\vep\kappa|^2)
d^2z\\
 =-\frac{2}{\pi}\iint\limits_{\Del}\iint\limits_{\Del}
\kappa(w)\left(\frac{g'(w)^2g'(z)}{(g(w)-g(z))^3}-\frac{1}{(w-z)^3}\right)
\ov{\mathcal{A}(g)(z)}d^2wd^2z \\
+\frac{2}{\pi}\iint\limits_{\Del}\iint\limits_{
\Del}\frac{\kappa(w)\mathcal{A}(g)(z)}{w(1-w\z)^3}d^2wd^2z=I_3+I_4.
\end{gather*}
Let $\log g'(z) = \sum_{n=0}^{\infty} a_n z^n$ be the power series
expansion of $\log g'(z)$. Then $\mathcal{A}(g) =
\sum_{n=1}^{\infty} na_nz^{n-1}$. Explicit computation gives
\begin{align*}
\frac{2}{\pi}\iint\limits_{\Del}\frac{\mathcal{A}(g)(z)}{w(1-w\z)^3}d^{2}z =\sum_{n=1}^{\infty}n(n+1)a_n w^{n-2}
=\mathcal{A}(g)'(w)+\frac{2}{w}\mathcal{A}(g)(w).
\end{align*}
Hence
\begin{align*}
I_4=
\iint\limits_{\Del}\kappa(w)\left( \mathcal{A}(g)'(w)+\frac{2}{w}\mathcal{A}(g)(w)\right)d^2w.
\end{align*}
To compute the other terms,
 we define the following holomorphic function on $\Del$,
\begin{align*}
h(w)& =\frac{1}{\pi}
\iint\limits_{\Del^*}\frac{g'(w)f'(z)}{(g(w)-f(z))^2}
\ov{\mathcal{A}(f)(z)}d^2z\\
& +\frac{1}{\pi}
\iint\limits_{\Del}\left(\frac{g'(w)g'(z)}{(g(w)-g(z))^2}-\frac{1}{(w-z)^2}\right)
\ov{\mathcal{A}(g)(z)}d^2z.\nonumber
\end{align*}
Then it is easy to check that
\begin{gather*}
\left.\frac{\pa}{\pa\vep}\right|_{\vep=0}\tilde\SSS(\gamma_{\vep})= I_1+I_2+I_3+I_4\\
=\iint\limits_{\Del} \kappa(w)\left(
h'(w)-\mathcal{A}(g)(w)h(w)-2\frac{g'(w)^2}{g(w)^2}+\frac{2}{w^2}
+\mathcal{A}(g)'(w)+\frac{2}{w}\mathcal{A}(g)(w)\right)d^2w.
\end{gather*}
To finish the proof, we claim that
\begin{align*}
h(w) = \mathcal{A}(g)(w) -2\frac{g'(w)}{g(w)}+\frac{2}{w},
\end{align*}
which is going to be proved in the next lemma. With this equation
for $h$, it is straightforward to compute that
\begin{align*}
\left.\frac{\pa}{\pa\vep}\right|_{\vep=0}\tilde\SSS(\gamma_{\vep}) &=\iint\limits_{\Del}
\left(2\mathcal{A}(g)'(w)-\mathcal{A}(g)(w)^2\right)\kappa(w)d^2w \\
&=2\iint\limits_{\Del}\mathcal{S}(g)(w)\kappa(w)d^2w
=2\iint\limits_{\Del}\mathcal{S}(\g_{\nu})(w)\mu(w)d^2w.
\end{align*}
Returning back to the model $T(1) = L^{\infty}(\Del^{*})_{1}/\sim$, we get the statement of the theorem.
\end{proof}
\begin{lemma}\label{important} In the model $T(1)\simeq L^{\infty}(\Del)_1/\sim$, let $\gamma=g^{-1}\circ f$ be the conformal welding corresponding to $\gamma\in\mathcal{T}_{0}(1)$. Then for $z\in\Del$ the following identity holds
\begin{gather*}
\mathcal{A}(g)(z) -2\frac{g'(z)}{g(z)}+\frac{2}{z}=\frac{1}{\pi}
\iint\limits_{\Del^*}\frac{g'(z)f'(w)}{(g(z)-f(w))^2}
\ov{\mathcal{A}(f)(w)}d^2w \\
+\frac{1}{\pi}
\iint\limits_{\Del}\left(\frac{g'(z)g'(w)}{(g(w)-g(w))^2}-\frac{1}{(z-w)^2}\right)
\ov{\mathcal{A}(g)(w)}d^2w. 
\end{gather*}
\end{lemma}
\begin{proof}
First we consider the case when $\mathcal{A}(g)$ and
$\mathcal{A}(f)$ are smooth functions on $S^1$. Specifically, we assume that the Beltrami differential $\mu$ corresponding to $\pi(\gamma)\in T_{0}(1)$, is smooth on $\C$ and $\left.\mu\right|_{S^1}=\left.\mu_{\bar{z}}\right|_{S^{1}}=0$.
Denote by $h(z)$ the right-hand side of the identity of the lemma. Changing the variables of integration and using Stokes' theorem, we obtain
\begin{align*}
h\circ g^{-1}(
g^{-1})'(z)=&-\frac{1}{\pi}\iint\limits_{\Omega^*}\frac{\ov{\mathcal{A}(f^{-1})(w)}}{(z-w)^2}d^2w
-\frac{1}{\pi}\iint\limits_{\Omega}\frac{\ov{\mathcal{A}(g^{-1})(w)}}{(z-w)^2}d^2w\\
=&\frac{1}{2\pi i}\oint_{\mathcal{C}}
\frac{1}{(z-w)}\left(\ov{\mathcal{A}(g^{-1})(w)}-\ov{\mathcal{A}(f^{-1})(w)}\right)d\bar{w},
\end{align*}
where $\Omega=g(\Del), \Omega^*=f(\Del^*)$ and
$\mathcal{C}=g(S^1)$. Next, consider the relation $\tilde\gamma\circ g^{-1}=f^{-1}$, where $\tilde\gamma=\gamma^{-1}$, and differentiate it twice with respect to $z$. Since $\tilde\gamma_{\bar{z}}$ vanishes on $S^{1}$, we get the following relations on $\mathcal{C}$,
\begin{align*}
\frac{\tilde\gamma_{z}}{\tilde\gamma}\circ g^{-1}(g^{-1})_{z} & =\frac{(f^{-1})_{z}}{f^{-1}}, \\
\mathcal{A}(\tilde\gamma)\circ g^{-1} (g^{-1})_{z}
& =\mathcal{A}(f^{-1})  - \mathcal{A}(g^{-1}).
\end{align*}
Hence
\begin{align*}
h\circ g^{-1}( g^{-1})'(z)=&-\frac{1}{2\pi i}\oint_{\mathcal{C}}
\frac{1}{(z-w)}\ov{\mathcal{(A}(\tilde\gamma)\circ
g^{-1})(w) (g^{-1})_{w}(w)}d\bar{w}\\
=&-\frac{1}{2\pi i}
\oint_{S^1}\frac{1}{z-g(w)}\ov{\mathcal{A}(\tilde\gamma)(w)}d\bar{w}.
\end{align*}
On the other hand, since $j\circ \tilde\gamma = \tilde\gamma\circ
j$, where $j$ is the inversion $z\mapsto \tfrac{1}{\z}$, we have
\begin{align*}
\ov{\mathcal{A}(\tilde\gamma)}
=&\mathcal{A}(\tilde\gamma)\circ j\,j_{\z}-2\,
\frac{\tilde\gamma_z}{\tilde\gamma}\circ j\,j_{\z} +\ov{\mathcal{A}(\bar{j})}.
\end{align*}
Hence
\begin{align*}
&h\circ g^{-1}( g^{-1})'(z)=\frac{1}{2\pi
i}\oint_{S^1}\frac{1}{z-g(w)}\left(\mathcal{A}(\tilde\gamma)\left(\frac{1}{\bar{w}}\right)
\frac{1}{\bar{w}^{2}}-2
\frac{\tilde\gamma_w}{\tilde\gamma}\left(\frac{1}{\bar{w}}\right)\frac{1}{\bar{w}^{2}}+\frac{2}{\bar{w}}\right)d\bar{w}\\
=&-\frac{1}{2\pi
i}\oint_{S^1}\frac{1}{z-g(w)}\left(\mathcal{A}(\tilde\gamma)(w)-2
\frac{\tilde\gamma_w(w)}{\tilde\gamma
\left(w\right)}+\frac{2}{w}\right)dw\\
=&\frac{1}{2\pi
i}\oint_{\mathcal{C}}\frac{1}{(z-w)}\left(\mathcal{A}(g^{-1})(w)
-\mathcal{A}(f^{-1})(w)+2\frac{(f^{-1})_{w}(w)}{f^{-1}(w)}
-2\frac{(g^{-1})_{w}(w)}{g^{-1}(w)}\right)dw.
\end{align*}
The functions
\begin{align*}
\mathcal{A}(f^{-1})(z)-2\frac{(f^{-1})_{z}(z)}{f^{-1}(z)}+\frac{2}{z}\hspace{1cm}\text{and}
\hspace{1cm}\mathcal{A}(g^{-1})(z)-2\frac{(g^{-1})_{z}(z)}{g^{-1}(z)}+\frac{2}{z}
\end{align*}
are holomorphic on $\Omega^*=f(\Del^*)$ and $\Omega=g(\Del)$ respectively and due to the normalization of $f$,
$$\mathcal{A}(f^{-1})(z)-2\frac{(f^{-1})_{z}(z)}{f^{-1}(z)}+\frac{2}{z}=O\left(\frac{1}{z^{2}}\right)\quad\text{as}\quad z\rightarrow\infty.$$
Thus we have by Cauchy formula
\begin{align*}
h\circ g^{-1}(
g^{-1})'(z)=-\left(\mathcal{A}(g^{-1})(z)-2\frac{(g^{-1})'(z)}{g^{-1}(z)}+\frac{2}{z}\right)
\end{align*}
or equivalently,
\begin{align*}
h(z)= \mathcal{A}(g)(z) -2\frac{g'(z)}{g(z)}+\frac{2}{z}.
\end{align*}
For a general point $\gamma=g^{-1}\circ f$ in $\mathcal{T}_0(1)$,
we let $f_n =r_n^{-1}\circ f\circ r_n$, where 
$r_n$ is the dilation $z\mapsto \tfrac{n+1}{n} z$. 
Since $f_n$ is a normalized
univalent function on $|z|>\tfrac{n}{n+1}$, 
corresponding $\gamma_n^{-1}=g_n^{-1}\circ f_n\in\mathcal{T}_0(1)$
satisfies the assumptions made in the beginning of the proof. Since 
$\mathcal{A}(f)\in A_{2}^{1}(\Del^{*})$, we see that
\begin{align*}
&\bigl\Vert \mathcal{A}(\iota\circ f_n\circ
\iota)-\mathcal{A}(\iota\circ
f\circ\iota)\bigr\Vert_{A^{1}_{2}(\Del)}\\=&\left\Vert
\left(\mathcal{A}(f_n)-2\frac{f_n'}{f_n}+\frac{2}{z}\right)-\left(
\mathcal{A}(f)-2\frac{f'}{f}+\frac{2}{z}\right)\right\Vert_{A_{2}^{1}(\Del^{*})}\rightarrow
0 \quad\text{as}\quad n\rightarrow \infty. 
\end{align*}
 By
Corollary \ref{convergea} and Corollary \ref{convergeb} in Appendix A we also have
\begin{gather*}\lim_{n\rightarrow\infty}\Vert
\mathcal{A}(g_n)-\mathcal{A}(g)\Vert_{A_{2}^{1}(\Del)}  = 0\\
\intertext{and}
\lim_{n\rightarrow\infty}\left\Vert
\left(\mathcal{A}(g_n)-2\frac{g_n'}{g_n}+\frac{2}{z}\right)-\left(
\mathcal{A}(g)-2\frac{g'}{g}+\frac{2}{z}\right)\right\Vert_{A_{2}^{1}(\Del)}  = 0.
\end{gather*}
In particular, since convergence in $A_2^1(\Del)$ implies
convergence in $A_{\infty}^1(\Del)$, 
we get
\[
\lim_{n\rightarrow\infty}\left(\mathcal{A}(g_n)(z)-2\frac{g_n'}{g_n}(z)+\frac{2}{z}\right)=\mathcal{A}(g)(z)-2\frac{g'}{g}(z)+\frac{2}{z},
\]
uniformly on compact subsets of $\Del$. 
Since we have already shown that
\[
h_n(z)=\mathcal{A}(g_n)(z)-2\frac{g_n'}{g_n}(z)+\frac{2}{z},
\]
to finish the proof of the lemma we need to verify that 
$\lim_{n\rightarrow\infty}h_n(z)=h(z)$ uniformly on compact
subsets  of $\Del$.

We denote by $K_1[n]$ and $K_2[n]$ the operators associated with the disjoint pair of univalent functions $(g_n, f_n)$, and by $K_1$ and $K_2$ --- the operators associated with the pair $(g,f)$. Then
\begin{align*}
h_n(z)-h(z) & =-\Bigr(K_1[n]\ov{\mathcal{A}(g_n)}\Bigr)(z) +\Bigl(K_1\ov{\mathcal{A}(g)}\Bigr)(z) \\
&\;\;\;+ \Bigl(K_2[n]
\ov{\mathcal{A}(f_n)}\Bigr)(z)-\Bigl(K_2
\ov{\mathcal{A}(f)}\Bigr)(z).
\end{align*}
Now using Theorem B.1 from Appendix B, and
the fact that the inverse map is continuous on $T_0(1)$, we get
that
\begin{align*}
\lim_{n\rightarrow\infty}\Vert K_1[n]-K_1\Vert =0,
\end{align*}
where $\Vert ~\Vert$ stands for the norm of the Banach space $\cB(\ov{A_{2}^{1}(\Del)}, A_{2}^{1}(\Del))$.
Since $\Vert K_1[n]\Vert\leq 1$, we have
\begin{gather*}
\left\Vert K_1[n]\ov{\mathcal{A}(g_n)}-
K_1\ov{\mathcal{A}(g)}\right\Vert_{A_{2}^{1}(\Del)} \\
\leq \left\Vert
K_1[n]\bigl(\ov{\mathcal{A}(g_n)-\mathcal{A}(g)}\bigr)\right\Vert_{A_{2}^{1}(\Del)} +
\left\Vert
\bigl(K_1[n]-K_1\bigr)\ov{\mathcal{A}(g)}\right\Vert_{A_{2}^{1}(\Del)}\\
\leq \Vert \mathcal{A}(g_n)- \mathcal{A}(g)\Vert_{A_{2}^{1}(\Del)}+\Vert
K_1[n]-K_1\Vert \Vert \mathcal{A}(g)\Vert_{A_{2}^{1}(\Del)},
\end{gather*}
which tends to $0$ as $n\rightarrow \infty$.  Consequently,
$\lim_{n\rightarrow\infty}(K_1[n]\ov{\mathcal{A}(g_n)})(z)=(K_1\ov{\mathcal{A}(g)})(z)$, uniformly on compact subsets of $\Del$.

To prove the convergence of the other term in $h_{n}(z)-h(z)$, we let $\mathfrak{g}_n = r_n\circ g_n$ and
$\mathfrak{f}_n=r_n\circ f_n=f\circ r_n$. Let
$\Omega_n^*=\mathfrak{f}_n(\Del^*)=f(\{|z|> \tfrac{n+1}{n}\})$. Since 
$\Omega_{n}^{*}\subseteq
\Omega_{n+1}^*$, the sequence of domains $\Omega_n=\mathfrak{g}_n(\Del)$
is a decreasing sequence that contains $0$ and $\bigcap
\mathfrak{g}_n(\Del)=g(\Del)=\Omega$. By Caratheodory kernel
theorem (see, e.g., \cite{Pom}), the sequence of univalent
functions $\mathfrak{g}_n:\Del\rightarrow \C$ converges uniformly
on compact sets to the univalent function $g:\Del\rightarrow \C$.
By Weierstrass theorem, $\lim_{n\rightarrow\infty}\mathfrak{g}_n'(z)=g'(z)$, uniformly on compact subsets of $\Del$. Using that the operator $K_2$ is unaffected by a simultaneous post-composition of $f$ and $g$ with $\alpha\in\PSL(2,\CC)$ and that
$\mathcal{A}(\mathfrak{f}_n)=\mathcal{A}(r_n\circ
f_n)=\mathcal{A}(f_n)$, we have
\begin{gather*}
\Bigl(K_2[n] \ov{\mathcal{A}(f_n)}\Bigr)(z)-\Bigl(K_2
\ov{\mathcal{A}(f)}\Bigr)(z)\\=
\frac{1}{\pi}\iint\limits_{\Del^*}
\frac{\mathfrak{g}_n'(z)\mathfrak{f}_n'(w)}{(\mathfrak{g}_n(z)-\mathfrak{f}_n(w))^2}
\ov{\mathcal{A}(\mathfrak{f_n})(w)}d^2w
-\iint\limits_{\Del^*}\frac{g'(z)f'(w)}{(g(z)-f(w))^2}\ov{\mathcal{A}(f)(w)}d^2w\\
 = u_n(\mathfrak{g}_n(z))\mathfrak{g}_n'(z)-u(g(z))g'(z).
\end{gather*}
Here for $z\in \Omega_{n}=\mathfrak{g}_{n}(\Del)$ we set
\begin{align*}
u_n(z) & =
\frac{1}{\pi}\iint\limits_{\Del^*}\frac{\mathfrak{f}_n'(w)}{(z
-\mathfrak{f}_n(w))^2}\ov{\mathcal{A}(\mathfrak{f_n})(w)}
d^2w=-\frac{1}{\pi}
\iint\limits_{\Omega_n^*} \frac{\ov{\mathcal{A}(\mathfrak{f}_n^{-1})(w)}}{(z-w)^2}d^2w,\\
\intertext{and for $z\in\Omega=g(\Del)$,}
u(z)& =\frac{1}{\pi}\iint\limits_{\Del}\frac{f'(w)}{(z
-f(w))^2}\ov{\mathcal{A}(f)(w)}
d^2w=-\frac{1}{\pi}\iint\limits_{\Omega^*}\frac{\ov{\mathcal{A}(f^{-1})(w)}}{(z-w)^2}d^2w.
\end{align*}
Let $\tilde{u}_{n}=\left.u_{n}\right|_{\Omega}$. Using
$\mathcal{A}(\mathfrak{f}_n^{-1})=\mathcal{A}(r_n^{-1}\circ
f^{-1})=\mathcal{A}(f^{-1})$ and $\Omega_{n}^*\subset\Omega^*$, we get
\begin{align*}
\tilde{u}_n(z)-u(z)=\frac{1}{\pi}\iint\limits_{\Omega^*\setminus
\Omega_n^*}\frac{\ov{\mathcal{A}(f^{-1})(w)}}{(z-w)^2}d^2w.
\end{align*}
Since Hilbert transform is an isometry, we obtain
\begin{align*}
&\Vert \tilde{u}_n\circ g\,g'- u\circ g\,g'\Vert_{A_{2}^{1}(\Del)}^2=\iint\limits_{\Del}
\left|\tilde{u}_n(g(z))g'(z)-u(g(z))g'(z)\right|^2
d^2z\\=&\iint\limits_{\Omega}|\tilde{u}_n(z)-u(z)|^2 d^2z
\leq \iint\limits_{\Omega^*\setminus
\Omega_n^*}|\mathcal{A}(f^{-1})(z)|^2d^2z=\iint\limits_{1<|z|<\tfrac{n+1}{n}}
|\mathcal{A}(f)(z)|^2d^2z.
\end{align*}
Since $\mathcal{A}(f)\in A_2^1(\Del^*)$, we get
$$\lim_{n\rightarrow\infty}\Vert \tilde{u}_n\circ g\,g'- u\circ g\,g'\Vert_{A_{2}^{1}(\Del)}=0,$$
and, consequently, $\lim_{n\rightarrow\infty}\tilde{u}_n(z)=u(z)$, uniformly on compact subsets of $\Omega$. 
For every compact subset $E\subset \Del$, $\mathfrak{g}_{n}(E)\subset\Omega$ for $n$ sufficiently large, and
it follows that
$$\lim_{n\rightarrow\infty}u_{n}(\mathfrak{g}_{n}(z))\mathfrak{g}_{n}^{\prime}(z)=u(g(z))g'(z),$$
uniformly on $E$.
\end{proof}
\begin{corollary} \label{1-form2} On $T_{0}(1)$,
$$\partial\tilde{\SSS}_{1}=2\boldsymbol\vartheta.$$
\end{corollary}
\begin{theorem}\label{equalityS} The functions $\SSS_{1}, \tilde\SSS_{1}, \SSS_{2}$ on $T_{0}(1)$ satisfy the following relations,
\begin{align*}
\SSS_2 = -\frac{1}{12\pi}\SSS_1 =-\frac{1}{12\pi}\tilde\SSS_1.
\end{align*}
In particular, in Bers coordinates $\vep_{\mu}$ on the chart $V_{\nu}$ at $[\nu]\in T_{0}(1)$,
\begin{displaymath}
\frac{\partial\SSS_{1}}{\partial\vep_{\mu}}([\nu])=2
\iint\limits_{\Del^*}\mathcal{S}(\g_{\nu})(z)\mu(z)d^2z,
\end{displaymath}
where $w_{\nu}=\g_{\nu}^{-1}\circ \f^{\nu}$ is the conformal welding corresponding to $[\nu]\in T_{0}(1)$.
\end{theorem}
\begin{proof} Since $\SSS_{2}(0)=\tilde\SSS_{1}(0)=0$,
Theorems \ref{varyS2} and  \ref{varyS3} immediately give
\begin{align*}
\SSS_2 =-\frac{1}{12\pi}\tilde\SSS_1.
\end{align*}
Since the function $\SSS_{2}$ is symmetric, 
\begin{align*}
\SSS_2([\mu])=\SSS_2([\mu]^{-1}),
\end{align*}
the function $\tilde\SSS_1$ is also symmetric, so that $\tilde\SSS_{1}=\SSS_{1}$.
\end{proof}
\begin{corollary} \label{1-form3} On $T_{0}(1)$,
$$\pa\SSS_{1}=2\boldsymbol\vartheta.$$
\end{corollary}
\begin{remark} Returning to the model $T(1)=L^{\infty}(\Del^{*})/\sim$, let $\gamma=g^{-1}\circ f\in \mathcal{T}_{0}(1)$. Introducing the operator $\mathbf{K}: \ov{A_{2}^{1}(\Del)}\oplus \ov{A_{2}^{1}(\Del^{*})}\rightarrow A_{2}^{1}(\Del)\oplus A_{2}^{1}(\Del^{*})$,
$$\mathbf{K}=\begin{pmatrix} K_{1} & K_{2} \\
K_{3} & K_{4} \end{pmatrix},
$$
and the vectors 
$$\mathbf{u}=\begin{pmatrix}u_{1}\\ u_{2}
\end{pmatrix},\; \mathbf{v}=\begin{pmatrix} v_{1}
\\ v_{2}\end{pmatrix} \in A_{2}^{1}(\Del)\oplus A_{2}^{1}(\Del^{*}),$$
where $u_{1} = \mathcal{A}(\imath\circ f\circ\imath )\circ\imath\,\imath^{\prime},\,u_{2}= -\mathcal{A}(\imath\circ g\circ\imath )\circ\imath\,\imath'$ and $v_{1}=\mathcal{A}(f),\,v_{2}=-\mathcal{A}(g)$. Applying Lemma \ref{important} to $\gamma$ and $\gamma^{-1}$ and using generalized Grunsky equality, we can succinctly rewrite the two identities as a single equation
\begin{equation*}
\mathbf{K}\,\bar{\mathbf{u}}  = -\mathbf{v}.
\end{equation*}
Indeed, Lemma \ref{important} applied to $\gamma$ and $\gamma^{-1}$ gives
$$K_{3}\bar{u}_{1} + K_{4}\bar{u}_{2}=-v_{2}\quad\text{and}\quad 
K_{1}\bar{v}_{1} + K_{2}\bar{v}_{2} =-u_{1},$$
and from generalized Grunsky equality it follows that the functions
\begin{align*}
w_{1}(z)& =\left(\log\frac{f(z)}{z}\right)'=-\sum_{n=1}^{\infty}nb_{-n,0}z^{n-1} \\
\intertext{and} w_{2}(z)& =-\left(\log\frac{g(z)}{z}\right)'= -\sum_{n=1}^{\infty}nb_{n,0}z^{-n-1}
\end{align*}
satisfy the equations
$$K_{1}\bar{w}_{1} + K_{2}\bar{w}_{2} = w_{1}\quad\text{and}\quad
K_{3}\bar{w}_{1} + K_{4}\bar{w}_{2} = w_{2}.$$
Since $u_{1} =v_{1} -2w_{1}$ and $u_{2}=v_{2} -2w_{2}$, we get the equation 
$\mathbf{K}\bar{\mathbf{u}}=-\mathbf{v}$. Similarly, we get the equation
$\mathbf{K} \bar{\mathbf{v}} = -\mathbf{u}$.
\end{remark}
\begin{remark}
For $C^{3}$ curves the result of Theorem \ref{varyS3} was obtained by
Schiffer and Hawley in \cite{Schiffer3}. They have used a completely different approach which can not be generalized
to quasi-circles for $T_{0}(1)$.
\end{remark}

The equality $\SSS_{2}=-\tfrac{1}{12\pi}\SSS_{1}$ can be also interpreted as a surgery type formula for determinants of elliptic operators (see \cite{burghelea,hassel-zelditch}). Namely, let $\Delta_{\varphi}$ be the Laplace operator of the conformal metric $e^{2\varphi(z)}|dz|^{2}$ on $\Del$ with Dirichlet boundary condition. Its zeta-function regularized determinant $\det\Delta_{\varphi}$ is given by the Polyakov-Alvarez formula
\begin{equation} \label{polyakov} 
\log\det\Delta_{\varphi}=-
\frac{1}{3\pi}\iint\limits_{\Del}|\varphi_{z}|^{2}d^{2}z
-\frac{1}{6\pi}\oint_{S^{1}}\varphi(e^{i\theta})d\theta +\log\det\Delta_{0}.
\end{equation}
Now let $\gamma=\g^{-1}\circ\f\in T_{0}(1)$ and set, as before, $\tilde{\g}=\imath\circ\g\circ\imath,\, \tilde\f=\imath\circ\f\circ\imath$. The metric $|\tilde\g'(z)|^{2}|dz|^{2}$ is a pull-back of the Euclidean metric $|dw|^{2}$ on $\tilde{\Omega}=\tilde\g(\Del)$ by the conformal mapping $\tilde\g$. Assume that $\phi(z)=\tfrac{1}{2}\log|\tilde\g'(z)|^{2}$ is of $C^{1}$ class on $S^{1}$, and denote by $\Delta_{\tilde{\Omega}}$ the Laplace operator of the Euclidean metric on $\tilde{\Omega}$ with Dirichlet boundary condition. From \eqref{polyakov} we immediately get
$$\log\det\Delta_{\tilde{\Omega}} =-\frac{1}{12\pi}\iint\limits_{\Del}|\mathcal{A}(\tilde\g)|^{2}d^{2}z -\frac{1}{3}\log|\tilde\g'(0)| +\log\det\Delta_{\Del}.$$
Now consider the metric $|\tilde\f'(\tfrac{1}{z})|^{2}|dz|^{2}$ on $\Del$ --- a pull-back of the flat metric
$$ds^{2}=\frac{|dw|^{2}}{|\tilde\f^{-1}(w)|^{4}}$$
on $\tilde{\Omega}^{*}=\tilde\f(\Del^{*})$
by the conformal mapping $\tilde\f\circ\imath$. Denoting by $\Delta_{\tilde{\Omega}^{*}}$ the Laplace operator of the metric $ds^{2}$ on $\tilde{\Omega}^{*}$ with Dirichlet boundary condition, we get  from \eqref{polyakov},
$$\log\det\Delta_{\tilde{\Omega}^{*}} =-\frac{1}{12\pi}\iint\limits_{\Del^{*}}|\mathcal{A}(\tilde\f)|^{2}d^{2}z +\log\det\Delta_{\Del^{*}},$$
where we again assumed that $\varphi(z)=\tfrac{1}{2}\log |\tilde\f'(\tfrac{1}{z})|^{2}$
is of $C^{1}$ class on $S^{1}$. Here $\Delta_{\Del^{*}}$ is the Laplace operator of the metric $\tfrac{|dw|^{2}}{|w|^{4}}$ on $\Del^{*}$. Note that the metric $ds^{2}$ is regular at $\infty$, so that $\Delta_{\tilde{\Omega}^{*}}$ is an elliptic operator (cf.~\cite{hassel-zelditch}). The following result now follows from Theorem \ref{equalityS} and the symmetry property $\SSS_{1}([\mu])=\SSS_{1}([\mu^{-1}])$. 
\begin{corollary} \label{surgery-det} 
Let $\gamma=\g^{-1}\circ \f\in T_{0}(1)$ be of $C^{3}$ class on $S^{1}$. Then for $\mathcal{C}=\f(S^{1})$, 
$$\Det_{F}(\mathcal{C})=\frac{\det\Delta_{\tilde{\Omega}}
\det\Delta_{\tilde{\Omega}^{*}}}{\det\Delta_{\Del}\det\Delta_{\Del^{*}}}.$$
\end{corollary}
\begin{remark} The statement of Corollary \ref{surgery-det} can be interpreted as a surgery type formula in the spirit of \cite{burghelea} for the Laplace operator of a conformal metric on the Riemann sphere $\PP^{1}$, which is the Euclidean metric on the interior domain $\tilde{\Omega}=\imath(\Omega^{*})$ and is the metric $ds^{2}=\tfrac{|dw|^{2}}{|\tilde{\f}^{-1}(w)|^{4}}$ on the exterior domain $\tilde{\Omega}^{*}=\imath(\Omega)$ (and thus is continuous on $\PP^{1}$). The Fredholm determinant $\Det_{F}(\mathcal{C})$ is the inverse of the determinant of the Neumann jump operator which corresponds to cutting of $\PP^{1}$ along the contour $\mathcal{C}$ and considering Dirichlet boundary conditions for interior and exterior Laplace operators (cf. \cite{hassel-zelditch}). 
\end{remark}
\section{Weil-Petersson potential}
\subsection{Weil-Petersson potential on $T_{0}(1)$}
As in the case of finite dimensional Teichm\"uller spaces \cite{LT}, it follows from the results of the previous section that the function $\SSS_1$ is a potential for the Weil-Petersson metric on $T_{0}(1)$. For the convenience of the reader, here we give the details.
\begin{theorem}\label{secondvary} In terms of the Bers coordinates on the chart $V_{\kappa}$ at $\kappa\in T_{0}(1)$,
\begin{displaymath}
\frac{\pa^2\SSS_{1} }{\pa \vep_{\mu}\bar{\vep}_{\nu}}([\kappa]) =
\iint\limits_{\Del^*} \mu(z) \ov{\nu(z)}\rho(z)d^2z.
\end{displaymath}
\end{theorem}
\begin{proof} We have
$$\frac{\pa^2\SSS_{1} }{\pa \vep_{\mu}\bar{\vep}_{\nu}}([\kappa])=\left.\frac{\pa}{\pa\bar{\vep}}
\right|_{\vep=0}\frac{\pa\SSS_{1}}{\pa\vep_{\mu}}
([\vep\nu\ast\kappa]).$$
Using Theorem \ref{equalityS} and the fact that at the point $\vep\nu\ast\kappa\in T_{0}(1)$ the vector field $\tfrac{\pa}{\pa\vep_{\mu}}$  on the chart $V_{\kappa}$ is represented by $P(R(\mu,\vep\nu))\in H^{-1,1}(\Del^{*})$ on the chart $V_{\vep\nu\ast\kappa}$ (see Section 3.3 in Part I), we get
\begin{align*}
\frac{\pa\SSS_{1}}{\pa \vep_{\mu}}(\vep\nu*\kappa) & =
2\iint\limits_{\Del^*} \mathcal{S}(g_{\vep\nu*\kappa}) P(R(\mu,
\vep\nu))d^2z\\& =2\iint\limits_{\Del^*}
\left(\mathcal{S}(g_{\vep})\circ w_{\vep\nu}
(w_{\vep\nu})_{z}^2\right)  Q(R(\mu, \vep\nu))(1-|\vep\nu|^2)d^2z,
\end{align*}
where $g_{\vep}^{-1}\circ
f^{\vep}= w_{\vep\nu}\circ w_{\kappa}$, $v_{\vep}
=f^{\vep}\circ f^{-1}$, and $Q(R(\mu,\vep\nu))$ was defined in Section 7.1 in Part I. 
Since $g_{\vep}\circ w_{\vep\nu} =
v_{\vep}\circ g$, we have
\begin{align}\label{identityschwarz}
\mathcal{S}(g_{\vep}) \circ w_{\vep\nu} (w_{\vep\nu})_{z}^2 +
\mathcal{S}(w_{\vep\nu})= \mathcal{S}(v_{\vep})\circ g
(g')^2+\mathcal{S}(g),
\end{align}
and it follows from the standard variational formula that
\begin{align*}
\left.\frac{\pa}{\pa\bar{\vep}}\right\vert_{\vep=0}\mathcal{S}(g_{\vep}) \circ
w_{\vep\nu} (w_{\vep\nu})_{z}^2=
\frac{6}{\pi}\iint\limits_{\Del^*}\frac{\ov{\nu(w)}}{(1-\bar{w}z)^4}d^2w.
\end{align*}
Since according to Theorem 7.4 in Part I
\begin{align*}
\left.\frac{\pa}{\pa \bar{\vep}}\right\vert_{\vep=0}Q(R(\mu,
\vep\nu))
\end{align*}
is an infinitesimally trivial Beltrami differential, we have
\begin{align*}
\frac{\pa^2 }{\pa \vep_{\mu}\bar{\vep}_{\nu}}S_1([\kappa])
=\frac{12}{\pi}\iint\limits_{\Del^*}\iint\limits_{\Del^*}\frac{\mu(z)\ov{\nu(w)}}{(1-\bar{w}z)^4}d^2wd^2z
=\iint\limits_{\Del^*} \mu(z)\ov{\nu(z)}\rho(z)d^2z.
\end{align*}
\end{proof}
\begin{corollary}\label{potential} On $T_{0}(1)$,
$$ \pa\bar{\pa}\SSS_1 = -2i \omega_{WP},$$
where $\omega_{WP}$ is the symplectic form of 
the Weil-Petersson metric. In other words, $\SSS_1$ is a potential of
the Weil-Petersson metric on $T_{0}(1)$.
\end{corollary}
\begin{remark} It follows from Corollary \ref{det=period} and Theorem \ref{equalityS} that on $T_{0}(1)$,
\begin{align*}
\pa\bar{\pa}\log \det N_{\Omega} = \frac{i}{6\pi}\omega_{WP}.
\end{align*}
In the spirit of the last remark in Section 8 of Part I, this result should be compared to the local index theorem for families of $\bar{\partial}$-operators on compact Riemann surfaces,
\begin{align*}
\pa\bar{\pa}\log \det \Delta_0 -\pa\bar{\pa}\log \det N_1 =
-\frac{i}{6\pi}\omega_{WP},
\end{align*}
where  $N_1$ is the period
matrix of $1$-forms on a compact Riemann surface $X$ and $\Delta_0$ is the Laplace operator of the hyperbolic metric on $X$ (see, e.g., \cite{ZT}). 
\end{remark}
\begin{remark} It follows from Corollary \ref{1-form3} that on $T_{0}(1)$,
\begin{align*}
\pa \boldsymbol{\vartheta}=0.
\end{align*}
Here is a direct proof of this result, following our work \cite{LT}.
From equation \eqref{identityschwarz}, we have at $[\kappa]\in T_{0}(1)$,
\begin{align*}
(L_{\nu}\boldsymbol{\vartheta})(z)
&=\left.\frac{\pa}{\pa\vep_{\mu}}\right\vert_{\vep=0}
\mathcal{S}(g_{\vep}) \circ w_{\vep\nu}
(w_{\vep\nu})_{z}^2(z)\\
& =-\frac{12}{\pi} \iint\limits_{\Del^*} \nu(w)
\left(\frac{g'(w)^2g'(z)^2}{(g(w)-g(z))^4}-\frac{1}{(w-z)^4}
\right)d^2w\\
&=-\frac{12}{\pi} \iint\limits_{\Del^*}
\nu(w)\frac{g'(w)^2g'(z)^2}{(g(w)-g(z))^4}d^2w.
\end{align*}
Hence,
\begin{align*}
\pa \boldsymbol{\vartheta}(\mu, \nu)
&=L_{\mu}\boldsymbol{\vartheta}(\nu)-L_{\nu}\boldsymbol{\vartheta}(\mu)\\
&=-\frac{12}{\pi} \iint\limits_{\Del^*}\iint\limits_{\Del^*}
\mu(w)\frac{g'(w)^2g'(z)^2}{(g(w)-g(z))^4}\nu(z)d^2wd^2z \\& \;\;\;\;+\frac{12}{\pi}
\iint\limits_{\Del^*}\iint\limits_{\Del^*}
\nu(w)\frac{g'(w)^2g'(z)^2}{(g(w)-g(z))^4}\mu(z)d^2wd^2z\\
& =0.
\end{align*}
\end{remark}
\subsection{Weil-Petersson potential on $T(1)$} 
The $1$-form $\boldsymbol{\vartheta}$ does not naturally extend to the whole
Hilbert manifold $T(1)$ (since $\left.\boldsymbol{\vartheta}\right|_{[\mu]}\in A_{2}(\Del^{*})$ if and only if $[\mu]\in T_{0}(1)$).  From Theorem \ref{hilbertspaces} we also see that $T_{0}(1)$ is the maximal subset of $T(1)$ on which the function $\SSS_{1}$ is well-defined. However, it is easy to construct a Weil-Petersson potential on $T(1)$ by using right translations. Namely, we
index the components of the Hilbert manifold $T(1)$ by the set $I$
(uncountable) and for every $\alpha\in I$ choose $[\mu_{\alpha}]\in T_{\alpha}(1)=R_{\mu_{\alpha}}T_{0}(1)$ such that $\mu_{0}=0$ for the component $T_{0}(1)$. This represents
$T(1)$ as a disjoint union
\begin{displaymath}
T(1) = \bigsqcup_{\alpha\in I} T_{\alpha}(1).
\end{displaymath}
Define
\begin{displaymath}
\SSS([\nu]) = \SSS_1( [\nu*\mu_{\alpha}^{-1}])\quad\text{for}\quad
[\nu]\in T_{\alpha}(1).
\end{displaymath}
It follows from the right-invariance of the Weil-Petersson metric that
the function $\SSS$ is a Weil-Petersson potential on $T(1)$.
\section{The period mapping} 
The generalization of the classical period mapping to the homogeneous space $\Mob(S^1)\bk\Diff_+(S^1)$ was outlined by
Kirillov and Yuriev in \cite{KY2} and developed by Nag \cite{Nag92}. In particular, in \cite{Nag92} it is explained in what sense  
this is a generalization of classical period mapping as an association between the complex structures and corresponding spaces of holomorphic $1$-forms. Subsequently in \cite{NS}, Nag and Sullivan extended the period mapping to the universal Teichm\"uller space $T(1)$. Here we prove that the Kirillov-Yuriev-Nag-Sullivan (KYNS) period mapping coincides with the mapping $\hat{\cP}$ defined in Remark \ref{Gperiod}.
\subsection{KYNS period mapping}
Following \cite{NS}, let $\mathcal{H}$ be the real Hilbert
space
\begin{align*}
\mathcal{H} = & H^{1/2}(S^1, \R) /\R\\
=& \left\{ f:S^1\rightarrow \R\;\Bigr|\; f(e^{i\theta})=\sideset{}{^{\prime}}\sum_{n=-\infty}^{\infty} c_n e^{in\theta},\;\;\;
\sum_{n=1}^{\infty}n|c_n|^2<\infty\right\},
\end{align*}
and let $\Theta$ be the symplectic form\footnote{We use a different sign convention since our complex structure has a different sign compared to \cite{KY2,Nag92,NS}.} on $\mathcal{H}$:
\begin{displaymath}
\Theta( f,g) =\frac{1}{2\pi} \oint_{S^1} g df .
\end{displaymath}
By complex linearity, the symplectic form $\Theta$ extends to the complexification of $\mathcal{H}$ --- the complex Hilbert space $\mathcal{H}_{\C}$,
\begin{align*}
\mathcal{H}_{\C} = &H^{1/2}(S^1, \C) /\C\\
= &\left\{ f:S^1\rightarrow \C\;\Bigr| f(e^{i\theta})=\sideset{}{^{\prime}}\sum_{n=-\infty}^{\infty}c_n e^{in\theta},\;\;\;\sideset{}{^{\prime}}\sum_{n=-\infty}^{\infty}|n||c_n|^2<\infty\right\}.
\end{align*}
With respect to this symplectic form, the Hilbert space $\mathcal{H}_{\C}$ has a canonical decomposition into two closed isotropic subspaces
\begin{align*}
\mathcal{H}_{\C}= W_+\oplus W_-,
\end{align*}
where
\begin{align*}
W_+= &\left\{ f: S^1\rightarrow \C\; \Bigr|\;
f(e^{i\theta})=\sum_{n=1}^{\infty}a_n e^{in\theta},\;\;\;
\sum_{n=1}^{\infty} n|a_n|^2<\infty\right\},\\
W_-= &\left\{ g: S^1\rightarrow \C\; \Bigr|\;
g(e^{i\theta})=\sum_{n=1}^{\infty}b_n e^{-in\theta},\;\;\;
\sum_{n=1}^{\infty} n|b_n|^2<\infty\right\}.
\end{align*}

 Let $\mathfrak{D}_{\infty}$ be the infinite dimensional analog
of Siegel disk \cite{Siegel},
\begin{gather*}
\mathfrak{D}_{\infty}=\Bigl\{Z\in\cB(W_-,W_+)\; : \; \Theta(Zf,g)=\Theta(Zg,f)\quad\text{and}\quad I-Z\bar{Z}>0\Bigr\}.
\end{gather*}
Here $\cB(W_ -,W_+)$ is the Banach space of all bounded linear operators from $W_-$ to $W_+$, and $\bar{Z}=JZJ: W_{+}\rightarrow W_{-}$, where $J$ is the standard conjugation operator on $\mathcal{H}_{\CC}$ defined by $JW_{+}=W_{-}$.
With respect to the standard bases 
$$\left\{ e_n = \tfrac{1}{\sqrt{n}}
e^{in\theta}\right\}_{n\in\mathbb{N}}\quad\text{and}\quad \left\{f_n=\tfrac{1}{\sqrt{n }}e^{-in\theta}\right\}_{n\in\mathbb{N}}$$ of the subspaces
$W_+$ and $W_-$, an operator $Z\in \mathfrak{D}_{\infty}$ is represented by an
infinite matrix, and the condition $\Theta(Zf,g)=\Theta(Zg,f)$
translates as $Z=Z^t$. Let $\Sp(\mathcal{H})$ be the group
of bounded symplectomorphisms on $\mathcal{H}$. Elements of $\Sp(\mathcal{H})$, extended complex linearly to $\mathcal{H}_{\CC}$, in the basis  $\{e_{n}\}_{n\in\mathbb{N}}\cup\{f_{n}\}_{n\in\mathbb{N}}$ of $\mathcal{H}_{\CC}$ can
be represented by matrices
\begin{align} \label{symplectic}
\begin{pmatrix} A & B \\
\bar{B} & \bar{A}
\end{pmatrix},\quad\text{where}\quad  AA^* - BB^* =I,\;\; AB^t=BA^t.
\end{align}
The group $\Sp(\mathcal{H})$ acts transitively on $\mathfrak{D}_{\infty}$ by
 \begin{align*}
Z \mapsto (AZ+B)(\bar{B}Z+\bar{A})^{-1},
 \end{align*}
and the stabilizer of the point $Z=0$ is the unitary subgroup $U$ of 
$\Sp(\mathcal{H})$
consisting of bounded symplectomorphisms with $B=0$. Thus
the canonical quotient map
$Q: \Sp(\mathcal{H})\rightarrow
\mathfrak{D}_{\infty}$,
\[
Q\left(\begin{pmatrix} A & B \\
\bar{B} & \bar{A}
\end{pmatrix}\right)= \left.(AZ+B)(\bar{B}Z+\bar{A})^{-1}\right|_{Z=0}=B\bar{A}^{-1},
\]
induces the isomorphism
\begin{align*}
\Sp(\mathcal{H})/U \simeq \mathfrak{D}_{\infty}.
\end{align*}

In \cite{NS}, Nag and Sullivan proved that the
assignment
$$\mathrm{Homeo}_{qs}(S^1)\ni\gamma\mapsto \hat{\Pi}(\gamma)\in 
\cB(\mathcal{H}_{\CC}),$$
where 
$$\hat{\Pi}(\gamma)(f)= f\circ
\gamma - \frac{1}{2\pi} \oint_{S^1} f\circ \gamma d\theta, \;\; f\in\mathcal{H}_{\CC},$$
defines a right action of the group $\mathrm{Homeo}_{qs}(S^1)$ on the Hilbert space $\mathcal{H}_{\CC}$ by symplectomorphisms. Thus the mapping 
$$\hat{\Pi}: \mathrm{Homeo}_{qs}(S^1)\rightarrow\Sp(\mathcal{H})$$
satisfies $\hat{\Pi}(\gamma_{1}\circ\gamma_{2})=\hat{\Pi}(\gamma_{2})\hat{\Pi}(\gamma_{1})$.
On the other hand, an operator $\hat{\Pi}(\gamma)$ preserves the subspaces $W_+$ and $W_-$, i.e., $\hat{\Pi}(\gamma)\in U$,  if and only if $\gamma\in\Mob(S^{1})$. The induced mapping 
$$\Pi=Q\circ\hat{\Pi}:
T(1)=\Mob(S^{1})\bk\text{Homeo}_{qs}(S^1)\rightarrow \Sp(\mathcal{H})/U\simeq
\mathfrak{D}_{\infty}$$ 
is what we call KYNS period mapping of $T(1)$.

With respect to the basis $\{e_n\}_{n\in\mathbb{N}}\cup\{f_n\}_{n\in\mathbb{N}}$ of
$\mathcal{H}_{\CC}$, the mapping $\hat{\Pi}:
\text{Homeo}_{qs}(S^1)\rightarrow \Sp(\mathcal{H})$ 
is given by the matrix
\begin{align*}
\hat{\Pi}(\gamma) & = \begin{pmatrix} \mathfrak{A} & \mathfrak{B}, \\
\bar{\mathfrak{B}} & \bar{\mathfrak{A}}
\end{pmatrix},\quad\gamma\in\text{Homeo}_{qs}(S^1), \\
\intertext{where}
\mathfrak{A}_{mn}(\gamma) & = \frac{1}{2\pi }
\sqrt{\frac{m}{n}}\oint_{S^1}
(\gamma(e^{i\theta}))^n e^{-im\theta} d\theta,\\
 \mathfrak{B}_{mn}(\gamma) & =
\frac{1}{2\pi } \sqrt{\frac{m}{n}}\oint_{S^1}
(\gamma(e^{i\theta}))^{-n} e^{-im\theta} d\theta.
\end{align*}
As a result, the KYNS period matrix $\Pi: T(1)\rightarrow \mathfrak{D}_{\infty}$ is given by the matrix 
$$ \Pi([\mu])=\mathfrak{B}\bar{\mathfrak{A}}^{-1},\quad [\mu]\in 
T(1),$$
where $\mathfrak{A}=\mathfrak{A}(w_{\mu})$ and $\mathfrak{B} =\mathfrak{B}(w_{\mu})$.
On the other hand, it follows from \eqref{symplectic} that
\begin{align*}
\hat{\Pi}(\gamma^{-1}) = \hat{\Pi}(\gamma)^{-1} =\begin{pmatrix} \mathfrak{A}(\gamma)^* & -\mathfrak{B}(\gamma)^t\\
-\mathfrak{B}(\gamma)^* & \mathfrak{A}(\gamma)^t
\end{pmatrix}.
\end{align*}
\begin{proposition}\label{Grunsky} 
The Grunsky matrices $B_{l},\,l=1,2,3,4,$ corresponding to $\gamma\in S^1\bk\text{Homeo}_{qs}(S^1)$, and
the elements of the matrix $\hat{\Pi}(\gamma)$ are related by
\begin{alignat*}{2}
B_1 & =\mathfrak{B}\bar{\mathfrak{A}}^{-1}, &\quad
B_2 & =(\mathfrak{A}^*)^{-1},\\
B_3 & = \bar{\mathfrak{A}}^{-1}, &\quad
B_4 & =-\mathfrak{B}^*(\mathfrak{A}^*)^{-1}.
\end{alignat*}
\end{proposition}
\begin{proof} Let $\gamma=g^{-1}\circ f$ be the conformal welding of $\gamma\in S^1\bk\text{Homeo}_{qs}(S^1)$, and 
let $P_n$ and $Q_n$ be the Faber polynomials associated to the pair $(f,g)$. Denoting by $\mathrm{P}_+:\mathcal{H}_{\C}\rightarrow W_+$ and
$\mathrm{P}_-:\mathcal{H}_{\C}\rightarrow W_-$ the orthogonal
projection operators, we get 
$$\mathfrak{A}^{\ast}=\left.\mathrm{P}_+
\hat{\Pi}(\gamma^{-1})\right\vert_{W_+}.$$ 
By definition of Faber polynomials, $\left.(P_{n}^{0}\circ f)\right|_{S^{1}}\in W_{+}$,
where $P^{0}_{n}(z)=P_n(z)-P_n(0)$ and
$n\geq 1$. We have on $S^{1}$,
$$\mathrm{P}_+\hat{\Pi}(\gamma^{-1})(P_n^0\circ f)
=\mathrm{P}_+\left((P_n^0\circ f)\circ
f^{-1}\circ g\right)=\mathrm{P}_+(P_n^0\circ g).$$
Since $\mathrm{P}_+(P_n^0\circ g)(e^{i\theta})=e^{in\theta}$ and
$(P_n^0\circ f)(e^{i\theta})= n\sum_{m=1}^{\infty} b_{-m,n}e^{im\theta}$, we obtain
\begin{align*}
\sum_{k=1}^{\infty}\mathfrak{A}^*_{mk} (B_2)_{kn}=\delta_{mn},
\end{align*}
i.e. $\mathfrak{A}^*B_2=\text{Id}$. Similarly, let $Q_{n}^{0}=Q_{n}-Q_{n}(\infty)$. By definition of Faber polynomials,
$(Q_{n}^{0}\circ f)(e^{i\theta})-e^{-in\theta}\in W_{+}$ for $n\geq 1$.
We have on $S^{1}$, 
$$\mathrm{P}_+\hat{\Pi}(\gamma^{-1})((Q_n^0(f(e^{i\theta})) - e^{-in\theta})
=\mathrm{P}_+\left((Q_n^0\circ g)-(\gamma^{-1})^{-n}
\right)=-\mathrm{P}_+((\gamma^{-1})^{-n}).$$
Since $-\mathfrak{B}(\gamma^{-1})=\mathfrak{B}^t$, we have
$-\mathrm{P}_+((\gamma^{-1})^{-n})(e^{i\theta})=\sum_{k=1}^{\infty}\sqrt{\tfrac{n}{k}}
\mathfrak{B}_{nk}e^{ik\theta}$. Using $Q_n^0(f(e^{i\theta})) - e^{-in\theta}=n\sum_{m=1}^{\infty}b_{-m,-n}e^{im\theta}$, we obtain
\begin{align*}
\sum_{k=1}^{\infty}\mathfrak{A}^*_{mk} (B_1)_{kn}=\mathfrak{B}_{nm},
\end{align*}
i.e., $\mathfrak{A}^*B_1=\mathfrak{B}^t$, which is equivalent to
$B_1=B_1^t=\mathfrak{B}\bar{\mathfrak{A}}^{-1}$. Using
$B_3=B_2^t$ and $B_4(\gamma)=\ov{B_1(\gamma^{-1})}$ concludes the proof.
\end{proof}
\begin{corollary}\label{coincidence}
The KYNS period mapping $\Pi$ coincides
with our period mapping $\hat{\cP}$ defined in Remark \ref{Gperiod}.
\end{corollary}
\begin{proof}
Due to Proposition \ref{Grunsky}, $B_{1}=\mathfrak{B}\bar{\mathfrak{A}}^{-1}$.
\end{proof}
\begin{remark} In \cite{NS} it was stated that the period mapping $\Pi:T(1)\rightarrow \mathfrak{D}_{\infty}$ is a holomorphic mapping of Banach manifolds. However, it was only shown that the induced mapping $D\Pi$ of tangent spaces is  complex linear injection, which is not enough to claim holomorphy for infinite dimensional manifolds. In Appendix B we prove that the mapping $\hat{\cP}: T(1)\rightarrow \cB(\ell^{2})$ is holomorphic, which completes the proof in \cite{NS}.
\end{remark}
Following G.~Segal~\cite{Segal}, we introduce the subgroup $\Sp_{0}(\mathcal{H})$ of the symplectic group $\Sp(\mathcal{H})$
for which $B\in\cS_{2}(W_-,W_+)$ --- the Hilbert space of
Hilbert-Schmidt operators from $W_{-}$ to $W_{+}$. The group 
$\Sp_{0}(\mathcal{H})$ acts transitively on the restricted
Siegel disc
\begin{align*}
\mathfrak{D}_{\infty}^0=\mathfrak{D}_{\infty}\cap\cS_{2}(W_{-},W_{+}).
\end{align*}
Corollaries \ref{coincidence} and \ref{trace-grunsky} immediately imply the following result.
\begin{corollary}
For $[\mu]\in T(1)$, $\Pi([\mu])\in\mathfrak{D}_{\infty}^0$ if and only
if $[\mu]\in T_0(1)$.
\end{corollary}
\begin{remark}
In view of the above corollary, define the restricted period mapping $$\Pi_0=\left.\Pi\right|_{T_{0}(1)}:T_0(1)\rightarrow \mathfrak{D}_{\infty}^0.$$ 
Since by Corollary \ref{coincidence}
$\Pi_{0}=\cP$, by Theorem \ref{holomorphicembedding} $\Pi_{0}$ is a holomorphic mapping of Hilbert manifolds. The homogenuous space $\mathfrak{D}_{\infty}^0$ carries a natural $\Sp_{0}(\mathcal{H})$-invariant K\"{a}hler metric with the K\"{a}hler potential $\Phi(Z)=\log\Det (1- Z\bar{Z})$.
It was first shown by Kirillov and Yuriev~\cite{KY2} and later by Nag \cite{Nag92} that the pullback of this metric to $\Mob(S^1)\bk\Diff_+(S^1)$ by the period mapping coincides, up to a constant, with the
Weil-Petersson metric. It immediately follows from Corollary \ref{coincidence} that
$$\SSS_{2}=\log\Det(I-Z\bar{Z}),$$
so that the pullback of the natural K\"{a}hler metric on $\mathfrak{D}_{\infty}^{0}$ by the restricted period mapping to $T_{0}(1)$ coincides, up to a constant, with the Weil-Petersson metric
on $T_{0}(1)$. Thus we have established the relations between all natural potential functions on
$T_{0}(1)$: up to a constant factor,  they are indeed all equal!
\end{remark}
\subsection{Embeddings into the Segal-Wilson universal Grassmannian} 
Let $\cV$ be an infinite-dimensional separable complex Hilbert space and let
$$\cV=V_{+}\oplus V_{-}$$
be its decomposition into the direct sum of infinite-dimensional closed subspaces $V_{+}$ and $V_{-}$. The Segal-Wilson universal Grassmannian $\mathrm{Gr}(\cV)$ \cite{Segal-Wilson,Pressley-Segal} is defined as a set of closed subspaces $W$ of $\cV$ satisfying the following conditions.
\begin{itemize}
\item[\textbf{UG1.}] The orthogonal projection $\mathrm{pr}_{+}: W\rightarrow V_{+}$ is a Fredholm operator.
\item[\textbf{UG2.}] The orthogonal projection $\mathrm{pr}_{-}: W\rightarrow V_{-}$ is a Hilbert-Schmidt operator.
\end{itemize}
Equivalently, $W\in\mathrm{Gr}(\cV)$, if $W$ is the image of an operator
$\mathrm{w}: V_{+}\rightarrow W$ such that $\mathrm{pr}_{+} \mathrm{w}$
is Fredholm and $\mathrm{pr}_{-} \mathrm{w}$ is Hilbert-Schmidt.
The Segal-Wilson Grassmannian $\mathrm{Gr}(\cV)$ is a Hilbert manifold modeled on
the Hilbert space $\cS_{2}(V_+, V_-)$ of Hilbert-Schmidt operators
from $V_+$ to $V_-$. 

For our purposes, let $\cV=\mathcal{H}_{\C}$ and $V_+=W_-$, $V_-=W_+$. To every $[\mu]\in T_{0}(1)$  we associate a closed subspace $W_{\mu}\subset \mathcal{H}_{\C}$ spanned by the functions $w_n(e^{i\theta})=\frac{1}{\sqrt{n}}Q_n(\f^{\mu}(e^{i\theta}))$, where
$w_{\mu} = \g_{\mu}^{-1}\circ \f^{\mu}$ is the corresponding conformal welding. Explicitly, in terms of the basis $\{e_{n}\}_{n\in\mathbb{N}}\cup\{f_{n}\}_{n\in\mathbb{N}}$ of
$\mathcal{H}_{\C}$,
\begin{align*}
w_n=f_{n}+\sum_{m=1}^{\infty}\sqrt{nm}
b_{-n,-m}e_{m}, \; n\in \mathbb{N}.
\end{align*}
We have $W_{\mu}=\mathrm{w}(V_{+})$, where $\mathrm{w}(f_{n})=w_{n}$,
$n\in\mathbb{N}$. Thus the mapping $\mathrm{pr}_{+}\mathrm{w}=I$ --- the identity operator on $V_{+}$, is obviously Fredholm, and the mapping $\mathrm{pr}_-\mathrm{w}=B_{1}(\f^{\mu})$ 
is Hilbert-Schmidt since $[\mu]\in T_{0}(1)$. According to Theorem \ref{holomorphicembedding}, the mapping
$$\cE: T_{0}(1)\rightarrow \mathrm{Gr}(\mathcal{H}_{\C})$$
given by $\cE([\mu])=W_{\mu}$ is a holomorphic inclusion of
$T_{0}(1)$ into the Segal-Wilson universal Grassmannian. For the homogeneous space $\Mob(S^{1})\bk\Diff_{+}(S^{1})$ this mapping was first considered in \cite{KY2}. 
\begin{remark} Seemingly another mapping of $\Mob(S^{1})\bk\Diff_{+}(S^{1})$ into the Segal-Wilson Grassmannian was considered in \cite{Todorov}. Namely, extend $\mu\in L^{\infty}(\Del^{\ast})$ by zero to $\Del$ and 
let $V_{\mu}$ be the space of distrubitional solutions
of the Beltrami equation $w_{\z}=\mu w_{z}$ on $\C$ having a single pole at $0$. The mapping in \cite{Todorov} was defined by the assignment
$[\mu]\rightarrow \left.V_{\mu}\right|_{S^{1}}$. It is easy to see that the space $V_{\mu}$ is spanned by the functions  
$w_{n}(z)=Q_{n}(\f^{\mu})(z)$, $n\in\mathbb{N}$, so that
$\left.V_{\mu}\right|_{S^{1}}=W_{\mu}$ and the mapping in \cite{Todorov} coincides with the Kirillov-Yuriev mapping \cite{KY2}.
\end{remark}
\begin{remark}The inclusion $\cE: T_{0}(1)\rightarrow \mathrm{Gr}(\mathcal{H}_{\C})$ is a holomorphic mapping due to the holomorphic dependence of $\f^{\mu}$ on $\mu$. 
Since $\f^{\mu}$ is holomorphic on $\Del$, the subspaces $W_{\mu}$
correspond to the different uniformizations of the same Riemann surface $\Omega=\f^{\mu}(\Del)\simeq \Del$. However, one can consider another mapping where the associated subspaces in the universal Grassmannian correspond to Riemann surfaces of different complex structure. 
Namely, set, as before, $\cV=\mathcal{H}_{\C}$ and let
$V_+=W_+$ and $V_-=W_-$. We denote the corresponding Segal-Wilson Grassmannian by $\widetilde{\mathrm{Gr}}(\mathcal{H}_{\C})$, and define the mapping
$$\tilde{\cE} : T_{0}(1)\rightarrow \widetilde{\mathrm{Gr}}(\mathcal{H}_{\C})$$
by assigning to every point $[\mu]\in T_{0}(1)$ the closed subspace
$\tilde{W}_{\mu}\subset\mathcal{H}_{\C}$ spanned by the functions
$\tilde{w}_{n}(e^{i\theta})=\tfrac{1}{\sqrt{n}}P_{n}(\g_{\mu}(e^{i\theta}))$, $n\in\mathbb{N}$. We have $\tilde{W}_{\mu}=\tilde{\mathrm{w}}(V_{+})$, where $\tilde{\mathrm{w}}(e_{n})=\tilde{w}_{n}$,
$n\in\mathbb{N}$, and $\mathrm{pr}_{+}\tilde{\mathrm{w}}=I$ --- the identity operator on $V_{+}$ and $\mathrm{pr}_-\tilde{\mathrm{w}}=B_{4}(\g_{\mu})$ 
is Hilbert-Schmidt since $[\mu]\in T_{0}(1)$. The mapping $\tilde{\cE}$ is not holomorphic. However, since $JW_{+}=W_{-}$, where $J$ is the standard conjugation operator on $\mathcal{H}_{\C}$, we have $\widetilde{\mathrm{Gr}}(\mathcal{H}_{\C})=J(\mathrm{Gr}(\mathcal{H}_{\C}))$, so that $\widetilde{\mathrm{Gr}}(\mathcal{H}_{\C})$ is a mirror image of $\mathrm{Gr}(\mathcal{H}_{\C})$. Denoting by $\mathcal{I}$ the inversion on the topological group $T_{0}(1)$, we get
$$\tilde{\cE}=J\circ\cE\circ\mathcal{I}.$$
\end{remark}
\begin{remark} One can describe the ``Schottky locus'', i.e., the image $\cE(T_{0}(1))$ in the Segal-Wilson Grassmannian $\mathrm{Gr}(\mathcal{H}_{\C})$. Indeed, since the corresponding points in $\mathrm{Gr}(\mathcal{H}_{\C})$ are associated with the Grunsky operators $B_{1}$, it is equivalent to the characterization of the image of the restricted period map $\Pi_{0}: T_{0}(1)\rightarrow \mathfrak{D}^{0}_{\infty}$. Let $C=\{C_{mn}\}_{m,n\in\mathbb{N}}\in\mathfrak{D}^{0}_{\infty}$, which we realized as symmetric, Hilbert-Schmidt operator on $\ell^{2}$ satisfying $I-C\bar{C}>0$. Then $C\in\Pi_{0}(T_{0}(1))$ if and only if the the following conditions are satisfied.
\begin{itemize}
\item[\textbf{S1.}] 
\begin{align*}
1+ \sum_{m=1}^{\infty} \frac{C_{m1}}{\sqrt{m}}\frac{z_1^{m}-
z_2^{m}}{z^{-1}_1-z^{-1}_2} = \exp\left( -\sum_{m,n=1}^{\infty}\frac{
C_{mn}}{\sqrt{mn}} z_1^{m} z_2^{n}\right).
\end{align*}
\item[\textbf{S2.}] There exist $D=\{D_{mn}\}_{m,n\in\mathbb{N}}\in\mathfrak{D}^{0}_{\infty}$ and $B\in\cB(\ell^{2})$ such that
$$I-C\bar{C}=BB^{\ast}\quad\text{and}\quad I-D\bar{D}=\bar{B}^{\ast}\bar{B}.$$
\item[\textbf{S3.}] 
\begin{align*}
1 + \sum_{m=1}^{\infty} \frac{D_{m1}}{\sqrt{m}}\frac{z_1^{-m}-
z_2^{-m}}{z_1-z_2} = \exp\left( -\sum_{m,n=1}^{\infty}\frac{
D_{mn}}{\sqrt{mn}} z_1^{-m} z_2^{-n}\right).
\end{align*}
\end{itemize}
Equations \textbf{S1} and \textbf{S3} are understood as infinite sequence of relations between elements of the matrices $C$ and $D$ obtained by comparing coefficients of $z^{m}_{1}z^{n}_{2}$ and $z^{-m}_{1}z^{-n}_{2}$ respectively.  Equations \textbf{S1} and \textbf{S3} are nothing but 
dispersionless Hirota equations (see, e.g., \cite{Teo2}).
They are just a reformulation of the definition of
the Grunsky coefficients of the univalent functions $f$ and $g$ and the identities
\begin{align*}
\frac{1}{f(z)}+c & =Q_1(f(z))=\frac{1}{z}+\sum_{m=1}^{\infty}b_{-1,-m}z^{m},\\
\frac{g(z)}{b} +d & =P_{1}(g(z)) =z + \sum_{m=1}^{\infty}b_{1m}z^{-m},
\end{align*}
where $c$ and $d$ are constants. See \cite{Teo2} for details.
\end{remark}
\appendix
\section{Hilbert manifold structure of $\mathcal{T}_{0}(1)$}
Here we show that the Hilbert manifold $\mathcal{T}_{0}(1)$, modeled on the Hilbert 
space $A_2(\Del)\oplus \C$, can
also be modeled on the Hilbert space $A^{1}_{2}(\Del)$, which
induce the same Hilbert manifold structure. This result is
parallel to the one in the Appendix of \cite{Teo}. 

Let $\beta$ be the Bers embedding
$\mathcal{T}(1)\hookrightarrow A_\infty(\Del)\oplus \C$, 
$$\mathcal{T}(1)\ni \gamma=g^{-1}\circ f\mapsto \left(\mathcal{S}(f), \tfrac{1}{2}\mathcal{A}(f)(0)\right),$$
and let
$\hat{\beta}: \mathcal{T}(1)\rightarrow A^{1}_{\infty}(\Del)$ be the pre-Bers embedding of $\mathcal{T}(1)$ into $A^{1}_{\infty}(\Del)$,
$$ \mathcal{T}(1)\ni \gamma=g^{-1}\circ f\mapsto \mathcal{A}(f).$$
By Theorem \ref{hilbertspaces}, $\hat{\beta}(\gamma)\in A_{2}^1(\Del)$ if and only 
if $\gamma\in \mathcal{T}_0(1)$.
\begin{lemma}\label{holomorphichilbert}
The map $\widehat{\Psi} : A_{2}^1(\Del)\rightarrow A_2(\Del)\oplus
\C$, 
\begin{displaymath}
\widehat{\Psi}(\psi) = \left(\Psi(\psi),\tfrac{1}{2}\psi(0)\right),
\end{displaymath}
where $\Psi(\psi)=\psi_{z}-\tfrac{1}{2}\psi^{2}$,
is a one to one holomorphic mapping on Hilbert spaces.
\end{lemma}
\begin{proof}
Firstly, the map $\Psi: A_{2}^1(\Del)\rightarrow A_2(\Del)$
is holomorphic. That is, 
for every $\psi, \varphi\in A_{2}^1(\Del)$, the map $\C\ni
t\mapsto \Psi(\psi+t\varphi)$ is holomorphic in a neighbourhood of
$0\in \C$. Indeed,
\begin{align*}
&\left\Vert
\frac{\Psi(\psi+t\varphi)-\Psi(\psi+t_0\varphi)}{t-t_0}-\frac{d}{dt}\Bigr\vert_{t=t_0}
\Psi(\psi+t\varphi)\right\Vert_{A_{2}(\Del)}\\
=&\frac{|t-t_0|}{2}\Vert\varphi^2\Vert_{A_{2}(\Del)}\leq\frac{|t-t_0|}{4}\Vert
\varphi \Vert_{A^{1}_{\infty}(\Del)} \Vert \varphi\Vert_{A^{1}_{2}(\Del)}
=O(|t-t_{0}|).
\end{align*}
Secondly, by Lemma \ref{subspace1},
\begin{align*}
|\psi(0)|\leq \sqrt{\frac{1}{\pi}}\Vert \psi\Vert_{A_{2}^{1}(\Del)},
\end{align*}
so that $\psi \mapsto \tfrac{1}{2}\psi(0)$ 
is a bounded complex-linear map. The injectivity of $\widehat{\Psi}$ has been proved in the Appendix of \cite{Teo}. 
\end{proof}
\begin{corollary}\label{open1}
The set $\hat{\beta}(\mathcal{T}_0(1))\subset A_{2}^1(\Del)$ is open in $A_{2}^{1}(\Del)$.
\end{corollary}
\begin{proof}
 It readily follows from the results in Section 3.3 of Part I that $\beta(\mathcal{T}_0(1))$ is open in $A_2(\Del)\oplus \C$.
The assertion now follows from the lemma above since
$\hat{\beta}(\mathcal{T}_0(1))=\widehat{\Psi}^{-1}\bigl(\beta(\mathcal{T}_0(1))\bigr)$. 
\end{proof}
\begin{theorem}\label{biholomorphic}
The embeddings $\beta: \mathcal{T}_{0}(1)\hookrightarrow A_2(\Del)\oplus \C$ and $\hat{\beta}: \mathcal{T}_{0}(1)\hookrightarrow A_{2}^1(\Del)$ induce the same Hilbert manifold structure on $\mathcal{T}_0(1)$.
\end{theorem}
\begin{proof}
The map $\widehat{\Psi}: \hat{\beta}(\mathcal{T}_0(1))\rightarrow
\beta(\mathcal{T}_0(1))$ is a holomorphic bijection between complex manifolds. To show that $\widehat{\Psi}$ is biholomorphic, by inverse
function theorem (see, e.g., \cite{Lang2}) we need to prove that for every $\psi\in \hat{\beta}(\mathcal{T}_0(1))$ the linear map
$D_{\psi} \widehat{\Psi}$ is a topological isomorphism between the
Hilbert spaces $A_{2}^1(\Del)$ and $A_{2}(\Del)\oplus\C$.
Let  $\gamma=g^{-1}\circ f\in\mathcal{T}_{0}(1)$ and $\psi= \mathcal{A}(f)\in\hat{\beta}(\mathcal{T}_{0}(1))$. The linear map $D_{\psi}\widehat{\Psi}:
A_2^1(\Del)\rightarrow A_2(\Del)\oplus\C$ is given by
\begin{align*}
\varphi \mapsto \left(\varphi_z -\psi \varphi,
\tfrac{1}{2}\varphi(0)\right).
\end{align*}
 For every $(\phi, c)\in A_{2}(\Del)\oplus
\C$, the holomorphic function $\varphi$ on $\Del$, defined by
\begin{align*}
\varphi(z) = f'(z) \left(\int_{0}^z \frac{\phi(u)}{f'(u)} du +
2c\right),
\end{align*}
satisfies
\begin{align*}
\varphi_z -\psi \varphi=\phi\quad\text{and}\quad \tfrac{1}{2}\varphi(0)=c.
\end{align*}
We claim that $\varphi\in A_{2}^1(\Del)$, so that the map $D_{\psi}\hat{\Psi}$ is onto. Indeed, repeating the proof of  Lemma \ref{sub2}, we get for $z\in\Del$,
$$|\phi(z)|^{2}\geq \tfrac{1}{2}|\varphi_{z}|^{2} -|\psi(z)|^{2}|\varphi(z)|^{2}.$$
By Becker-Pommerenke theorem, there exists $r'>0$ such that
$$|(1-|z|^{2})\psi(z)|\leq \frac{1}{2\sqrt{2}}\quad\text{for all}\quad r'<|z|<1.$$
Thus for $r'<|z|<1$,
$$2(1-|z|^{2})^{2}|\phi(z)|^{2}\geq (1-|z|^{2})^{2}|\varphi_{z}(z)|^{2}-\tfrac{1}{4}|\varphi(z)|^{2},$$
and the result follows as in the proof of Lemma \ref{sub2}.
Uniqueness theorem for differential equations guarantees that the map
$D_{\psi}\hat{\Psi}$ is one-to-one. Finally, by using the same arguments as in the proof of Lemma \ref{sub1} and Lemma \ref{subspace1}, there exists $C>0$ such that for all $\varphi\in A_{2}^{1}(\Del)$, 
\begin{align*}
\Vert D_{\psi}\hat{\Psi}(\varphi)\Vert_{A_{2}(\Del)} \leq C \Vert
\varphi\Vert_{A_{2}^{1}(\Del)}.
\end{align*}
Hence $D_{\psi}\hat{\Psi}$ is a bounded linear bijection between
Hilbert spaces.
\end{proof}
\begin{corollary}\label{convergea} Let $\{\gamma_{n}
\}_{n=1}^{\infty}$ be a sequence of points in $\mathcal{T}_0(1)$, 
$\gamma_n=g_n^{-1}\circ f_n$,
 and let $\gamma=g^{-1}\circ f\in\mathcal{T}_{0}(1)$. 
Then the following conditions are equivalent. 
\begin{itemize} 
\item[(i)] In $\mathcal{T}_{0}(1)$ topology,
$$\lim_{n\rightarrow\infty}\gamma_n=\gamma.$$ 
\item[(ii)] In $A^{1}_{2}(\Del)$ topology,
$$\lim_{n\rightarrow\infty}\mathcal{A}(f_n)(z) = \mathcal{A}(f)(z).$$
\item[(iii)] In $A^{1}_{2}(\Del^{*})$ topology,
\begin{displaymath}
\lim_{n\rightarrow\infty}\left(\mathcal{A}(g_n)(z)-2\frac{g_n'(z)}{g_n(z)}+\frac{2}{z}\right)
=\mathcal{A}(g)(z)-2\frac{g'(z)}{g(z)}+\frac{2}{z}.
\end{displaymath}
\end{itemize}
\end{corollary}
\begin{proof} The equivalence (i)$\Leftrightarrow$(ii) follows from Theorem \ref{biholomorphic}. Since $\mathcal{T}_{0}(1)$ is a topological group, $\lim_{n\rightarrow\infty}\gamma_{n}=\gamma$ if and only if $\lim_{n\rightarrow\infty}\gamma_{n}^{-1}=\gamma^{-1}$. Now let $j(z)=\tfrac{1}{\bar{z}}$ and let $r$ be the dilation $z\mapsto \overline{g'(\infty)}\,z$.
We have $\gamma^{-1}=\tilde{g}^{-1}\circ \tilde{f}$, where $\tilde{f} =r\circ j\circ g \circ j$ and
\begin{align*}
\mathcal{A}(\tilde{f}) =\mathcal{A}(r\circ j\circ
g \circ j)=\ov{\left(
\mathcal{A}(g)-2\frac{g'}{g}+2\bar{j}\right)\circ
j j_{\z}},
\end{align*}  
so that the equivalence (i)$\Leftrightarrow$(iii) follows from the equivalence (i)$\Leftrightarrow$(ii).
\end{proof}
Next, consider the mappings $\beta^{*}: T_{0}(1)\rightarrow A_{2}(\Del^{*})$,
$$\beta^{*}(\gamma)=\ov{\beta(\gamma^{-1})\circ j j_{\bar{z}}^{2}}=\mathcal{S}(g),$$
and $\hat{\beta}^{*}: \mathcal{T}_{0}(1) \rightarrow A^{1}_{2}(\Del^{*})$,
$$\hat{\beta}^{*}(\gamma)=\mathcal{A}(g),$$
where $\gamma=g^{-1}\circ f\in\mathcal{T}_{0}(1)$. Also,
consider the mapping $\Psi^{*}: A^{1}_{2}(\Del^{*})\rightarrow A_{2}(\Del^{\ast})$, defined by
$$\Psi^{*}(\psi)=\psi_{z} - \tfrac{1}{2}\psi^{2},$$
and let 
$$\widetilde{A^{1}_{2}}(\Del^{\ast}) =\left\{\psi\in A_{2}^{1}(\Del^{\ast}) : \psi(z)=O\left(\frac{1}{z^{3}}\right)\;\;\text{as}\;\;z\rightarrow\infty\right\}.$$
We have the following result.
\begin{lemma} \label{biholomorphic2}
\noindent
\begin{itemize}
\item[(i)] The map $\Psi^{\ast}: A^{1}_{2}(\Del^{*})\rightarrow A_{2}(\Del^{\ast})$ is a holomorphic mapping on Hilbert spaces and its restriction to the subspace $\widetilde{A^{1}_{2}}(\Del^{\ast})$ is injective.
\item[(ii)] The set $\hat{\beta}^{*}(\mathcal{T}_{0}(1))=\hat{\beta}^{*}(T_{0}(1))$ is open in $A^{1}_{2}(\Del^{*})$.
\end{itemize}
\end{lemma}
\begin{proof} Holomorphy of $\Psi^{\ast}$ is proved along the same lines as Lemma \ref{holomorphichilbert}. From the proof of Theorem A.5 in \cite{Teo} it follows that the restriction of the map  $\Psi^{\ast}$ to the subspace $\widetilde{A^{1}_{2}}(\Del^{\ast})$ is one to one. To prove part (ii),  observe that $\beta^{\ast} = \Psi^{\ast}\circ \hat{\beta}^{\ast}$ and $\beta^{\ast}(\mathcal{T}_{0}(1))=\beta^{\ast}(T_{0}(1))$. Since $\hat{\beta}^{*}(\mathcal{T}_{0}(1)), \,\hat{\beta}^{*}(T_{0}(1))\in \widetilde{A^{1}_{2}}(\Del^{\ast})$ and the restriction of $\Psi^{\ast}$ to 
$\widetilde{A^{1}_{2}}(\Del^{\ast})$ is injective, we have the equality $\hat{\beta}^{*}(\mathcal{T}_{0}(1))=\hat{\beta}^{*}(T_{0}(1))$. The proof that this set is open in $A_{2}^{1}(\Del^{\ast})$ is analogous to the proof of Corollary \ref{open1}.
\end{proof}
\begin{corollary} \label{convergeb}
Let $\{\gamma_{n}
\}_{n=1}^{\infty}$, $\gamma_n=g_n^{-1}\circ f_n$, be a sequence of points in $\mathcal{T}_0(1)$ such that
$$\lim_{n\rightarrow\infty}\gamma_n=\gamma=g^{-1}\circ f\in\mathcal{T}_{0}(1).$$ 
Then the following statements hold. 
\begin{itemize} 
\item[(i)] In $A_{2}(\Del^{*})$ topology,
$$\lim_{n\rightarrow\infty}\mathcal{S}(g_n) = \mathcal{S}(g).$$
\item[(ii)] In $A^{1}_{2}(\Del^{*})$ topology,
\begin{displaymath}
\lim_{n\rightarrow\infty}\mathcal{A}(g_n)
=\mathcal{A}(g).
\end{displaymath}
\end{itemize}
\end{corollary}
\begin{proof} Since $\mathcal{T}_{0}(1)$ is a topological group,
$\lim_{n\rightarrow\infty}\gamma_{n}^{-1}=\gamma^{-1}=\tilde{g}^{-1}\circ \tilde{f}$. We have $\tilde{f}=r\circ j\circ g \circ j$, so that
$$\mathcal{S}(\tilde{f}) =\ov{\mathcal{S}(g)\circ j j_{\z}^{2}},$$
which proves part (i). Part (ii) follows from Lemma \ref{biholomorphic2}. 
\end{proof}
\section{The period mapping $\hat{\cP}$}
Let $\cS_{\infty}$ be the closed ideal of compact operators in the Banach algebra $\cB(\ell^2)$ of bounded operators on $\ell^2$. Here we prove that the period mapping $\hat{\cP}:T(1)\rightarrow
\cB(\ell^2)$, defined in Remark \ref{Gperiod}, is a holomorphic mapping of complex Banach manifolds and that  
$$\hat{\cP}^{-1}(\cS_{\infty})= S=\Mob(S^1)\bk\text{Homeo}_{s}(S^1).$$
\begin{theorem}\label{holembed}
The inclusion $\hat{\cP}: T(1)\rightarrow \cB(\ell^{2})$ is a holomorphic mapping of Banach manifolds.
\end{theorem}
\begin{proof}
As in the proof of Theorem \ref{holomorphicembedding}, we will show that for every
$[\nu]\in T(1)$ and $\mu\in \Omega^{-1,1}(\Del^*)$, the map $\C\ni
t \mapsto B_1(t)=B_1(f^{\nu+t\mu})$ is holomorphic in a
neighborhood of $t=0$ in $\C$.  Choose $\delta>0$ so that
$\Vert \nu+t\mu\Vert_{\infty}<1$ for all $|t|<\delta$. For every
$t_0$ such that $|t_0|<\delta$, let $\delta_1$ be such that
$0<\delta_1< \delta-|t_0|$. Then for all $|t-t_0|<\delta_1$, we
have as in Theorem \ref{holomorphicembedding},
\begin{gather*}
\left(K_1^{\nu+t\mu}-K_1^{\nu+t_0\mu} -(t-t_0)\left.
\frac{d}{dt}\right\vert_{t=t_0}K_1^{\nu+t\mu}\right)(z,w)\\
=\frac{(t-t_0)^2}{2\pi i}
\oint_{|\zeta-t_0|=\delta_1}\frac{K_1^{\nu+\zeta\mu}(z,w)}{(\zeta-t)(\zeta-t_0)^2}d\zeta.
\end{gather*}
This gives
\begin{gather*}
\left\Vert \frac{B_1(\f^{\nu+t\mu})-B_1(\f^{\nu+t_0\mu})}{t-t_0} -
\left.\frac{d}{dt}\right\vert_{t=t_0}B_1(\f^{\nu+t\mu})\right\Vert\\
=\sup_{\Vert u \Vert_{2}=1}\left( \iint\limits_{\Del}
\left|\iint\limits_{\Del}\frac{t-t_0}{2\pi i}
\oint_{|\zeta-t_0|=\delta_1}\frac{K_1^{\nu+\zeta\mu}(z,w)\ov{u(w)}}{(\zeta-t)(\zeta-t_0)^2}d\zeta
d^2w\right|^2 d^2z\right)^{1/2}\\
\leq \frac{|t-t_0|}{2\pi}\sup_{\Vert u \Vert_{2}=1}\Biggl(
\iint\limits_{\Del}\Biggl(\oint_{|\zeta-t_0|=\delta_1}\frac{|d\zeta|}{|\zeta-t|^2|\zeta-t_0|^4}\Biggr)\\
\Biggl(\oint_{|\zeta-t_0|=\delta_1}\left|\iint\limits_{\Del}K_1^{\nu+\zeta\mu}(z,w)\ov{u(w)}
d^2w\right|^2 |d\zeta| \Biggr)d^2z\Biggr)^{1/2}\\
=\frac{|t-t_0|}{2\pi}\Biggl(
\oint_{|\zeta-t_0|=\delta_1}\frac{|d\zeta|}{|\zeta-t|^2|\zeta-t_0|^4}\Biggr)^{1/2}
\sup_{\Vert u \Vert_{2}=1}\left(\oint_{|\zeta-t_0|=\delta_1}\Vert
K_1^{\nu+\zeta\mu}\bar{u}\Vert_{2}^2 |d\zeta|\right)^{1/2}.
\end{gather*}
Since $ \Vert K_1\Vert <1$, we obtain
\begin{gather*}
\left\Vert \frac{B_1(\f^{\nu+t\mu})-B_1(\f^{\nu+t_0\mu})}{t-t_0} -
\left.\frac{d}{dt}\right\vert_{t=t_0}B_1(\f^{\nu+t\mu})\right\Vert\\
\leq 
\frac{|t-t_0|}{2\pi}\left(\oint_{|\zeta-t_0|=\delta_1}\frac{|d\zeta|}{|\zeta-t|^2|\zeta-t_0|^4}
\oint_{|\zeta-t_0|=\delta_1}|d\zeta|\right)^{1/2}\\ =O(t-t_{0})\quad\text{as}\quad t\rightarrow t_{0}.
\end{gather*}
\end{proof}

To prove that $\hat{\cP}(S)\subset\cS_{\infty}$, we first give a characterization of
the submanifold $S=\Mob(S^1)\bk\text{Homeo}_s(S^1)$ of $T(1)$. It has been shown by Gardiner and Sullivan \cite{GaSu} that $\beta(S)=A_{\infty}^0(\Del)\cap\beta(T(1))$, where 
$\beta:T(1)\rightarrow A_{\infty}(\Del)$  is the Bers embedding and $A_{\infty}^0(\Del)$ is the subspace of the Banach space $A_{\infty}(\Del)$, defined by
\[
A_{\infty}^0(\Del)=\left\{\phi\in A_{\infty}(\Del) \, : \;\,
\lim\limits_{|z|\rightarrow 1^-} (1-|z|^2)^2\phi(z)=0 \right\}.
\]
Analogous to Theorem A.1 in Part I, we have the following result.
\begin{lemma}\label{close}\label{closure}
The closure of the homogeneous space $\Mob(S^1)\bk\Diff_+(S^1)\subset T(1)$ in the Banach manifold topology is the Banach submanifold $S$ of $T(1)$.
\end{lemma}
\begin{proof}
For $\phi\in A_{\infty}^0(\Del)\cap \beta(T(1))$, let $\phi_{n}=\phi\circ r_{n}$, where $r_{n}$ is the dilation $z\mapsto\tfrac{n}{n+1}z,\, n\in\mathbb{N}$. Since $\phi\in A_{\infty}^{0}(\Del)$, for every $\vep>0$ there exists $0<r<1$ such that 
\[
\sup_{r\leq |z|\leq 1} (1-|z|^2)^2|\phi(z)|<\frac{\vep }{4}.
\]
Thus there exists $N'$ such that 
\[
\sup_{r'\leq |z|\leq 1} (1-|z|^2)^2|\phi_{n}(z)|<\frac{\vep }{4}
\]
for $n>N'$, where $r'=\tfrac{1+r}{2}$.
The sequence $\{\phi_{n}\}$ converges uniformly to $\phi$ on compact subsets of $\Del$, so that there exists $N''$ such that
\[
\sup_{|z|\leq r'}(1-|z|^2)^2|\phi_n(z)-\phi(z)|<\frac{\vep}{2}\quad\text{for}\quad n>N''.
\]
Thus $\Vert\phi_{n}-\phi\Vert_{\infty}<\vep$ for $n>N=\max\{N',N''\}$, so that
$$\lim_{n\rightarrow\infty}\phi_{n}=\phi$$
in the $A_{\infty}(\Del)$ topology. Since $\beta(T(1))$ is open in $A_{\infty}(\Del)$,
$\phi_{n}\in\beta(T(1))$ for large enough $n$. The functions $\phi_{n}$ are smooth on $S^{1}$ (in fact analytic), so that corresponding $\gamma_{n}\in\Mob(S^{1})\bk\mathrm{Homeo}_{qs}(S^{1})$ are also smooth on $S^{1}$. This proves that $\ov{\Mob(S^{1})\bk\Diff_{+}(S^{1})}=S$.
\end{proof}
\begin{remark} Together with Theorem A.1 in Part I, Lemma \ref{close}
explains the distinguished role of the embedded manifold $\Mob(S^1)\bk\Diff_+(S^1)\hookrightarrow T(1)$ in \Te theory.  
Its closure in $T(1)$ under the Banach manifold topology is the Banach submanifold $S$, whereas its closure under the Hilbert manifold topology is the Hilbert submanifold $T_0(1)$.
\end{remark}
\begin{theorem}\label{half} The image of the Banach submanifold $S$ under the KYNS period mapping 
$\hat{\cP}:T(1)\rightarrow \cB(\ell^2)$ is given by 
$$\hat{\cP}(S)=\cS_{\infty}\cap\hat{\cP}(T(1)),$$
where $\cS_{\infty}$ is the space of compact operators on $\ell^{2}$. 
\end{theorem}
\begin{proof}
It is easy to show that $\hat{\cP}(S)\subset\cS_{\infty}$. Indeed, 
by Theorem \ref{traceclass}, $\hat{\cP}(\Mob(S^1)\bk\Diff_+(S^1))\subset\cS_{2}\subset\cS_{\infty}$.  Since the mapping $\hat{\cP}$ is continuous (actually, holomorphic), using Lemma \ref{closure} proves the claim. 

To prove the converse inclusion $\hat{\cP}^{-1}(\cS_{\infty}\cap\hat{\cP}(T(1)))\subset S$, we use methods developed by Brazilevic in \cite{Bra}. Let $\mathcal{U}$ be the space of univalent functions on $\Del$. Following \cite{Bra}, consider
the following function $\FFF:\mathcal{U}\times\mathcal{U}\times\Del \rightarrow\RR$,
\begin{align*}
\FFF(f_1, f_2)(z)=\sqrt{\pi}
(1-|z|^2)\left(\iint\limits_{\Del}\left|K_1(f_1)(z,w)-K_1(f_2)(z,w)\right|^2
d^2w\right)^{1/2}.
\end{align*}
When $f_1=f$ and $f_2=\text{id}$ --- the identity mapping, we denote
\begin{align*}
\FFF(f)( z)=\FFF(f, \mathrm{id})(
z)=\sqrt{\pi}(1-|z|^2)\KKK_1(z,z)^{1/2}.
\end{align*}
In \cite{Bra}, Brazilevic has introduced a new metric on $\mathcal{U}$,
\begin{align*}
d(f_1,f_2)=\sup_{z\in \Del} \FFF(f_1,f_2)( z),
\end{align*}
and has shown that
\begin{align*}
\Vert \mathcal{S}(f_1)-\mathcal{S}(f_2)\Vert_{\infty}\leq 6 d(f_1,
f_2).
\end{align*}
For fixed $\zeta\in \Del$, consider the kernel 
\begin{align*}
K_1(f)(z,\zeta)=\frac{1}{\pi}\sum_{n=1}^{\infty}n\left(\sum_{m=1}^{\infty}m b_{-n,-m}\zeta^{m-1}\right)z^{n-1}
\end{align*}
as a holomorphic function on $\Del$.
By Grunsky inequality,
\begin{align*}
\Vert K_{1}(f)(\;\cdot\;,\zeta)\Vert_2^2&=\KKK_1(f)(\zeta,\zeta)=\frac{1}{\pi}
\sum_{n=1}^{\infty}n
\left|\sum_{m=1}^{\infty}mb_{-n,-m}\zeta^{m-1}\right|^2 \\
&\leq\frac{1}{\pi}\sum_{n=1}^{\infty} 
n|\zeta|^{2n-2}=\frac{1}{\pi(1-|\zeta|^2)^2}<\infty,
\end{align*}
so that $K_{1}(f)(\;\cdot\;,\zeta)\in A_{2}^{1}(\Del)$.
For fixed $\zeta\in\Del$ and $f_1, f_2\in\mathcal{U}$ we define
\[
\psi(f_1,f_2; \zeta)(z)=K_{1}(f_1)(z,\zeta)-K_{1}(f_2)(z, \zeta).
\]
Then $\psi(f_{1},f_{2},\zeta)\in A_2^1(\Del)$. For $\psi(f_1,f_2;\zeta)\neq 0$ we set
\[
\mathrm{u}(f_1, f_2; \zeta)=\frac{\psi(f_1,f_2;\zeta)}{\Vert
\psi(f_1,f_2;\zeta)\Vert_2},\]
and for $\psi(f_{1},f_{2},\zeta)=0$ we set $\mathrm{u}(f_1, f_2; \zeta)=0$.
The following lemma  generalizes Brazilevic's result \cite{Bra}.
\begin{lemma} \label{bra} For $f_1, f_2\in\mathcal{U}$ and $z\in\Del$,
\begin{gather*}
(1-|z|^2)^2\left|\mathcal{S}(f_1)(z)-\mathcal{S}(f_2)(z)\right|\leq
6\FFF(f_1,f_2)(z)\leq6\Vert K_1(f_1)-K_1(f_2)\Vert.
\end{gather*}
\end{lemma}
\begin{proof}
We use the same approach as in \cite{Bra}. Since for $\lambda_{1},\lambda_{2}\in\PSL(2,\CC)$, $\mathcal{S}(\lambda_1\circ f)=\mathcal{S}(f)$,  $K_{1}(\lambda_{1}\circ f)=K_{1}(f)$, and $\FFF(\lambda_1\circ f_1, \lambda_2\circ f_2)(z)=
\FFF(f_1, f_2)(z)$, it is sufficient to consider only $f\in\mathcal{U}$ normalized by
$f(0)=0$ and $f'(0)=1$. We have for  fixed $z\in\Del$,
\begin{align*}
(1-|z|^2)^2\mathcal{S}(f)(z) &= \mathcal{S}(\lambda(f)_{z}\circ
f\circ \sigma_z)(0)\\
\intertext{and}
\FFF(f_1,f_2)(z)
&=\FFF(\lambda(f_{1})_{z}\circ f_1\circ \sigma_z, \lambda(f_2)_{z}\circ f_2\circ \sigma_z)(0),
\end{align*}
where $\sigma_{z}\in\Mob(S^{1})$ and $\lambda(f)_{z}\in\PSL(2,\CC)$ are given by\footnote{Here subscript $z$ does not denote a derivative.}
\[
\sigma_z(w)=\frac{w+z}{1+w\z}\quad \text{and}\quad\lambda(f)_{z}(w)=
\frac{w- f(z)}{f'(z)(1-|z|^{2})}.
\]
Since for a normalized $f\in\mathcal{U}$ the univalent function $\lambda(f)_{z}\circ f\circ \sigma_z$ is also normalized,
for the first inequality we need only to show that
for any normalized $f\in\mathcal{U}$,
\[
|\mathcal{S}(f_1)(0)-\mathcal{S}(f_2)(0)|\leq 6 \FFF(f_1,
f_2)(0).
\]
Since
\begin{align*}
\mathcal{S}(f)(z)= 6\lim_{w\rightarrow z}\left(
\frac{f'(z)f'(w)}{(f(z)-f(w))^2}-\frac{1}{(z-w)^2}\right)=-6\sum_{n,m=1}^{\infty}
nmb_{-n,-m}z^{n+m-2},
\end{align*}
we have
\begin{align*}
|\mathcal{S}(f_1)(0)-\mathcal{S}(f_2)(0)|=6|b_{-1,-1}(f_1)-b_{-1,-1}(f_2)|.
\end{align*}
On the other hand, it is straightforward to compute that
\begin{align*}
\FFF(f_1,f_2)^{2}(0) & = \pi
\iint\limits_{\Del}|K_1(f_1)(0,w)-K_1(f_2)(0,w)|^2d^2w
\\
&=\sum_{m=1}^{\infty}m|b_{-1,-m}(f_1)-b_{-1,-m}(f_2)|^2,
\end{align*}
and the first inequality follows. 

Next we observe that 
\begin{equation} \label{f=k}
\FFF(f_{1},f_{2})(z)=\sqrt{\pi}\Vert(K_{1}(f_{1})-K_{1}(f_{2}))\ov{\mathrm{u}(f_{1}, f_{2};z)}\Vert_{A^{1}_{\infty}(\Del)}.
\end{equation}
Indeed, by Cauchy-Schwarz inequality,
\begin{gather*}
\left(\left(K_1(f_1)-K_1(f_2)\right)\ov{\psi(f_1,f_2;z)}\right)(w) \\
=\iint\limits_{\Del}\left(K_1(f_1)(w,\zeta)-K_1(f_2)(w,\zeta)\right)
\ov{\left(K_1(f_1)(\zeta, z)-K_1(f_2)(\zeta,z)\right)}d^2\zeta\\
\leq \Vert\psi(f_{1},f_{2};z)\Vert_{2}\, \Vert\psi(f_{1},f_{2};w)\Vert_{2},
\end{gather*}
with the equality for $w=z$. Hence
\begin{gather*}
\left\Vert \left(K_1(f_1)-K_1(f_2)\right)\ov{\psi(f_1,
f_2;z)}\right\Vert_{A^{1}_{\infty}(\Del)} \\
=(1-|z|^2)\iint\limits_{\Del}\left|K_1(f_1)(\zeta,z)-K_1(f_2)(\zeta,z)\right|^2d^2\zeta\\
=(1-|z|^2)\Vert \psi(f_1,f_2;z) \Vert_2^2 =\FFF(f_{1},f_{2})(z)
\frac{\Vert \psi(f_1,f_2;z) \Vert_2}{\sqrt{\pi}}.
\end{gather*}

Finally, using \eqref{f=k} and the estimate in Lemma \ref{subspace1}, we get 
\begin{equation} \label{Estimate}
\FFF(f_{1},f_{2})(z)\leq \Vert(K_{1}(f_{1})-K_{1}(f_{2}))\ov{\mathrm{u}(f_{1},f_{2};z)}\Vert_{2}\leq \Vert K_{1}(f_{1})-K_{1}(f_{2})\Vert.
\end{equation}
\end{proof}
\begin{remark} It immediately follows from Lemma \ref{bra} that
\begin{align*}
\Vert \mathcal{S}(f_1)-\mathcal{S}(f_2)\Vert_{\infty}\leq 6 d(f_1,
f_2)\leq 6\Vert K_1(f_1)-K_1(f_2)\Vert,
\end{align*}
which is a stronger version of Brazilevic's result \cite{Bra}. In
case $f_1=f$ and $f_2=\mathrm{id}$ we have \[ \Vert
\mathcal{S}(f)\Vert_{\infty} \leq 6d(f) \leq 6\Vert
K_1(f)\Vert,
\] 
where $d(f)=d(f,\mathrm{id})$.
Since $ \Vert K_1(f)\Vert\leq 1$, where equality holds if and only if $\CC\setminus f(\Del)$ has Lebesgue measure zero, this
recovers another result in \cite{Bra} that $d(f)\leq 1$ for $f\in\mathcal{U}$, and  $d(f)=1$ implies that  $\CC\setminus f(\Del)$ has Lebesgue measure zero.
\end{remark}

Given a normalized univalent function $f:\Del\rightarrow \C$, let
$f_n:\Del\rightarrow \C$ be the normalized univalent function
defined by $f_n=r_{n}^{-1}\circ f\circ r_n$, where $r_{n}$ is the dilation $z\mapsto\tfrac{n}{n+1}z$.
Since $f_n$ is analytic on $S^1$, we have
\begin{align*}
\lim_{|z|\rightarrow 1^-} (1-|z|^2)^2\mathcal{S}(f_n)(z) &=0, \\
\lim_{|\zeta|\rightarrow 1^-} (1-|\zeta|^2)K_{1}(f_n)(z, \zeta) & =0,
\end{align*}
and also
\begin{align*}
\lim_{|\zeta|\rightarrow 1^-}\left\Vert (1-|\zeta|^2)K_{1}(f_{n})(\;\cdot\;,\zeta)
\right\Vert_2^2 &= \lim_{|\zeta|\rightarrow
1^-}(1-|\zeta|^2)^2\KKK_1(f_{n})(\zeta, \zeta)=0.
\end{align*}
\begin{lemma}\label{strong}
Let $f:\Del\rightarrow \C$ be a normalized univalent function and let
$\{f_n\}_{n=1}^{\infty}$ be the sequence of normalized univalent functions defined above. Then
\begin{align*}
\lim_{n\rightarrow \infty}K_1(f_n)=K_{1}(f)
\end{align*}
in the strong operator topology. 
\end{lemma}
\begin{proof}
For $\psi\in A_{2}^{1}(\Del)$ set $\psi_{n}=r_{n}\circ\psi\circ r_{n}$. It is elementary to show that
$$\lim_{n\rightarrow\infty}\Vert \psi-\psi_{n}\Vert_{2}=0.$$
For $(K_{1}(f)\bar{\psi})_{n} =r_{n}\circ (K_{1}(f)\bar{\psi})\circ r_{n}$ we have,
$$K_{1}(f_{n})\bar{\psi}_{n} =(K_{1}(f)\bar{\psi})_{n}-r_{n}\circ (K_{1}(f)\ov{\psi(1-\chi_{n})})\circ r_{n},$$
where $\chi_{n}$ is the characteristic function of the disk $\Del_{n}=r_{n}(\Del)$. Using this identity and the inequalities $\Vert K_{1}(f)\Vert\leq 1,\,\Vert\psi_{n}\Vert_{2}\leq \Vert\psi\Vert_{2}$, we obtain
\begin{gather*}\Vert (K_{1}(f)-K_{1}(f_{n}))\bar{\psi}\Vert_{2}\leq\Vert K_{1}(f) \bar{\psi} - K_{1}(f_{n})\bar{\psi}_{n}\Vert_{2}
+\Vert K_{1}(f_{n})\ov{(\psi_{n}-\psi)}\Vert_{2}\\
\leq \Vert K_{1}(f) \bar{\psi} - (K_{1}(f)\bar{\psi})_{n}\Vert_{2} +\Vert K_{1}(f)\ov{(\psi(1-\chi_{n})})\Vert_{2}
+\Vert \psi - \psi_{n} \Vert_{2}\\
\leq \Vert K_{1}(f) \bar{\psi} - (K_{1}(f)\bar{\psi})_{n}\Vert_{2} +\Vert \psi(1-\chi_{n})\Vert_{2}
+\Vert \psi - \psi_{n} \Vert_{2}.
\end{gather*}
Since $\psi\in A^{1}_{2}(\Del)$,
$$\lim_{n\rightarrow\infty}\Vert \psi(1-\chi_{n})\Vert_{2}=0,$$
and we get the assertion of the lemma.
\end{proof}
\begin{lemma}\label{con} Let $\gamma=g^{-1}\circ f\in\mathcal{T}(1)$ be such that $K_1(f)$ is a compact operator. Then for every sequence $\{\zeta_m\}_{m=1}^{\infty}$ of points in $\Del$, the corresponding sequence of functions $\{u_{m}\}_{m=1}^{\infty}$ in
$A^{1}_{2}(\Del)$, where
\begin{align*}
u_{m}(z)= (1-|\zeta_{m}|^2)K_{1}(f)(z,\zeta_{m}),\;z\in\Del,
\end{align*}
contains a convergent subsequence in $A_2^1(\Del)$.
\end{lemma}
\begin{proof}
Consider the following sequence of functions,
$$v_{m}(z) =z^{-2}(1-|\zeta_{m}|^{2})K_{3}(f)(z^{-1},\zeta_{m})\in A_{2}^{1}(\Del).$$
Using the formula
$$\KKK_{3}(\zeta,\zeta) +\KKK_{4}(\zeta,\zeta)=\frac{1}{\pi(1-|\zeta|^{2})^{2}},$$
which follows from the operator identity $\KKK_{3} +\KKK_{4}=I$, and the inequality $\KKK_{4}(\zeta,\zeta)\geq 0$, we get
$$\Vert v_{m}\Vert^{2}_{2}=(1-|\zeta_{m}|^{2})^{2}\KKK_{3}(\zeta_{m},\zeta_{m})\leq \frac{1}{\pi}.$$
Now consider
the operator $\tilde{K}_{3}(f): \ov{A_{2}^{1}(\Del)}\rightarrow A_{2}^{1}(\Del)$, defined by the kernel
$$\tilde{K}_{3}(f)(z,w)=z^{-2}K_{3}(f)\left(z^{-1},w\right).$$
In the standard basis for $A_{2}^{1}(\Del)$ it is given by the matrix $B_{3}(f)$ and, therefore, is a topological isomorphism. Setting
$K(f)=K_{1}(f)\tilde{K}_{3}(f)^{-1}$, we get
$$u_{m}=K(f)v_{m}.$$
Since the operator $K(f)$ is compact and the sequence $\{v_{m}\}_{m=1}^{\infty}$ is bounded, the statement follows.
\end{proof}
Now we can finish the proof of the Theorem. 
Suppose that for $[\mu]\in T(1)$ the corresponding operator $K_1(f)$ is compact but $[\mu]\notin S$. According to Remark \ref{ga-su}, this implies that there exist $\vep>0$ and a sequence $\zeta_m\in \Del$ satisfying
\begin{align*}
|\zeta_m|>1-\frac{1}{m}\quad\text{and}\quad
(1-|\zeta_m|^2)^2|\mathcal{S}(f)(\zeta_m)|\geq \vep.
\end{align*}
By Lemma \ref{con}, there exists a subsequence $\zeta_{m_k}$  such
that the sequence of functions
\begin{align*}
u_{m_{k}}(z)=(1-|\zeta_{m_k}|^2) K_{1}(f)(z,\zeta_{m_k})
\end{align*}
converges to $u\in A_2^1(\Del)$ in $A_{2}^{1}(\Del)$. Since
\begin{align*}
\lim_{|\zeta|\rightarrow 1^-}(1-|\zeta|^2)K_{1}(f_n)(z,\zeta)=0,
\end{align*}
for any $n\in \N$, the sequence of functions
\begin{align*}
(1-|\zeta_{m_k}|^2)\psi(f, f_n;
\zeta_{m_k})=(1-|\zeta_{m_k}|^2)\left(K_{1}(f)(\;\cdot\;,\zeta_{m_{k}})-K_{1}(f_n)(\;\cdot\;,
\zeta_{m_k})\right)
\end{align*}
also converges to $u$ as $k\rightarrow\infty$. 
From Lemma \ref{bra} and \eqref{Estimate} we get the following inequality
\begin{align*}
(1-|\zeta_{m_k}|^2)^2\left|\mathcal{S}(f)(\zeta_{m_k})-\mathcal{S}(f_n)(\zeta_{m_k})\right|
\leq  6\left\Vert \left(K_1(f)-K_1(f_n)\right) \ov{\mathrm{u}(f,
f_n;\zeta_{m_k})}\right\Vert_2,
\end{align*}
which for $\psi(f,f_{n},\zeta_{m_{k}})=0$ is an equality. Now
passing to the limit $k\rightarrow \infty$ for fixed $n\in \N$, we
obtain
\begin{align*}
\vep\leq 6\left\Vert
\left(K_1(f)-K_1(f_n)\right)\bar{\mathrm{u}}\right\Vert_2,
\end{align*}
where
$$\mathrm{u}=\frac{u}{\;\;\Vert u\Vert_{2}}\neq 0.$$
However, according to Lemma \ref{strong},
\[
\lim_{n\rightarrow
\infty}\left\Vert\left(K_1(f)-K_1(f_n)\right)\bar{\mathrm{u}}\right\Vert_2=0.
\]
This contradiction proves that $[\mu]\in S$.
\end{proof}
\begin{remark}
For $[\mu]\in S$ the proof of Lemma \ref{close} shows that
$$\lim_{n\rightarrow\infty}\mathcal{S}(f_n)=\mathcal{S}(f)$$
in $A_{\infty}(\Del)$ topology. Since the period mapping $\hat{\cP}$ is continuous, 
$$\lim_{n\rightarrow\infty}K_1(f_n)=K_1(f)$$
in the norm topology on $\cB(\ov{A_{2}^{1}(\Del)},A_{2}^{1}(\Del))$.
\end{remark}
The following commutative diagram displays the properties of the tower of embedded manifolds $T_0(1)\hookrightarrow
S\hookrightarrow T(1)$ under the KYNS period mapping $\hat{\cP}$,
the pre-Bers embedding $\hat{\beta}$ and the Bers embedding $\beta=\Psi\circ\hat{\beta}$,
\medskip
\medskip
\begin{displaymath}
\begin{CD} 
\cS_{2}(\ell^{2}) @>>>\cS_{\infty}(\ell^{2}) @>>> \cB(\ell^{2}) \\
@AA{\cP}A @AA{\hat{\cP}}A @AA{\hat{\cP}}A\\
T_{0}(1)@>>>S @>>> T(1) \\
@VV{\hat{\beta}}V @VV{\hat{\beta}}V @VV{\hat{\beta}}V \\
A^{1}_{2}(\Del) @>>>A^{1,0}_{\infty}(\Del) @>>> A^{1}_{\infty}(\Del) \\
@VV{\Psi}V @VV{\Psi}V @VV{\Psi}V \\
A_{2}(\Del) @>>>A^{0}_{\infty}(\Del) @>>> A_{\infty}(\Del)
\end{CD}
\end{displaymath}

\medskip
\medskip
\noindent
Here $A_{\infty}^{1,0}(\Del)$ is the closed subspace of $A^{1}_{\infty}(\Del)$, defined by
\begin{align*}
A_{\infty}^{1,0}(\Del)=\left\{ \psi\in
A_{\infty}^1(\Del)\; : \; \lim_{|z|\rightarrow
1^-}(1-|z|^2)\psi(z)=0\right\}.
\end{align*}
All horizontal maps  are embeddings, and all vertical maps are holomorphic mappings of Banach and Hilbert manifolds respectively. All these properties have been proved already, except for the simple fact $\Psi(A^{1,0}_{\infty}(\Del))\subset A^{0}_{\infty}(\Del)$, which easily follows from Cauchy integral formula. 
\bibliographystyle{amsalpha}
\bibliography{potential}

\def\cprime{$'$}
\providecommand{\bysame}{\leavevmode\hbox to3em{\hrulefill}\thinspace}
\providecommand{\MR}{\relax\ifhmode\unskip\space\fi MR }
\providecommand{\MRhref}[2]{%
  \href{http://www.ams.org/mathscinet-getitem?mr=#1}{#2}
}
\providecommand{\href}[2]{#2}
\begin{thebibliography}{Hum72}

\bibitem[BFK92]{burghelea}
D.~Burghelea, L.~Friedlander, and T.~Kappeler, \emph{Meyer-{V}ietoris type
  formula for determinants of elliptic differential operators}, J. Funct. Anal.
  \textbf{107} (1992), no.~1, 34--65.

\bibitem[BP78]{BecPom}
Jochen Becker and Christian Pommerenke, \emph{\"{U}ber die quasikonforme
  {F}ortsetzung schlichter {F}unktionen}, Math. Z. \textbf{161} (1978), no.~1,
  69--80.

\bibitem[Bra65]{Bra}
J.~E. Brazilevi{\v{c}}, \emph{On dispersion of coefficients of univalent
  functions}, Mat. Sb. (N.S.) \textbf{68 (110)} (1965), 549--560.

\bibitem[Dur83]{Duren}
P.~L. Duren, \emph{Univalent functions}, Grundlehren der Mathematischen
  Wissenschaften [Fundamental Principles of Mathematical Sciences], vol. 259,
  Springer-Verlag, New York, 1983.

\bibitem[GK69]{gohberg-krein}
I.~C. Gohberg and M.~G. Kre\u{\i}n, \emph{Introduction to the theory of linear
  nonselfadjoint operators}, Translated from the Russian by A. Feinstein.
  Translations of Mathematical Monographs, Vol. 18, American Mathematical
  Society, Providence, R.I., 1969.

\bibitem[GS92]{GaSu}
Frederick~P. Gardiner and Dennis~P. Sullivan, \emph{Symmetric structures on a
  closed curve}, Amer. J. Math. \textbf{114} (1992), no.~4, 683--736.

\bibitem[Ham02]{Hamilton}
D.~H. Hamilton, \emph{Conformal welding}, Handbook of complex analysis:
  geometric function theory, Vol.\ 1, North-Holland, Amsterdam, 2002,
  pp.~137--146.

\bibitem[Hum72]{Hummel}
J.~A. Hummel, \emph{Inequalities of {G}runsky type for {A}haronov pairs}, J.
  Analyse Math. \textbf{25} (1972), 217--257.

\bibitem[HZ99]{hassel-zelditch}
Andrew Hassell and Steve Zelditch, \emph{Determinants of {L}aplacians in
  exterior domains}, Internat. Math. Res. Notices (1999), no.~18, 971--1004.

\bibitem[KY88]{KY2}
A.~A. Kirillov and D.~V. Yuriev, \emph{Representations of the {V}irasoro
  algebra by the orbit method}, J. Geom. Phys. \textbf{5} (1988), no.~3,
  351--363.

\bibitem[Lan95]{Lang2}
Serge Lang, \emph{Differential and {R}iemannian manifolds}, third ed.,
  Springer-Verlag, New York, 1995.

\bibitem[Nag92]{Nag92}
Subhashis Nag, \emph{A period mapping in universal {T}eichm\"uller space},
  Bull. Amer. Math. Soc. (N.S.) \textbf{26} (1992), no.~2, 280--287.

\bibitem[NS95]{NS}
Subhashis Nag and Dennis Sullivan, \emph{Teichm\"uller theory and the universal
  period mapping via quantum calculus and the {$H\sp {1/2}$} space on the
  circle}, Osaka J. Math. \textbf{32} (1995), no.~1, 1--34.

\bibitem[Pom75]{Pom}
C.~Pommerenke, \emph{Univalent functions}, Vandenhoeck \& Ruprecht,
  G\"ottingen, 1975, With a chapter on quadratic differentials by Gerd Jensen,
  Studia Mathematica/Mathematische Lehrb\"ucher, Band XXV.

\bibitem[PS86]{Pressley-Segal}
Andrew Pressley and Graeme Segal, \emph{Loop groups}, Oxford Mathematical
  Monographs, The Clarendon Press Oxford University Press, New York, 1986.

\bibitem[Sch57]{Schiffer1}
M.~Schiffer, \emph{The {F}redholm eigen values of plane domains}, Pacific J.
  Math. \textbf{7} (1957), 1187--1225.

\bibitem[Sch59]{Schiffer-multiple}
\bysame, \emph{Fredholm eigen values of multiply-connected domains}, Pacific J.
  Math. \textbf{9} (1959), 211--269.

\bibitem[Sch81]{Schiffer2}
Menahem Schiffer, \emph{Fredholm eigenvalues and {G}runsky matrices}, Ann.
  Polon. Math. \textbf{39} (1981), 149--164.

\bibitem[Seg81]{Segal}
Graeme Segal, \emph{Unitary representations of some infinite-dimensional
  groups}, Comm. Math. Phys. \textbf{80} (1981), no.~3, 301--342.

\bibitem[SH62]{Schiffer3}
M.~Schiffer and N.~S. Hawley, \emph{Connections and conformal mapping}, Acta
  Math. \textbf{107} (1962), 175--274.

\bibitem[Sie64]{Siegel}
Carl~Ludwig Siegel, \emph{Symplectic geometry}, Academic Press, New York, 1964.

\bibitem[STZ99]{Todorov}
Maria~E. Schonbek, Andrey~N. Todorov, and Jorge~P. Zubelli, \emph{Geodesic
  flows on diffeomorphisms of the circle, {G}rassmannians, and the geometry of
  the periodic {K}d{V} equation}, Adv. Theor. Math. Phys. \textbf{3} (1999),
  no.~4, 1027--1092.

\bibitem[SW85]{Segal-Wilson}
Graeme Segal and George Wilson, \emph{Loop groups and equations of {K}d{V}
  type}, Inst. Hautes \'Etudes Sci. Publ. Math. (1985), no.~61, 5--65.

\bibitem[Teo02]{Teo}
Lee-Peng Teo, \emph{Velling-{K}irillov metric on the universal
  {T}eichm\"{u}ller curve,}, Preprint arXiv: math.CV/0206202 (2002).

\bibitem[Teo03]{Teo2}
\bysame, \emph{Analytic functions and integrable hierarchies---characterization
  of tau functions}, Lett. Math. Phys. \textbf{64} (2003), no.~1, 75--92.

\bibitem[TT03a]{LT}
Leon~A. Takhtajan and Lee-Peng Teo, \emph{Liouville action and
  {W}eil-{P}etersson metric on deformation spaces, global {K}leinian
  reciprocity and holography}, Comm. Math. Phys. \textbf{239} (2003), no.~1-2,
  183--240.

\bibitem[TT03b]{TT-I}
\bysame, \emph{{W}eil-{P}etersson metric on the universal {T}eichm\"{u}ller
  space {I}: {C}urvature properties and {C}hern forms}, Preprint arXiv:
  math.CV/0312172 (2003).

\bibitem[ZT87]{ZT}
P.~G. Zograf and L.~A. Takhtadzhyan, \emph{A local index theorem for families
  of $\overline\partial$-operators on {R}iemann surfaces}, Uspekhi Mat. Nauk
  \textbf{42} (1987), no.~6(258), 133--150 (Russian); English. transl. in:
  Russian Math. Surveys, \textbf{42:6} (1987), 169--190.

\end{thebibliography}
\end{document}